\documentclass[reqno]{amsart}

\usepackage{lineno}
\usepackage[toc]{appendix}
\modulolinenumbers[5]
\usepackage{comment}
\usepackage[numbers,sort&compress]{natbib}
\usepackage{graphicx}
\usepackage{epsfig,epstopdf}
\usepackage{amsmath,amssymb}
\usepackage{mathrsfs}
\usepackage{paralist}
\usepackage{pdflscape}
\usepackage{tikz}
\usetikzlibrary{arrows.meta}
\usepackage{nicematrix}

\usepackage{graphicx}
\usepackage{amsmath,amssymb}
\usepackage{amsthm}
\usepackage[table,dvipsnames]{xcolor}
\usepackage{mathrsfs}
\usepackage{paralist}
\usepackage{multirow}
\usepackage[T1]{fontenc}
\usepackage[utf8]{inputenc}
\usepackage[english]{babel}
\DeclareMathAlphabet\mathbfcal{OMS}{cmsy}{b}{n}
\usepackage[figurename=Figure]{caption}
\usepackage{lineno}
\usepackage{bbm}
\usepackage[normalem]{ulem}
\usepackage[percent]{overpic}
\usepackage{color}
\usepackage{tikz}
\usepackage{yhmath}

\usepackage{wrapfig}
\usepackage{float}
\usepackage{booktabs}
\usepackage{makecell}
\usepackage{subcaption}
\usepackage{comment}

\definecolor{stablegreen}{RGB}{210,240,210}
\definecolor{unstablered}{RGB}{245,210,210}

\usepackage{hyperref}
\hypersetup{
    colorlinks=true,     
    linkcolor=red,       
    citecolor=green,     
    filecolor=magenta,   
    urlcolor=blue        
}

\newcommand{\ud}{\mathrm{d}}

\newcommand{\ie}{\emph{i.e.}} 
\newcommand{\eg}{\emph{e.g.}} 
\newcommand{\cf}{\emph{cf.}} 
\renewcommand{\Pr}{\mathrm{Pr}} 
\newcommand{\Rey}{\mathrm{Re}} 
\newcommand{\Ma}{\mathrm{Ma}} 
\newcommand{\mr}[1]{\mathrm{#1}}

\newcommand{\E}[1]{\times 10^{#1}}
\newcommand{\vect}[1]{\boldsymbol{\mathrm{#1}}}
\newcommand{\uvect}[1]{\hat{\vect{#1}}} 
\newcommand{\tens}[1]{\boldsymbol{\mathsf{#1}}}
\newcommand{\tensij}[1]{\mathsf{#1}}
\newcommand{\norm}[1]{\left\lVert #1 \right\rVert}
\newcommand{\abs}[1]{\left\lvert #1 \right\rvert}

\newcommand{\pd}[2]{\frac{\partial{#1}}{\partial{#2}}} 
\newcommand{\ilpd}[2]{\partial{#1}/\partial{#2}} 
\newcommand{\der}[2]{\frac{\ud{#1}}{\ud{#2}}} 
\newcommand{\ilder}[2]{\ud{#1}/\ud{#2}} 

\newcommand{\grad}[1]{\nabla #1}
\newcommand{\divg}[1]{\nabla\cdot{#1}}
\newcommand{\curl}[1]{\nabla\times{#1}}
\newcommand{\lapl}[1]{\nabla^2 #1}              

\newcommand{\del}[1]{}

\usepackage{authblk}

\usepackage{amsaddr}
\usepackage{siunitx}
        
\usepackage[a4paper,bindingoffset=0.0in,%
            left=1.0in,right=1.0in,top=1.in,bottom=1.in,%
            footskip=.4in]{geometry}

\title[A novel compact scheme for second-order fluxes applied to the SD method]{A novel compact scheme for second-order fluxes applied to the Spectral Difference method}

\author{Guido Lodato$^{1}$ and Niccol{\`o} Tonicello$^{2}$}
\address{$^1$ Normandie Université, INSA et Université de Rouen, St. Etienne du Rouvray (Rouen), France}
\address{$^2$ MathLab, Mathematics Area, SISSA, via Bonomea 265, I-34136 Trieste, Italy}

\begin{document}
\begin{abstract}
The discretization of second-order (viscous) terms within Discontinuous Spectral Element Methods (DSEMs) typically relies on the introduction of an auxiliary gradient variable, whose numerical treatment at element interfaces strongly affects both the accuracy and the stability of the resulting scheme. Among the possible choices, the Bassi-Rebay (BR1)~\citep{bassi:97} formulation is particularly attractive owing to its simplicity, parameter-free character and applicability to nonlinear fluxes and arbitrary grids. However, it is well known to suffer from sub-optimal convergence at even orders of approximation and to require an extended, five-element computational stencil in one dimension. In this work, inspired by the formulation proposed by Huynh for the Flux Reconstruction approach~\citep{huynh:09}, we develop a novel compact, fully-centered scheme for the discretization of second-order fluxes within the Spectral Difference (SD) method. The proposed approach modifies the reconstruction of the auxiliary gradient variable through the introduction of interface-dependent, one-sided continuous fluxes, which reduces the stencil required for the viscous discretization from five to three elements while preserving the centered and parameter-free nature of the original BR1 formulation. The methodology is presented for the one-dimensional case and subsequently extended to multiple dimensions. A temporal eigenanalysis, including a combined-mode formulation and the addition of interior penalty terms, is used to characterize the dissipation and dispersion properties of the new scheme in comparison with the standard extended-stencil formulation. A series of numerical experiments is then considered: a convergence study for the linear diffusion equation, the diffusion of an under-resolved localized Dirac's delta, the nonlinear one-dimensional porous medium equation, and implicit large-eddy simulations of the three-dimensional Taylor-Green vortex at $\Rey = 1\,600$ and $5\,000$. The results show that the compact formulation restores the expected convergence order for all polynomial degrees, including the even orders for which the standard BR1 scheme underperforms, while substantially suppressing spurious oscillations in under-resolved and nonlinear settings. In the three-dimensional turbulent test cases, the compact scheme is shown to remain stable in regimes where the standard formulation fails, owing to its improved damping of high-wavenumber, spurious numerical modes. These results indicate that the proposed compact formulation represents an attractive alternative to the standard BR1 approach for the discretization of second-order fluxes within the SD method, particularly for large-scale and implicit LES applications.
\end{abstract}

\maketitle

\section{Introduction}
The significant growth in computational resources over the last few decades has made Computational Fluid Dynamics (CFD) a widespread tool for many engineering design processes. Among the many scientific advancements, the development of innovative numerical schemes has experienced significant growth, seeking candidates for the next generation of CFD solvers. The broad family of Discontinuous Spectral Element Methods (DSEMs) has shown significant promise due to their favorable numerical dispersion/diffusion characteristics, geometrical flexibility, and algorithmic scalability~\citep{vincent2011facilitating}, even more so with the widespread usage of GPUs~\citep{fehn2018efficiency,gasparino2024sod2d,kurz2025galaexi}. These include the Discontinuous Galerkin (DG)~\citep{hesthaven2007nodal,cockburn:98,cockburn:98b}, the Flux Reconstruction (FR)~\citep{huynh2007flux,vincent2011new}, and the Spectral Difference (SD)~\citep{kopriva:96,liu:06a,wang:07,jameson:10} schemes.

High-order discontinuous spectral element methods have been applied to a large variety of Partial Differential Equations (PDEs) in fluid dynamics such as incompressible and compressible Navier-Stokes equations~\citep{bosnyakov2014high,mengaldo2021industry,tonicello2022analysis,moxey2020nektar++,lv2021discontinuous,dzanic2022positivity,ferrer2023high}, with extension to multi-component~\citep{marchal2023extension,lv2014discontinuous,ching2025positivity,rising2026simulation} and multi-phase~\citep{manzanero2020entropy,ntoukas2022entropy,tonicello2024high,orlando2024implicit,tonicello2025extension} formulations. In particular, such methodologies have been deeply studied for scale-resolving simulations of turbulence such as Direct Numerical Simulations (DNS)~\citep{chapelier2014evaluation,de2018use,bassi2016development,krank2018direct,tonicello2022turbulence,lodato:16,lodato:17} and Large-Eddy Simulations (LES)~\citep{lodato2021large,tonicello2021analysis,fernandez2018ability,moura2017eddy,de2019performance,mengaldo2021industry}. In this scenario, the high-order capabilities of these schemes represent a significant advantage over low-order methods, generally providing higher accuracy with a significantly smaller total number of degrees of freedom (DoF).

Within the specific context of LES, a large part of the community has investigated the use of so-called \emph{implicit} LES (ILES). Such approaches aim at designing the numerical discretization (in particular, its numerical dissipation) so that it acts as an implicit, built-in sub-grid scale model for turbulence. Consequently, a large amount of research has focused on simplified techniques to quantify and investigate the numerical dissipation and dispersion properties of DSEMs. Most of these techniques are based on linear eigenanalysis and its many variations, such as temporal~\citep{lele1992compact,bogey2004family,van2008stability,vincent2011insights,moura2015linear,vanharen2017revisiting}, spatial~\citep{hu2002eigensolution,mengaldo2018spatial,mengaldo2018spatial2,moura2020spatial,tonicello2021comparative}, non-modal analysis~\citep{Fernandez_2019}, and combined-mode analysis~\citep{alhawwary2020combined}. These techniques are generally based on the numerical discretization of the linear advection equation, considering a wave-like ansatz that makes it possible to analytically study how each wavenumber evolves in time and space in terms of propagation speed and damping.

As most relevant turbulence applications are, almost by definition, characterized by advection-dominated flows, the vast majority of these analyses have focused on the linear advection equation as a representative framework to study the evolution of spatio-temporal frequencies. Within the DSEM setting, the choice of polynomial reconstruction and numerical fluxes are, in particular, the most relevant tuning parameters for implicit LES. The latter, for the Navier-Stokes equations, can be computed using Riemann solvers for general hyperbolic systems of equations. Many different Riemann solvers have been proposed in the literature, and all of them, in the regime of under-resolved flows such as LES, behave differently depending on how they handle discontinuities between elements. As a simple example, it is well known that the classical Rusanov flux can lead to hyper-upwinding in low-Mach flows, whereas the Roe flux is capable of modulating the amount of numerical dissipation across different flow regimes and variables. Consequently, ILES performed using Rusanov or Roe fluxes behaves significantly differently in this regime (see, for example,~\cite{moura2017eddy}). Similarly, Large-Eddy Simulations performed with different polynomial orders but the same total number of degrees of freedom also behave significantly differently (see, for example,~\cite{chapelier2016spectral}).

Only recently attention has been devoted to the discretization of second-order terms within spectral element methods, with applications to LES~\cite{ferrer2017interior,kou2023jump,du2026assessment}. In fact, even in advection-dominated flows, the representation of viscous terms can be particularly important. This is certainly the case, for example, in the vicinity of wall boundaries or when using explicit eddy-viscosity models in significantly under-resolved LES. Consequently, in advection-diffusion problems, the numerical scheme should properly describe both the hyperbolic and parabolic features of the governing equations. While DSEMs are particularly well suited for hyperbolic problems, they require additional corrections to handle parabolic systems.

Within the framework of DSEMs, the main idea is to rewrite the second-order PDEs into an extended system of first-order PDEs by introducing the gradient as an additional unknown. The only difference in discretizing such a system using DSEMs is the appearance of additional interfacial contributions for the second-order terms, requiring an appropriate definition of numerical fluxes for both the solution and its gradients at the element interfaces.

The simplest approach, proposed by Bassi and Rebay (BR1)~\cite{bassi:97} for the DG formulation, consists of simply taking the arithmetic mean of both the solution and its gradients at the interface. This formulation is particularly convenient since it is simple to implement, parameter-free, and applicable to nonlinear fluxes and arbitrary grids. However, it is also affected by well-known deficiencies, such as sub-optimal convergence for even orders of approximation, a consistent but unstable character, and an extended computational stencil involving five elements in the one-dimensional case. Since the first work on the SD scheme by Kopriva~\cite{kopriva:98} employed an adapted version of the BR1 scheme, we will consider it as a relevant comparison throughout the remainder of the paper.

Since the introduction of the BR1 approach, several additional formulations have been proposed to overcome its deficiencies. Among these, the BR2~\cite{bassi1997high2}, Local Discontinuous Galerkin (LDG)~\cite{cockburn1998local}, and Interior Penalty (IP)~\cite{arnold1982interior} approaches are certainly worth mentioning. Many of these, although better performing than BR1 for non negligible viscous effects, are characterized by another set of drawbacks such as stiffer CFL restrictions, ad-hoc parameter dependencies or persistent suboptimal convergence~\cite{hartmann2008optimal}.

It is also relevant to note that the large majority of such methodologies have been applied primarily within a DG/FR setting, while very limited literature is currently available regarding the treatment of second-order terms within the SD scheme~\citep{marchal2023extension}.

In this work, inspired by the formulation proposed by Huynh~\cite{huynh:09} for Flux-Reconstruction schemes, we propose a compact, fully-centered formulation for the discretization of second-order operators for the SD scheme. The main idea is to modify the reconstruction of the auxiliary gradient variable by introducing interface-dependent continuous fluxes, which allows the information exchange between neighboring elements to be localized. Instead of first constructing a globally continuous auxiliary flux and then computing the gradient, the proposed formulation introduces left- and right-sided continuous fluxes within each element, which are corrected only at the corresponding interface. The interfacial gradient is then obtained by averaging the one-sided contributions from the two neighboring elements. This procedure preserves the centered nature of the BR1 formulation while reducing the stencil required for the computation of the auxiliary variable from five to three elements. As a result, the stencil required by the viscous discretization is reduced from five to three elements while preserving the centered and parameter-free character of the original BR1 formulation. The proposed approach retains the simplicity and generality of BR1, while providing a more compact discretization that is potentially better suited for large-scale simulations and implicit LES applications.

The paper is organized as follows. In section~\ref{sec:1}, we introduce the methodology for the one-dimensional pure diffusion equation within the SD scheme for both classical BR1 numerical flux  and for the newly proposed compact formulation. In the same section, we present the extension to multiple dimensions. In section~\ref{sec:2}, we use temporal eigenanalysis to inspect the numerical dissipation and dispersion properties of both approaches implemented within the SD scheme, including the possibility to augment both methodologies with interior penalty terms. In section~\ref{sec:3} we consider a series of numerical experiments of increasing complexity. We start by performing a convergence study for the one-dimensional diffusion equation. We then inspect the robustness of the proposed approach for an under-resolved case of pure diffusion of a localized Dirac's delta. We then move to non-linear diffusion problems by considering the one-dimensional porous medium equation. As last numerical experiment, we tested the compact formulation for the three-dimensional Taylor-Green Vortex problem. Finally, in section~\ref{sec:4} we outline the key conclusions of this work.

\section{Methodology for pure diffusion equation}\label{sec:1}
In order to introduce the notation which will be instrumental to perform the eigenanalysis in subsequent sections, the original SD method by~\citet{kopriva:96,kopriva:98} for second-order fluxes is here briefly outlined.
Let us consider the following three-dimensional conservation law,
\begin{equation}
\pd{u}{t} + \divg{ \vect{F}(u,\grad{u}) } = 0,
\end{equation} 
where $\vect{F}$ is a flux vector which depends on both the solution $u$ and its gradient $\nabla u$.
The above equation, which is second-order in space, is reduced to the following system of first-order equations via the introduction of the auxiliary variable $\vect{v}$:
\begin{equation}
\left\{
\begin{aligned}
\pd{u}{t} + \divg{ \vect{F}(u,\vect{v}) } &= 0,\\
\grad{u} &= \vect{v}.
\end{aligned}
\right.\label{eq:ns}
\end{equation}

\begin{figure}
\centering
\includegraphics*[trim=0 0 0 0,width=0.98\textwidth]{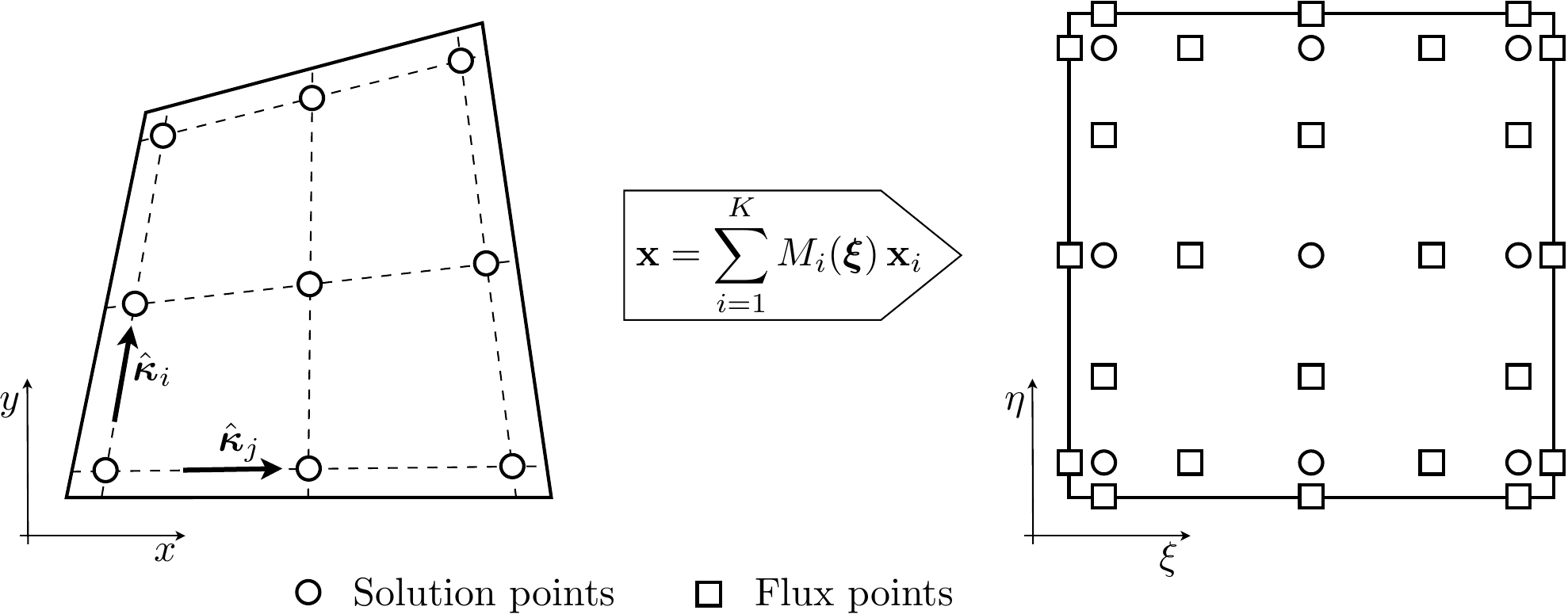}
\caption{Schematic representation of the two-dimensional distribution of solution and flux points within the SD element for $n=3$.}
\label{fig:2d:sketch}
\end{figure}

Equation~\eqref{eq:ns} is integrated over a physical domain which is subdivided into $N_{\rm e}$ non-overlapping hexahedral elements.
To achieve an efficient implementation, each element in the physical domain is transformed to a standard cubic element described by local coordinates  $\vect{\xi} = (\xi,\;\eta,\;\zeta) \in [-1:1]^3$ (\cf~figure~\ref{fig:2d:sketch}) via the introduction of the relevant transformation of coordinates
\begin{equation}
\vect{x}_e = \sum_{i=1}^K M_i(\vect{\xi}) \,\vect{x}_{e,i},
\end{equation}
where $M_i(\vect{\xi})$ are suitable shape functions and $K$ is the number of points $\vect{x}_{e,i}$ defining the $e$-th physical element (\eg, $K=8$ or $20$ for standard linear or quadratic hexahedra, respectively). 
After introducing the Jacobian of the transformation, $\tens{J}_e = \ilpd{\vect{x}_e}{\vect{\xi}}$, its determinant, $J_e = \det(\tens{J}_e)$, and its adjoint, $\tens{S}_e = \mr{adj}(\tens{J}_e)$,\footnote{Note that the $i$-th row of the adjoint evaluated at the element's interfaces orthogonal to $\xi_i$, namely $\ilpd{\vect{x}}{\xi_j} \times \ilpd{\vect{x}}{\xi_k}$ for $\xi_i = \pm 1$, relates to the elementary area vector of that interface. This motivates the use of the symbol $\tens{S}$ for the adjoint.} the governing equations in the computational domain take the form\footnote{\label{note:metric:identities}In order to obtain equations~\eqref{eq:ns:mapped} and~\eqref{eq:adj:f}, the \emph{metric identities}, namely $\nabla_\xi \cdot \mr{adj}(\tens{J}) = \vect{0}$, have been used.}
\begin{equation}
\left\{
\begin{aligned}
\pd{q_e}{t} + \nabla_\xi \cdot \vect{G}_e(u_e,\vect{v}_e)  &= 0,\\
\nabla_\xi \cdot \tens{H}_e(u_e) &= \vect{w}_e,
\end{aligned}
\right.
\label{eq:ns:mapped}
\end{equation}
with
\begin{equation}
\begin{aligned}
q_e = J_e u_e, &\quad 
\vect{G}_e(u_e,\vect{v}_e) = \tens{S}_e \cdot \vect{F}(u_e,\vect{v}_e),\\
\vect{w}_e = J_e \vect{v}_e, &\quad
\tens{H}_e(u_e) = \tens{S}_e^{\top} u_e. 
\end{aligned}\label{eq:adj:f}
\end{equation}

In the above relations, $u_e$ and $\vect{v}_e$ are polynomial approximations of the solution and the auxiliary variable, respectively. The relevant order depends on the number of supporting points as explained below.
Within each standard element, two sets of points are defined, namely, $n$ solution points and $n+1$ staggered flux points in each direction, as schematically illustrated in figure~\ref{fig:2d:sketch} for the two-dimensional case.
The solution points are used to support the polynomial approximation of the solutions $q_e$ and $\vect{w}_e$, whereas the flux points are used to support the higher-order polynomial approximation of the fluxes $\vect{G}_e$ and $\tens{H}_e$.
Along each direction, the solution points are selected according to the $n$-points Gauss-Legendre quadrature rule, whereas, the flux points are located at the Gauss-Legendre quadrature points of order $n-1$ plus the two end points at the element's interfaces~\citep{jameson:10,jameson:12}.

\subsection{The one-dimensional case}\label{sec:1d:std}

In the one-dimensional case, and for linear elements, the Jacobian of the transformation to the computational space reduces to a scaling factor between the physical element and the standard one,
\begin{equation}
x = \underbrace{\tfrac{1}{2}(1-\xi)}_{M_1(\xi)} x_{e} + \underbrace{\tfrac{1}{2}(1+\xi)}_{M_2(\xi)} x_{e+1},\quad\Rightarrow J_e = \tfrac{1}{2}(x_{e+1} - x_{e})\quad\text{and}\quad
\tensij{S}_e = 1,
\end{equation}

and equations~\eqref{eq:ns:mapped} and~\eqref{eq:adj:f} within the $e$-th element reduce to
\begin{equation}
\left\{
\begin{aligned}
\pd{q_e}{t} + \pd{}{\xi} G_e(u_e,v_e)  &= 0,\\
\pd{}{\xi} H_e(u_e) &= w_e,
\end{aligned}
\right.
\label{eq:ns:mapped:1d}
\end{equation}
with
\begin{equation}
q_e = J_e u_e, \quad 
G_e(u_e,v_e) = F(u_e,v_e), \quad
w_e = J_e v_e, \quad
H_e(u_e) = u_e. 
\label{eq:adj:f:1d}
\end{equation}

The main ingredients to build the SD scheme are:
\begin{inparaenum}[(a)]
\item the interpolation operator from solution points to flux points;
\item the differentiation operator from flux points to flux points.
\end{inparaenum}
Both operators are built from the polynomial representation of the solution or its flux using the Lagrange polynomial bases on the $n$ solution points or the $n+1$ flux points, respectively:
\begin{equation}
\ell_i(\xi) = \prod_{\substack{k=1\\k\neq i}}^n \frac{\xi-\xi_k}{\xi_i-\xi_k}, \quad\text{and}\quad
\ell_{i-1/2}(\xi) = \prod_{\substack{k=1\\k\neq i}}^{n+1} \frac{\xi-\xi_{k-1/2}}{\xi_{i-1/2}-\xi_{k-1/2}},
\label{eq:lagr:poly}
\end{equation}
where $\xi_i$, with $i=1,\dots,n$, and $\xi_{i-1/2}$, with $i=1,\dots,n+1$, are the locations of the solution and flux points, respectively (note that $\xi_{1/2} = -1$ and $\xi_{n+1/2} = 1$).
The above bases provide an order $p = n - 1$ polynomial approximation of the solution and an order $p + 1 = n$ polynomial approximation of the flux.

Starting from the polynomial approximation of the solution, this is written as
\begin{equation}
u_e(\xi,t) = \sum_{i=1}^n u_e(\xi_i,t) \ell_i(\xi),
\label{eq:sol:poly}
\end{equation}
where $u_e(\xi_i,t)$ are the nodal values of the solution on the solution points, which can be represented by the solution vector $\vect{u}_e^{\rm S}(t) = (u_e(\xi_1,t),\dots,u_e(\xi_n,t))$. Similar vectors can be defined for every nodal quantity. The time dependance of nodal values will be omitted for ease of notation hereafter.

Using equation~\eqref{eq:sol:poly}, the interpolation operator is easily obtained:
\begin{equation}
\mathcal{L}_{ij} = \ell_j(\xi_{i-1/2}) \in \mathbb{R}^{(n+1)\times n},\quad\Rightarrow
\vect{u}_e^{\rm F} = \tens{\mathcal{L}} \cdot \vect{u}_e^{\rm S},
\label{eq:Lmat}
\end{equation}
$\vect{u}_e^{\rm F} \in \mathbb{R}^{n+1}$ being the solution vector at the flux points.
Similarly, if $\vect{\phi}_e^{\rm F} = (\phi_e(\xi_{1/2}),\dots,\phi_e(\xi_{n+1/2}))$ is a flux vector within element $e$, the differentiation operator becomes:
\begin{equation}
\mathcal{M}_{ij} = \der{}{\xi}\ell_{j-1/2}(\xi_{i}) \in \mathbb{R}^{n \times (n+1)},\quad\Rightarrow
\der{\vect{\phi}_e^{\rm S}}{\xi} = \tens{\mathcal{M}} \cdot \vect{\phi}_e^{\rm F},
\label{eq:Mmat}
\end{equation}
where $\ud \vect{\phi}_e^{\rm S} / \ud \xi  \in \mathbb{R}^{n}$ is the derivative of the flux vector at the solution points.\footnote{
In the present case of staggered flux and solution points, the derivative of the Lagrange polynomials can be computed using the following identity,
$$\der{}{\xi}\prod_k a_k(\xi) = \sum_j\bigg[ \der{a_j}{\xi} \frac{1}{a_j} \prod_k a_k \bigg],\quad\Rightarrow \der{}{\xi}\ell_{i-1/2}(\xi) = \sum_{\substack{j=1\\j\neq i}}^{n+1} \bigg[ \frac{1}{\xi-\xi_{j-1/2}} \prod_{\substack{k=1\\k\neq i}}^{n+1} \frac{\xi-\xi_{k-1/2}}{\xi_{i-1/2}-\xi_{k-1/2}}\bigg] = \ell_{i-1/2}(\xi) \sum_{\substack{j=1\\j\neq i}}^{n+1} \bigg[ \frac{1}{\xi-\xi_{j-1/2}} \bigg].$$
}

For later use, correction vectors $\in \mathbb{R}^{n+1}$ are also defined as:
\begin{align}
\begin{aligned}
\uvect{c}_{\rm L} &= \{ \ell_{1/2}(\xi_{i-1/2}) \}_{i = 1,\dots,n+1} = (1,\;0,\dots,\;0)^\top,\\
\uvect{c}_{\rm R} &= \{ \ell_{n+1/2}(\xi_{i-1/2}) \}_{i = 1,\dots,n+1} = (0,\dots,\;0,\;1)^\top,
\end{aligned}
\label{eq:cl:cr}
\end{align}
which will be used to impose corrections at interface flux points.

Back to the equation~\eqref{eq:ns:mapped:1d}, the standard SD scheme, which corresponds to the BR1 formulation for the DG method~\citep{bassi:97,bassi:00}, is built as follows~\citep{kopriva:96,kopriva:98,sun:07}:
\begin{enumerate}
\item The solution $u_e$ is first interpolated to the flux points in each element:
$$\vect{u}_e^{\rm F} = \tens{\mathcal{L}} \cdot \vect{u}_e^{\rm S}.$$

Note that the above solution is generally discontinuous across elements' interfaces.

\item\label{item:2} The flux $H_e$ is readily computed at the interior flux points within each element.
On the other hand, at each element' interface, common numerical fluxes $\hat{H}$ must be defined. A common choice is the centered flux:
$$\hat{H}(u_{\rm L}, u_{\rm R}) = \tfrac{1}{2}(u_{\rm L} + u_{\rm R}).$$
Accordingly, the continuous flux vector at flux points can be assembled from the discontinuous one by including corrections at the interfaces:
\begin{equation}
\vect{H}_e^{\rm F,C} = \vect{u}_e^{\rm F,C} = \vect{u}_e^{\rm F} + [\hat{H}_{e}^{\rm L} - u_e(-1)]\uvect{c}_{\rm L} + [\hat{H}_{e}^{\rm R} - u_e(1)]\uvect{c}_{\rm R},
\label{eq:c:flux}
\end{equation} 
where $u_e(-1) = \uvect{c}_{\rm L} \cdot \vect{u}_e^{\rm F}$ and $u_e(1) = \uvect{c}_{\rm R} \cdot \vect{u}_e^{\rm F}$ are the first and last components of $\vect{u}_e^{\rm F}$,  and
\begin{equation}
\begin{aligned}
\hat{H}_{e}^{\rm L} &= \hat{H}(u_{e-1}(1), u_{e}(-1)) = \tfrac{1}{2}(\uvect{c}_{\rm R}^\top \cdot \vect{u}_{e-1}^{\rm F} + \uvect{c}_{\rm L}^\top \cdot \vect{u}_e^{\rm F}),\\
\hat{H}_{e}^{\rm R} &= \hat{H}(u_{e}(1), u_{e+1}(-1)) = \tfrac{1}{2}(\uvect{c}_{\rm R}^\top \cdot \vect{u}_{e}^{\rm F} + \uvect{c}_{\rm L}^\top \cdot \vect{u}_{e+1}^{\rm F}).
\end{aligned}
\end{equation}

It is worthwhile noting that equation~\eqref{eq:c:flux} involves a stencil that spans over three elements. 

\item\label{item:3} The auxiliary variable $\vect{v}_e^{\rm S}$ at solution points is readily obtained via the differentiation operator,
\begin{equation}
\vect{v}_e^{\rm S} = J_e^{-1}\tens{\mathcal{M}} \cdot \vect{H}_e^{\rm F,C},
\label{eq:w:e}
\end{equation}
and interpolated at the flux points,
\begin{equation}
\vect{v}_e^{\rm F} = \tens{\mathcal{L}} \cdot \vect{v}_e^{\rm S},
\label{eq:v:e}
\end{equation}
which is again a discontinuous quantity across elements' interfaces.

\item Concerning the computation of the flux $G_e$, again, its evaluation at interior flux points is straightforward using interior values of $\vect{u}_e^{\rm F}$ and $\vect{v}_e^{\rm F}$.
At the element interfaces, the already computed common flux from equation~\eqref{eq:c:flux} can be used, but some care must be taken for the auxiliary variable $v$ which is double-valued.
Here, two typical possibilities are:
\begin{inparaenum}[(a)]
\item use the average interface value of $\hat{v} = (v_{\rm L} + v_{\rm R})/2$ before computing $\hat{G} = G_e(\hat{u},\hat{v})$, with $\hat{u} = \hat{H}$ from step~\ref{item:2} above~\citep{sun:07};
\item compute the average interface value of $G_e$ from its left/right values obtained from discontinuous solutions, namely, $\hat{G} = (G_{e,{\rm L}} + G_{e,{\rm R}})/2$, with $G_{e,{\rm L/R}} = G_e(u_{\rm L/R},v_{\rm L/R})$.
\end{inparaenum}

Both approaches can be augmented by an interior penalty (IP) term~\citep{hesthaven:book,shahbazi:05,manzanero:18,arnold:02}, the latter having been the choice to couple IP terms in recent works with the SD method~\citep{lodato:19,lodato:19corr,lodato:22b}.

Focusing on the former approach, the continuous auxiliary variable is defined as
\begin{equation}
\vect{v}_e^{\rm F,C} = \vect{v}_e^{\rm F} + [\hat{v}_{e}^{\rm L} - v_e(-1)]\uvect{c}_{\rm L} + [\hat{v}_{e}^{\rm R} - v_e(1)]\uvect{c}_{\rm R},\quad\Rightarrow \vect{G}_e^{\rm F,C} = G_e(\vect{u}_e^{\rm F,C},\vect{v}_e^{\rm F,C}),
\label{eq:c:aux:v}
\end{equation}
where $v_e(-1) = \uvect{c}_{\rm L}^\top \cdot \vect{v}_e^{\rm F}$ and $v_e(1) = \uvect{c}_{\rm R}^\top \cdot \vect{v}_e^{\rm F}$ and
\begin{equation}
\begin{aligned}
\hat{v}_{e}^{\rm L} &= \tfrac{1}{2}(\uvect{c}_{\rm R}^\top \cdot \vect{v}_{e-1}^{\rm F} + \uvect{c}_{\rm L}^\top \cdot \vect{v}_e^{\rm F}),\\
\hat{v}_{e}^{\rm R} &= \tfrac{1}{2}(\uvect{c}_{\rm R}^\top \cdot \vect{v}_{e}^{\rm F} + \uvect{c}_{\rm L}^\top \cdot \vect{v}_{e+1}^{\rm F}).
\label{eq:aux:v:int}
\end{aligned}
\end{equation}

A very important point is that the stencil for $\vect{v}_e^{\rm F,C}$ now spans over five elements from $e-2$ to $e+2$. This is due to the fact that the stencil to obtain $\vect{v}_{e-1}^{\rm F}$ from $\vect{H}_{e-1}^{\rm F,C}$ involves elements $e-2$, $e-1$ and $e$ and, analogously, the stencil to obtain $\vect{v}_{e+1}^{\rm F}$  involves elements $e$, $e+1$ and $e+2$.

\item The final steps involve the differentiation of  $\vect{G}_e^{\rm F,C}$ at solution points via the differentiation operator to advance the solution $\vect{q}_e^{\rm S}$ in time,
\begin{equation}
\der{\vect{q}_e^{\rm S}}{t} = -\tens{\mathcal{M}} \cdot \vect{G}_e^{\rm F,C},
\end{equation}
and obtain the updated value of $\vect{u}_e^{\rm S} = J_e^{-1}\vect{q}_e^{\rm S}$.

\end{enumerate}

\subsection{The compact one-dimensional flux}\label{sec:1d:cmp}

As pointed out in the previous section, the stencil of the continuous auxiliary variable $v$ in equation~\eqref{eq:c:aux:v}, \ie~the gradient of the solution $u$, extends over five elements.
Following what is suggested by~\citet{huynh:09} in the framework of the Flux Reconstruction scheme, a compact stencil over three elements can be obtained by a suitable definition of the continuous flux $H_e$ within the elements sharing the same interface. This can be particularly useful in terms of memory storage for implicit time integrations in DSEMs~\cite{renac2012time,watkins2016multi,bassi2016development,ortleb2020comparative} or for adjoint solves based on the SD scheme~\cite{ou2011unsteady,clinco2026end}.
This approach is here extended for the SD scheme. Notice that for specific choices of correction functions in FR setting, the SD scheme can be recovered. However, this is true only for linear problems such as linear advection or pure diffusion. In the following sections, we will also consider nonlinear problems to fully assess the novel combination of compact viscous fluxes applied to the SD method.

In practice, the reduced stencil is achieved via a modification of the steps~\ref{item:2} and~\ref{item:3}, namely, in the way the continuous flux $\vect{H}_e^{\rm F,C}$ is defined and how is then used to compute the auxiliary solution $\vect{w}_e^{\rm S}$ within the elements sharing any given interface.

Considering equation~\eqref{eq:c:flux}, we now introduce the following continuous fluxes for the left and right interfaces of the $e$-th element:
\begin{align}
\vect{H}_e^{\rm C,L} = \vect{u}_e^{\rm C,L} = \vect{u}_e^{\rm F} + [\hat{H}_{e}^{\rm L} - u_e(-1)]\uvect{c}_{\rm L} = \vect{H}_e^{\rm F,C} - [\hat{H}_{e}^{\rm R} - u_e(1)]\uvect{c}_{\rm R},
\label{eq:c:flux:l}\\
\vect{H}_e^{\rm C,R} = \vect{u}_e^{\rm C,R} = \vect{u}_e^{\rm F} + [\hat{H}_{e}^{\rm R} - u_e(1)]\uvect{c}_{\rm R} = \vect{H}_e^{\rm F,C} - [\hat{H}_{e}^{\rm L} - u_e(-1)]\uvect{c}_{\rm L}.
\label{eq:c:flux:r}
\end{align} 

Note that $\vect{H}_e^{\rm C,L}$ (resp.~$\vect{H}_e^{\rm C,R}$) is not corrected to the right (resp.~to the left) and, as a consequence, depends now on data from elements $e-1$ and $e$ (resp.~elements $e$ and $e+1$).

Concerning the auxiliary variable, equations~\eqref{eq:w:e} and~\eqref{eq:v:e} are now used for the interior flux points only, whereas, for the interface flux points we first set
\begin{align}
\vect{v}_e^{\rm S,L} &= J_e^{-1}\tens{\mathcal{M}} \cdot \vect{H}_e^{\rm C,L} = \vect{v}_e^{\rm S} - J_e^{-1}\tens{\mathcal{M}} \cdot \uvect{c}_{\rm R}[\hat{H}_{e}^{\rm R} - u_e(1)],\label{eq:v:s:l}\\
\vect{v}_e^{\rm S,R} &= J_e^{-1}\tens{\mathcal{M}} \cdot \vect{H}_e^{\rm C,R} = \vect{v}_e^{\rm S} - J_e^{-1}\tens{\mathcal{M}} \cdot \uvect{c}_{\rm L}[\hat{H}_{e}^{\rm L} - u_e(-1)],\label{eq:v:s:r}
\end{align} 
and then interpolate at the interface flux points only:
\begin{equation}
\uvect{c}_{\rm L}^\top \cdot \vect{v}_e^{\rm F,L} = \uvect{c}_{\rm L}^\top \cdot \tens{\mathcal{L}} \cdot \vect{v}_e^{\rm S,L},\quad\text{and}\quad
\uvect{c}_{\rm R}^\top \cdot \vect{v}_e^{\rm F,R} = \uvect{c}_{\rm R}^\top \cdot \tens{\mathcal{L}} \cdot \vect{v}_e^{\rm S,R}.
\label{eq:cmp:int:fp}
\end{equation}

It is worth noting that the above quantities inherit the two-elements stencil of $\vect{H}_e^{\rm C,L}$ and $\vect{H}_e^{\rm C,R}$. Moreover, for implementation purposes, is useful to mention that $\tens{\mathcal{M}} \cdot \uvect{c}_{\rm L}$ and $\tens{\mathcal{M}} \cdot \uvect{c}_{\rm R}$ are the first and last columns of $\tens{\mathcal{M}}$, whereas $\uvect{c}_{\rm L}^\top \cdot \tens{\mathcal{L}}$ and $\uvect{c}_{\rm R}^\top \cdot \tens{\mathcal{L}}$ are the first and last rows of $\tens{\mathcal{L}}$.

Finally, the continuous auxiliary variable is still defined via equation~\eqref{eq:c:aux:v}, but now the common interface values from equation~\eqref{eq:aux:v:int} are replaced with
\begin{equation}
\begin{aligned}
\hat{v}_{e}^{\rm L} &= \tfrac{1}{2}(\uvect{c}_{\rm R}^\top \cdot \vect{v}_{e-1}^{\rm F,R} + \uvect{c}_{\rm L}^\top \cdot \vect{v}_e^{\rm F,L}),\\
\hat{v}_{e}^{\rm R} &= \tfrac{1}{2}(\uvect{c}_{\rm R}^\top \cdot \vect{v}_{e}^{\rm F,R} + \uvect{c}_{\rm L}^\top \cdot \vect{v}_{e+1}^{\rm F,L}).
\label{eq:cmp:int:avg}
\end{aligned}
\end{equation}

In this case, the interface value $\hat{v}_{e}^{\rm L}$ has a two-elements stencil spanning across elements $e-1$ and $e$, and the interface value $\hat{v}_{e}^{\rm R}$ has a two-elements stencil spanning across elements $e$ and $e+1$. Therefore, the resulting stencil for $\vect{v}_e^{\rm F,C}$ now spans over only three elements from $e-1$ to $e+1$. 

It is worth highlighting that, from the implementation point of view, the modified quantities in equations~\eqref{eq:c:flux:l} and~\eqref{eq:c:flux:r} do not need to be computed, as the only terms which are actually needed are the relevant corrections  in equations~\eqref{eq:v:s:l} and~\eqref{eq:v:s:r}. 

\subsection{The compact flux in multiple dimensions}
The extension to multiple dimensions of the compact flux is here detailed for completeness.
The starting equations are~\eqref{eq:ns:mapped} and~\eqref{eq:adj:f}.
In order to keep the notation as clear and simple as possible, we assume a two-dimensional domain whose typical elements are those depicted in figure~\ref{fig:2d:sketch}.
The extension to three dimensions will follow similar steps in a straightforward way.

It is important to point out that the use of tensor elements (\ie, quadrilateral and hexahedral elements) allows applying the relevant operators direction by direction.
As such, if $\tens{u}_e^{\rm S} \in \mathbb{R}^{n \times n}$ is the solution matrix at solution points, in which $\tensij{u}_{e,ij}^{\rm S}$ represents the solution at location $\vect{\xi} = (\xi_i,\;\eta_j)^\top$, then the interpolation operator becomes
\begin{equation}
\tens{u}_e^{\rm F\xi} = \tens{u}_e^{\rm S} \cdot \tens{\mathcal{L}}^\top,\quad\text{and}\quad
\tens{u}_e^{\rm F\eta} = \tens{\mathcal{L}} \cdot \tens{u}_e^{\rm S},\label{eq:int:2d}
\end{equation}
with $\tens{u}_e^{\rm F\xi} \in \mathbb{R}^{n \times (n+1)}$ (resp.~$\tens{u}_e^{\rm F\eta} \in \mathbb{R}^{(n+1) \times n}$) and $\tensij{u}_{e,ij}^{\rm F\xi}$ (resp.~$\tensij{u}_{e,ij}^{\rm F\eta}$) the interpolated solution at location $\vect{\xi} = (\xi_{i-1/2},\;\eta_j)^\top$ (resp.~$\vect{\xi} = (\xi_i,\;\eta_{j-1/2})^\top$).

As far as the differentiation operator is concerned, analogous rules apply:
\begin{equation}
\der{\tens{\phi}_e^{\rm S}}{\xi} = \tens{\phi}_e^{\rm F\xi} \cdot \tens{\mathcal{M}}^\top,\quad\text{and}\quad
\der{\tens{\phi}_e^{\rm S}}{\eta} = \tens{\mathcal{M}} \cdot \tens{\phi}_e^{\rm F\eta}
\end{equation}
where the $ij$-th element of the matrices $\ilder{\tens{\phi}_e^{\rm S}}{\xi}$ and $\ilder{\tens{\phi}_e^{\rm S}}{\eta}$, both $\in \mathbb{R}^{n \times n}$, are the components of the gradient, in computational space, of the quantity $\phi_e$ at location $\vect{\xi} = (\xi_i,\;\eta_j)^\top$.

Finally, we define generic arrays at solution points, $\xi$-aligned flux points and $\eta$-aligned flux points, respectively,
$$
[\;\cdot\;]_e^{\rm S} \in \mathbb{R}^{n \times n},\quad
[\;\cdot\;]_e^{\rm F\xi} \in \mathbb{R}^{n \times (n+1)},\quad\text{and}\quad
[\;\cdot\;]_e^{\rm F\eta}\in \mathbb{R}^{(n+1) \times n},
$$
where it is understood that the operation within brackets is evaluated at the corresponding solution or flux points. For instance, $[J^{-1}]_e^{\rm S}$ is the reciprocal of $J_e = \det(\tens{J}_e)$ evaluated at each solution point, whereas $[\tens{S}]_e^{\rm F\xi}$ is the adjoint of $\tens{J}_e$ evaluated at $\xi$-aligned flux points. In some cases the Hadamard entry-wise product $(\;\odot\;)$ will also be used.

Following the steps described in section~\ref{sec:1d:std}, with the modifications detailed in section~\ref{sec:1d:cmp}, the scheme is built as follows:

\begin{enumerate}
\item The discontinuous solution at flux points is first computed using equation~\eqref{eq:int:2d}.

\item The discontinuous solution is then made continuous with the following relations for $\xi$-aligned flux points:
\begin{align}
\tens{u}_{e}^{\rm F\xi,C} &= \tens{u}_{e}^{\rm F\xi} + \uvect{c}_{\rm L} \cdot [\uvect{u}_{e}^{\rm L} - \tens{u}_{e}^{\rm F\xi}\cdot \uvect{c}_{\rm L}]^\top + \uvect{c}_{\rm R} \cdot [\uvect{u}_{e}^{\rm R} - \tens{u}_{e}^{\rm F\xi}\cdot \uvect{c}_{\rm R}]^\top,\label{eq:c:sol:2d:xi}\\
\tens{u}_{e}^{\rm F\xi,L} &= \tens{u}_{e}^{\rm F\xi,C} - \uvect{c}_{\rm R} \cdot [\uvect{u}_{e}^{\rm R} - \tens{u}_{e}^{\rm F\xi}\cdot \uvect{c}_{\rm R}]^\top,\\
\tens{u}_{e}^{\rm F\xi,R} &= \tens{u}_{e}^{\rm F\xi,C} - \uvect{c}_{\rm L} \cdot [\uvect{u}_{e}^{\rm L} - \tens{u}_{e}^{\rm F\xi}\cdot \uvect{c}_{\rm L}]^\top,
\end{align}
with
\begin{equation}
\uvect{u}_{e}^{\rm L} = \tfrac{1}{2}(\tens{u}_{e-}^{\rm F\xi} \cdot \uvect{c}_{\rm R} + \tens{u}_e^{\rm F\xi} \cdot \uvect{c}_{\rm L}),\quad
\uvect{u}_{e}^{\rm R} = \tfrac{1}{2}(\tens{u}_{e}^{\rm F\xi} \cdot \uvect{c}_{\rm R} + \tens{u}_{e+}^{\rm F\xi} \cdot \uvect{c}_{\rm L}),
\end{equation}
where the subscripts $e-$ and $e+$ identify the elements to the left and to the right of the $e$-th element, respectively.
Similar relations are obtained for $\eta$-aligned flux points, with corrections applied at the bottom and top interfaces, respectively:
\begin{align}
\tens{u}_{e}^{\rm F\eta,C} &= \tens{u}_{e}^{\rm F\eta} + \uvect{c}_{\rm L} \cdot [\uvect{u}_{e}^{\rm B} - \uvect{c}_{\rm L}^\top \cdot \tens{u}_{e}^{\rm F\eta}] + \uvect{c}_{\rm R} \cdot [\uvect{u}_{e}^{\rm T} - \uvect{c}_{\rm R}^\top \cdot \tens{u}_{e}^{\rm F\eta}],\label{eq:c:sol:2d:eta}\\
\tens{u}_{e}^{\rm F\eta,B} &= \tens{u}_{e}^{\rm F\eta,C} - \uvect{c}_{\rm R} \cdot [\uvect{u}_{e}^{\rm T} - \uvect{c}_{\rm R}^\top \cdot \tens{u}_{e}^{\rm F\eta}],\\
\tens{u}_{e}^{\rm F\eta,T} &= \tens{u}_{e}^{\rm F\eta,C} - \uvect{c}_{\rm L} \cdot [\uvect{u}_{e}^{\rm B} - \uvect{c}_{\rm L}^\top \cdot \tens{u}_{e}^{\rm F\eta}],
\end{align}
with
\begin{equation}
\uvect{u}_{e}^{\rm B} = \tfrac{1}{2}( \uvect{c}_{\rm R}^\top \cdot \tens{u}_{e-}^{\rm F\eta} + \uvect{c}_{\rm L}^\top \cdot \tens{u}_{e}^{\rm F\eta}),\quad
\uvect{u}_{e}^{\rm T} = \tfrac{1}{2}( \uvect{c}_{\rm R}^\top \cdot \tens{u}_{e}^{\rm F\eta} + \uvect{c}_{\rm L}^\top \cdot \tens{u}_{e+}^{\rm F\eta}),
\end{equation}
where, this time, the subscripts $e-$ and $e+$ identify the elements at the bottom and at the top of the $e$-th element, respectively.

Thus, from the above solution values and the adjoint $\tens{S}$, which, it is worth noting, is intrinsically continuous across elements' interfaces, the continuous fluxes are readily obtained:
\begin{equation}
\tens{H}_{m,e}^{\rm F\xi,C} = \tens{u}_{e}^{\rm F\xi,C} \odot [\tensij{S}_{1m}]_e^{\rm F\xi},\quad
\tens{H}_{m,e}^{\rm F\xi,L} = \tens{u}_{e}^{\rm F\xi,L} \odot [\tensij{S}_{1m}]_e^{\rm F\xi},\quad
\tens{H}_{m,e}^{\rm F\xi,R} = \tens{u}_{e}^{\rm F\xi,R} \odot [\tensij{S}_{1m}]_e^{\rm F\xi},
\end{equation}
\begin{equation}
\tens{H}_{m,e}^{\rm F\eta,C} = \tens{u}_{e}^{\rm F\eta,C} \odot [\tensij{S}_{2m}]_e^{\rm F\eta},\quad
\tens{H}_{m,e}^{\rm F\eta,B} = \tens{u}_{e}^{\rm F\eta,B} \odot [\tensij{S}_{2m}]_e^{\rm F\eta},\quad
\tens{H}_{m,e}^{\rm F\eta,T} = \tens{u}_{e}^{\rm F\eta,T} \odot [\tensij{S}_{2m}]_e^{\rm F\eta},
\end{equation}
where, in the two-dimensional case we are focusing on, the subscript $m = 1,\, 2$.

\item In order to obtain the auxiliary variable at the flux points, namely $\grad{u}$, the differentiation operator is first applied to the above continuous fluxes and the result is then interpolated at the flux points. Again, to obtain the compact scheme, distinction must be made between interior and interface flux points.

Concerning the former, we have
\begin{equation}
\tens{v}_{m,e}^{\rm S} = [J^{-1}]_e^{\rm S} \odot (\tens{H}_{m,e}^{\rm F\xi,C} \cdot \tens{\mathcal{M}}^\top + \tens{\mathcal{M}} \cdot \tens{H}_{m,e}^{\rm F\eta,C} ),
\end{equation}
which is then interpolated only to the interior flux points:
\begin{equation}
\tens{v}_{m,e}^{\rm F\xi} = \tens{v}_{m,e}^{\rm S} \cdot \tens{\mathcal{L}}^\top,\quad
\tens{v}_{m,e}^{\rm F\eta} = \tens{\mathcal{L}} \cdot \tens{v}_{m,e}^{\rm S}.\label{eq:aux:v:std:2d}
\end{equation}

Focusing instead on the element's interfaces, we set 
\begin{align}
\tens{v}_{m,e}^{\rm S,L} &= [J^{-1}]_e^{\rm S} \odot (\tens{H}_{m,e}^{\rm F\xi,L} \cdot \tens{\mathcal{M}}^\top + \tens{\mathcal{M}} \cdot \tens{H}_{m,e}^{\rm F\eta,C} ),\label{eq:aux:v:cmp:l}\\
\tens{v}_{m,e}^{\rm S,R} &= [J^{-1}]_e^{\rm S} \odot (\tens{H}_{m,e}^{\rm F\xi,R} \cdot \tens{\mathcal{M}}^\top + \tens{\mathcal{M}} \cdot \tens{H}_{m,e}^{\rm F\eta,C} ),\\
\tens{v}_{m,e}^{\rm S,B} &= [J^{-1}]_e^{\rm S} \odot (\tens{H}_{m,e}^{\rm F\xi,C} \cdot \tens{\mathcal{M}}^\top + \tens{\mathcal{M}} \cdot \tens{H}_{m,e}^{\rm F\eta,B} ),\\
\tens{v}_{m,e}^{\rm S,T} &= [J^{-1}]_e^{\rm S} \odot (\tens{H}_{m,e}^{\rm F\xi,C} \cdot \tens{\mathcal{M}}^\top + \tens{\mathcal{M}} \cdot \tens{H}_{m,e}^{\rm F\eta,T} ),\label{eq:aux:v:cmp:t}
\end{align}
which are then interpolated at the interfaces:
\begin{equation}
\tens{v}_{m,e}^{\rm F\xi,L} \cdot \uvect{c}_{\rm L} = \tens{v}_{m,e}^{\rm S,L} \cdot \tens{\mathcal{L}}^\top \cdot \uvect{c}_{\rm L},\quad
\tens{v}_{m,e}^{\rm F\xi,R} \cdot \uvect{c}_{\rm R} = \tens{v}_{m,e}^{\rm S,R} \cdot \tens{\mathcal{L}}^\top \cdot \uvect{c}_{\rm R},\label{eq:aux:cmp:lr:2d}
\end{equation}
\begin{equation}
\uvect{c}_{\rm L}^\top \cdot \tens{v}_{m,e}^{\rm F\eta,B} = \uvect{c}_{\rm L}^\top \cdot \tens{\mathcal{L}} \cdot \tens{v}_{m,e}^{\rm S,B},\quad 
\uvect{c}_{\rm R}^\top \cdot \tens{v}_{m,e}^{\rm F\eta,T} = \uvect{c}_{\rm R}^\top \cdot \tens{\mathcal{L}} \cdot \tens{v}_{m,e}^{\rm S,T},\label{eq:aux:cmp:bt:2d}
\end{equation}

It is worthwhile pointing out that, for each family of flux points, in equations~\eqref{eq:aux:v:cmp:l}--\eqref{eq:aux:v:cmp:t} the corrected compact fluxes are differentiated along the corresponding direction, whereas the standard fluxes are differentiated along the transverse direction. As a result, the global stencil span five elements, namely, the central element and its four neighbour.

\item To compute the flux $\vect{G}_e$, the auxiliary variable is made continuous across interfaces using simple interface averages involving values from equations~\eqref{eq:aux:v:std:2d}, \eqref{eq:aux:cmp:lr:2d} and~\eqref{eq:aux:cmp:bt:2d}:
\begin{equation}
\begin{aligned}
\tens{v}_{m,e}^{\rm F\xi,C} &= \tens{v}_{m,e}^{\rm F\xi} + \uvect{v}_{m,e}^{\rm L} \cdot \uvect{c}_{\rm L}^\top + \uvect{v}_{m,e}^{\rm R} \cdot \uvect{c}_{\rm R}^\top,\\
\tens{v}_{m,e}^{\rm F\eta,C} &= \tens{v}_{m,e}^{\rm F\eta} + \uvect{c}_{\rm L} \cdot (\uvect{v}_{m,e}^{\rm B})^\top + \uvect{c}_{\rm R} \cdot (\uvect{v}_{m,e}^{\rm T})^\top,
\end{aligned}\label{eq:c:aux:v:2d}
\end{equation}
with
\begin{equation}
\begin{aligned}
\uvect{v}_{m,e}^{\rm L} = \tfrac{1}{2}(\tens{v}_{m,e-}^{\rm F\xi,R} \cdot \uvect{c}_{\rm R} + \tens{v}_{m,e}^{\rm F\xi,L} \cdot \uvect{c}_{\rm L}),&\quad
\uvect{v}_{m,e}^{\rm R} = \tfrac{1}{2}(\tens{v}_{m,e}^{\rm F\xi,R} \cdot \uvect{c}_{\rm R} + \tens{v}_{m,e+}^{\rm F\xi,L} \cdot \uvect{c}_{\rm L}),\\
\uvect{v}_{m,e}^{\rm B} = \tfrac{1}{2}( \uvect{c}_{\rm R}^\top \cdot \tens{v}_{m,e-}^{\rm F\eta,T} + \uvect{c}_{\rm L}^\top \cdot \tens{v}_{m,e}^{\rm F\eta,B})^\top,&\quad
\uvect{v}_{m,e}^{\rm T} = \tfrac{1}{2}( \uvect{c}_{\rm R}^\top \cdot \tens{v}_{m,e}^{\rm F\eta,T} + \uvect{c}_{\rm L}^\top \cdot \tens{v}_{m,e+}^{\rm F\eta,B})^\top.
\end{aligned}
\end{equation}

Notice that equation~\eqref{eq:c:aux:v:2d} makes use of the fact that the interpolation in equation~\eqref{eq:aux:v:std:2d} was performed at interior flux points only, \ie, $\tens{v}_{m,e}^{\rm F\xi}$ (resp.~$\tens{v}_{m,e}^{\rm F\eta}$) has the first and last columns (resp.~rows) of zeros.
Eventually, using the continuous solution from equations~\eqref{eq:c:sol:2d:xi} and~\eqref{eq:c:sol:2d:eta} and the continuous auxiliary variable from equation~\eqref{eq:c:aux:v:2d}, the diffusive flux can be computed at all flux points:
\begin{equation}
\tens{G}_{e}^{\rm F\xi,C} = [\tens{S} \cdot \vect{F}(u^{\rm F\xi,C},\vect{v}^{\rm F\xi,C})]_e^{\rm F\xi},\quad\text{and}\quad
\tens{G}_{e}^{\rm F\eta,C} = [\tens{S} \cdot \vect{F}(u^{\rm F\eta,C},\vect{v}^{\rm F\eta,C})]_e^{\rm F\eta}.
\end{equation}

\item Finally, the flux divergence is computed at solution points to advance the solution in time,
\begin{equation}
\der{\tens{q}_e^{\rm S}}{t} = -(\tens{G}_e^{\rm F\xi,C} \cdot \tens{\mathcal{M}}^\top + \tens{\mathcal{M}} \cdot \tens{G}_e^{\rm F\eta,C}),
\end{equation}
and obtain the updated value of $\tens{u}_e^{\rm S} = [J^{-1}]_e^{\rm S} \odot \tens{q}_e^{\rm S}$.

\end{enumerate}

\section{Temporal Eigenanalysis}\label{sec:2}
As mentioned in the introduction, a large variety of numerical techniques to assess the intrinsic characteristics of high-order discontinuous spectral element methods have been proposed in the literature~\citep{lele1992compact,bogey2004family,van2008stability,vincent2011insights,moura2015linear,vanharen2017revisiting}. In this section we consider the classical standard temporal eigenanalysis applied to the newly proposed compact formulation of the SD scheme for the diffusion equation.
In the standard temporal analysis, the diffusion equation is discretized looking for wave-like solutions in order to study the temporal evolution of specific wavenumbers. Then, dispersion and dissipation properties of any scheme follow directly from the corresponding eigensolutions.

Let us now consider the case of a constant unitary viscosity: 
\begin{equation}
\pd{u}{t}-\frac{\partial^{2}u}{\partial x^{2}}=0.
\label{advection}
\end{equation}
This equation admits plane wave solutions of the form 
\begin{equation}
u(x,t)=e^{\iota(\theta x-\omega t)},\quad \text{with $\iota^{2} = -1$},
\label{wave}
\end{equation}
provided that the angular frequency $\omega=\omega(\theta)$ is such that
\begin{equation}
\textrm{Re}(\omega)=0\qquad \textrm{and} \qquad \textrm{Im}(\omega)=-\theta^{2}, 
\label{ReIm}
\end{equation}
where $\theta$ is a real-valued wavenumber chosen at the initial condition. Equation~\eqref{ReIm} provides what are respectively known as dispersion and diffusion relations for the exact plane wave solution.

Following exactly the same steps presented in the previous section we can obtain a semi-discrete temporal eigenanalysis formulation for the proposed scheme.

As a first step in the computation of the PDE residual at the solution points, we begin with the discrete representation of the solution within each spectral element, which is entirely determined by the nodal values of the solution evaluated at the solution points. The second step in the computation of the second-order flux consists of extrapolating these values to the flux points, which is performed using the matrix $\boldsymbol{\mathcal{L}}$:
\begin{equation}
\vect{u}_e^{\rm F} = \tens{\mathcal{L}} \cdot \vect{u}_e^{\rm S}.
\end{equation}
Since the goal is to correct the values of $\textbf{v}_{e}^{\mathrm{S}}$ in equations \eqref{eq:v:s:l} and \eqref{eq:v:s:r}, we follow the standard approach of first evaluating these values and then correcting them using information from neighboring elements, as shown in equations \eqref{eq:v:s:l} and \eqref{eq:v:s:r}.
So, we can write equation \eqref{eq:c:flux} as:
\begin{equation}
\begin{aligned}  
\vect{u}_e^{\rm F,C} &= \vect{u}_e^{\rm F} + [\hat{H}_{e}^{\rm L} - u_e(-1)]\uvect{c}_{\rm L} + [\hat{H}_{e}^{\rm R} - u_e(1)]\uvect{c}_{\rm R} \\
&= \tens{\mathcal{L}}\cdot\vect{u}_e^{\rm S}+\bigg[ \frac{1}{2}(\uvect{c}_{\rm R}^\top \cdot \vect{u}_{e-1}^{\rm F}+\uvect{c}_{\rm L}^\top \cdot \vect{u}_{e}^{\rm F})-\uvect{c}_{\rm L}^\top \cdot \vect{u}_{e}^{\rm F}\bigg]\uvect{c}_{\rm L}+\bigg[\frac{1}{2}(\uvect{c}_{\rm R}^\top \cdot \vect{u}_{e}^{\rm F}+\uvect{c}_{\rm L}^\top \cdot \vect{u}_{e+1}^{\rm F})-\uvect{c}_{\rm R}^\top \cdot \vect{u}_{e}^{\rm F}\bigg]\uvect{c}_{\rm R}  \\
&= \tens{\mathcal{L}}\cdot\vect{u}_e^{\rm S}+\frac{1}{2}\big(\uvect{c}_{\rm R}^\top \cdot \vect{u}_{e-1}^{\rm F}-\uvect{c}_{\rm L}^\top \cdot \vect{u}_{e}^{\rm F}\big)\uvect{c}_{\rm L}+\frac{1}{2}\big(\uvect{c}_{\rm L}^\top \cdot \vect{u}_{e+1}^{\rm F}-\uvect{c}_{\rm R}^\top \cdot \vect{u}_{e}^{\rm F}\big)\uvect{c}_{\rm R}  \\
&= \tens{\mathcal{L}}\cdot\vect{u}_e^{\rm S}+\frac{1}{2}\big(\uvect{c}_{\rm R}^\top \cdot \tens{\mathcal{L}} \cdot\vect{u}_{e-1}^{\rm S}-\uvect{c}_{\rm L}^\top \cdot \tens{\mathcal{L}} \cdot\vect{u}_{e}^{\rm S}\big)\uvect{c}_{\rm L}+\frac{1}{2}\big(\uvect{c}_{\rm L}^\top \cdot \tens{\mathcal{L}}\cdot\vect{u}_{e+1}^{\rm S}-\uvect{c}_{\rm R}^\top \cdot \tens{\mathcal{L}} \cdot\vect{u}_{e}^{\rm S}\big)\uvect{c}_{\rm R},
\end{aligned}
\end{equation}
and, by introducing the classical assumption of temporal eigenanalysis which consists of Bloch periodic solutions fulfilling the condition
\begin{equation}
    \vect{u}_{e\pm1}^{\rm S} = \vect{u}_{e}^{\rm S} e^{\pm \theta \iota},
\end{equation}
we can obtain the following expression:
\begin{equation}
\begin{aligned} 
\vect{u}_e^{\rm F,C} &= \tens{\mathcal{L}}\cdot\vect{u}_e^{\rm S}+\frac{1}{2}\big(\uvect{c}_{\rm R}^\top \cdot \tens{\mathcal{L}} \cdot\vect{u}_{e}^{\rm S} e^{-\theta \iota}-\uvect{c}_{\rm L}^\top \cdot \tens{\mathcal{L}} \cdot\vect{u}_{e}^{\rm S}\big)\uvect{c}_{\rm L}+\frac{1}{2}\big(\uvect{c}_{\rm L}^\top \cdot \tens{\mathcal{L}}\cdot\vect{u}_{e}^{\rm S} e^{+\theta \iota}-\uvect{c}_{\rm R}^\top  \cdot \tens{\mathcal{L}} \cdot\vect{u}_{e}^{\rm S}\big)\uvect{c}_{\rm R} \\
&= \bigg[\tens{\mathcal{L}} + \frac{1}{2} \uvect{c}_{\rm L}\cdot\uvect{c}_{\rm R}^{\top} \cdot\tens{\mathcal{L}}e^{-\theta \iota} - \frac{1}{2}\uvect{c}_{\rm L}\cdot\uvect{c}_{\rm L}^{\top}\cdot\tens{\mathcal{L}}+\frac{1}{2} \uvect{c}_{\rm R}\cdot\uvect{c}_{\rm L}^{\top} \cdot\tens{\mathcal{L}}e^{+\theta \iota} - \frac{1}{2}\uvect{c}_{\rm R}\cdot\uvect{c}_{\rm R}^{\top}\cdot\tens{\mathcal{L}}\bigg]\cdot\vect{u}_e^{\rm S}=\tens{\mathcal{H}}\cdot\vect{u}_e^{\rm S}.
\end{aligned}
\end{equation}

Now, let us continue to implement the standard approach knowing that in the computation of $\vect{v}_{e}^{\rm S}$ we will add a correction based on neighboring elements. The standard definition of $\vect{v}_{e}^{\rm S}$ can be written as:
\begin{equation}
    \vect{v}_{e}^{\rm S} = J_{e}^{-1} \tens{\mathcal{M}}\cdot\vect{u}_e^{\rm F,C}=J_{e}^{-1} \tens{\mathcal{M}} \cdot\tens{\mathcal{H}}\cdot\vect{u}_e^{\rm S}.
\end{equation}
We now need to add the corrections of equations \eqref{eq:v:s:l} and \eqref{eq:v:s:r}. We will also need to separate the treatment of left and right faces by subtracting the correction of the opposite face. Up to this point, in fact, the derivation is exactly the same as the one for the standard fully centered scheme.

We then apply the aforementioned corrections as:
\begin{equation}
\begin{aligned} 
\vect{v}_e^{\rm S,L} &= J_{e}^{-1} \tens{\mathcal{M}} \cdot\tens{\mathcal{H}}\cdot\vect{u}_e^{\rm S} - J_e^{-1}\tens{\mathcal{M}} \cdot \uvect{c}_{\rm R}[\hat{H}_{e}^{\rm R} - u_e(1)], \\
&= J_{e}^{-1} \tens{\mathcal{M}}\cdot\bigg[ \tens{\mathcal{H}} -\frac{1}{2}\uvect{c}_{\rm R}\cdot\uvect{c}_{\rm L}^{\top}\cdot\tens{\mathcal{L}}e^{+\theta \iota} + \frac{1}{2}\uvect{c}_{\rm R}\cdot\uvect{c}_{\rm R}^{\top}\cdot\tens{\mathcal{L}}\bigg]\cdot\vect{u}_e^{\rm S}=\tens{\mathcal{W}}^{\rm L}\cdot\vect{u}_e^{\rm S}.
\end{aligned}
\end{equation}
This expression provides a correction to the standard approach associated with the matrix $\tens{\mathcal{H}}$, obtained by subtracting the contribution coming from the opposite face.

In the same way it is possible to obtain the expression for $\vect{v}_e^{\rm S,R}$ as:
\begin{align} 
\vect{v}_e^{\rm S,R} = J_{e}^{-1} \tens{\mathcal{M}}\cdot\bigg[ \tens{\mathcal{H}} -\frac{1}{2}\uvect{c}_{\rm L}\cdot\uvect{c}_{\rm R}^{\top}\cdot\tens{\mathcal{L}}e^{-\theta \iota} + \frac{1}{2}\uvect{c}_{\rm L}\cdot\uvect{c}_{\rm L}^{\top}\cdot\tens{\mathcal{L}}\bigg]\cdot\vect{u}_e^{\rm S}=\tens{\mathcal{W}}^{\rm R}\cdot\vect{u}_e^{\rm S}.
\end{align}
Equation~\eqref{eq:cmp:int:fp} can then be written as
\begin{equation}
\uvect{c}_{\rm L}^\top \cdot \vect{v}_e^{\rm F,L} = \uvect{c}_{\rm L}^\top \cdot\tens{\mathcal{L}} \cdot\tens{\mathcal{W}}^{\rm L}\cdot\vect{u}_e^{\rm S},\quad\text{and}\quad
\uvect{c}_{\rm R}^\top \cdot \vect{v}_e^{\rm F,R} = \uvect{c}_{\rm R}^\top \cdot\tens{\mathcal{L}} \cdot\tens{\mathcal{W}}^{\rm R}\cdot\vect{u}_e^{\rm S}.
\end{equation}
whereas equation~\eqref{eq:cmp:int:avg} reads:
\begin{equation}
\begin{aligned}
\hat{v}_{e}^{\rm L} &= \tfrac{1}{2}(\uvect{c}_{\rm R}^\top \cdot \tens{\mathcal{L}} \cdot \tens{\mathcal{W}}^{\rm R} e^{-\theta \iota} + \uvect{c}_{\rm L}^\top \cdot \tens{\mathcal{L}} \cdot \tens{\mathcal{W}}^{\rm L}) \cdot \vect{u}_{e}^{\rm S},\\
\hat{v}_{e}^{\rm R} &= \tfrac{1}{2}(\uvect{c}_{\rm R}^\top \cdot \tens{\mathcal{L}} \cdot \tens{\mathcal{W}}^{\rm R} + \uvect{c}_{\rm L}^\top \cdot \tens{\mathcal{L}} \cdot \tens{\mathcal{W}}^{\rm L} e^{+\theta \iota}) \cdot \vect{u}_{e}^{\rm S}.
\end{aligned}
\label{eq:cmp:int:avg:2d}
\end{equation}
%
%
%
Equation~\eqref{eq:c:aux:v}, in the case of the compact formulation, is thus rewritten as:
\begin{equation}
\begin{aligned}
\vect{v}_e^{\rm F,C} &= \vect{v}_e^{\rm F} + [\hat{v}_{e}^{\rm L} - v_e(-1)]\uvect{c}_{\rm L} + [\hat{v}_{e}^{\rm R} - v_e(1)]\uvect{c}_{\rm R} \\
&= \bigg[ \tens{\mathcal{L}}\cdot (J_{e}^{-1}\tens{\mathcal{M}}\cdot\tens{\mathcal{H}}) - \uvect{c}_{\rm L}\cdot\uvect{c}_{\rm L}^{\top} \cdot\tens{\mathcal{L}}\cdot (J_{e}^{-1} \tens{\mathcal{M}}\cdot \tens{\mathcal{H}}) - \uvect{c}_{\rm R}\cdot\uvect{c}_{\rm R}^{\top} \cdot\tens{\mathcal{L}}\cdot (J_{e}^{-1} \tens{\mathcal{M}} \cdot\tens{\mathcal{H}})\\
&+ \tfrac{1}{2}\uvect{c}_{\rm L}\cdot\uvect{c}_{\rm R}^\top \cdot\tens{\mathcal{L}} \cdot\tens{\mathcal{W}}^{\rm R}e^{-\theta \iota} + \tfrac{1}{2}\uvect{c}_{\rm L}\cdot\uvect{c}_{\rm L}^\top \cdot\tens{\mathcal{L}} \cdot\tens{\mathcal{W}}^{\rm L}\\
&+ \tfrac{1}{2}\uvect{c}_{\rm R}\cdot\uvect{c}_{\rm L}^\top  \cdot\tens{\mathcal{L}} \cdot\tens{\mathcal{W}}^{\rm L}e^{+\theta \iota} + \tfrac{1}{2}\uvect{c}_{\rm R}\cdot\uvect{c}_{\rm R}^\top \cdot\tens{\mathcal{L}} \cdot\tens{\mathcal{W}}^{\rm R}\bigg] \cdot\vect{u}_{e}^{\rm S}
= \tens{\mathcal{Q}}\cdot\vect{u}_e^{\rm S}.
\end{aligned}
\label{eq:c:aux:v:cmp}
\end{equation}
%

In case the scheme includes interior penalty terms, these must be accounted for in $\hat{v}_{e}^{\rm L}$ and $\hat{v}_{e}^{\rm R}$. Accordingly, equations~\eqref{eq:cmp:int:avg:2d} are augmented with the following quantities, respectively:
\begin{equation}
\tau ( \uvect{c}_{\rm L}^\top-\uvect{c}_{\rm R}^\top e^{-\iota \theta})\cdot\tens{\mathcal{L}} \cdot\vect{u}_{e}^{\rm S},\quad\text{and}\quad
\tau (\uvect{c}_{\rm L}^\top e^{+\iota \theta} -\uvect{c}_{\rm R}^\top)\cdot\tens{\mathcal{L}} \cdot\vect{u}_{e}^{\rm S},
\end{equation}
where $\tau=\eta_{\rm IP} (p + 1)^{2}/h$ is the optimal IP coefficient for one-dimensional computations~\citep{shahbazi:05,manzanero:18}. In the standard temporal eigenanalysis, the value of $h$ is simply equal to $2$ as we work in the reference element with $\xi_{i} \in[-1,1]$. In the general, multi-dimensional case, $h$ is evaluated, at each interface, from the relation $1/h = \norm{\vect{S}} \max(J_{e-}^{-1},J_{e+}^{-1})$, where $J_{e-}^{-1}$ and $J_{e+}^{-1}$ are the integrated Jacobians within the elements sharing the interface and $\vect{S}$ the relevant surface vector obtained from the integrated adjoint $\tens{S}$ over the interface's flux points~\citep{manzanero:18}.
As a result, equation~\eqref{eq:c:aux:v:cmp} becomes:
\begin{equation}
    \vect{v}_e^{\rm F,C} = (\tens{\mathcal{Q}} + \tens{\mathcal{Q}}_{\rm IP})\cdot\vect{u}_e^{\rm S},
\end{equation}
with 
\begin{equation}
\tens{\mathcal{Q}}_{\rm IP} = \tau ( \uvect{c}_{\rm L}\cdot\uvect{c}_{\rm L}^\top - \uvect{c}_{\rm L}\cdot\uvect{c}_{\rm R}^\top e^{-\theta \iota} ) \cdot\tens{\mathcal{L}} 
 + \tau (  \uvect{c}_{\rm R}\cdot\uvect{c}_{\rm L}^\top e^{+\theta \iota} - \uvect{c}_{\rm R}\cdot\uvect{c}_{\rm R}^\top ) \cdot\tens{\mathcal{L}}.
\end{equation}

It is then possible to directly apply the differentiation matrix, thereby obtaining the discrete set of ODEs governing the dynamics of the SD discretization of the diffusion equation:
\begin{equation}
    \der{\vect{u}_e^{\rm S}}{t} = -\tens{\mathcal{M}} \cdot (\tens{\mathcal{Q}} + \tens{\mathcal{Q}}_{\rm IP}) \cdot\vect{u}_e^{\rm S}=\tens{\mathcal{B}}\cdot\vect{u}_e^{\rm S},
\end{equation}
which, combined with the semi-discrete eigenanalysis formulation for the time derivative, namely,
\begin{equation}
    \der{\vect{u}_e^{\rm S}}{t} = -\iota \widetilde{\omega} \vect{u}_e^{\rm S},
\end{equation}
leads to the final form of the eigenvalue problem to be studied:
\begin{equation}
-\iota \widetilde{\omega} \vect{u}_e^{\rm S} = \tens{\mathcal{B}}(\theta)\cdot\vect{u}_e^{\rm S}.
\label{eq:eigen}
\end{equation}
Notice that, in this simplified one-dimensional setting where eigenanalysis is commonly performed, the Jacobian reduces to $J_{e}^{-1}=1/2$, and one would typically encounter this scaling applied twice, leading to a factor of $4$. This factor commonly appears in most works performing temporal eigenanalysis for discontinuous spectral element methods (see, for example, ~\cite{watkins2016numerical}). Here, however, it is not needed, since the inverse Jacobian has already been incorporated into the formulation.

In standard temporal eigenanalysis, the objective is to prescribe a value of $\theta$, solve the eigenvalue problem \eqref{eq:eigen}, and determine $\widetilde{\omega}(\theta)$. We use the notation $\widetilde{(\;\cdot\;)}$ to emphasize that this quantity represents the numerical approximation of the exact value $\omega$ appearing in equation \eqref{ReIm}. Notice that, for each value of $\theta$, the dimensionality of the eigenvalue problem is $p+1$. As a consequence, for any given value of $\theta$, we obtain $p+1$ eigenvalues. In the literature on temporal eigenanalysis of DSEMs, there has been some debate regarding the interpretation of these $p+1$ eigenvalues. In the vast majority of studies, it is commonly assumed that only one mode represents the physical mode, while the others are disregarded as \emph{spurious}~\citep{hu1999analysis,van2007dispersion,vincent2011insights}. More recent research, instead, has introduced a more rigorous formalism in which the contribution of all modes can be agglomerated into a single combined mode~\citep{alhawwary2018fourier,alhawwary2019study,alhawwary2020combined,moura2024joint}. In order to provide the most general picture, in this work we show all the modes obtained from the eigenvalue problem \eqref{eq:eigen} without imposing any hierarchy. However, for completeness, the same results obtained through the combined-mode analysis are reported in Appendix~\ref{sec:5}.

\begin{figure}[h!]
 \centering  
 \begin{subfigure}{0.48\textwidth}
     \includegraphics[width=\textwidth]{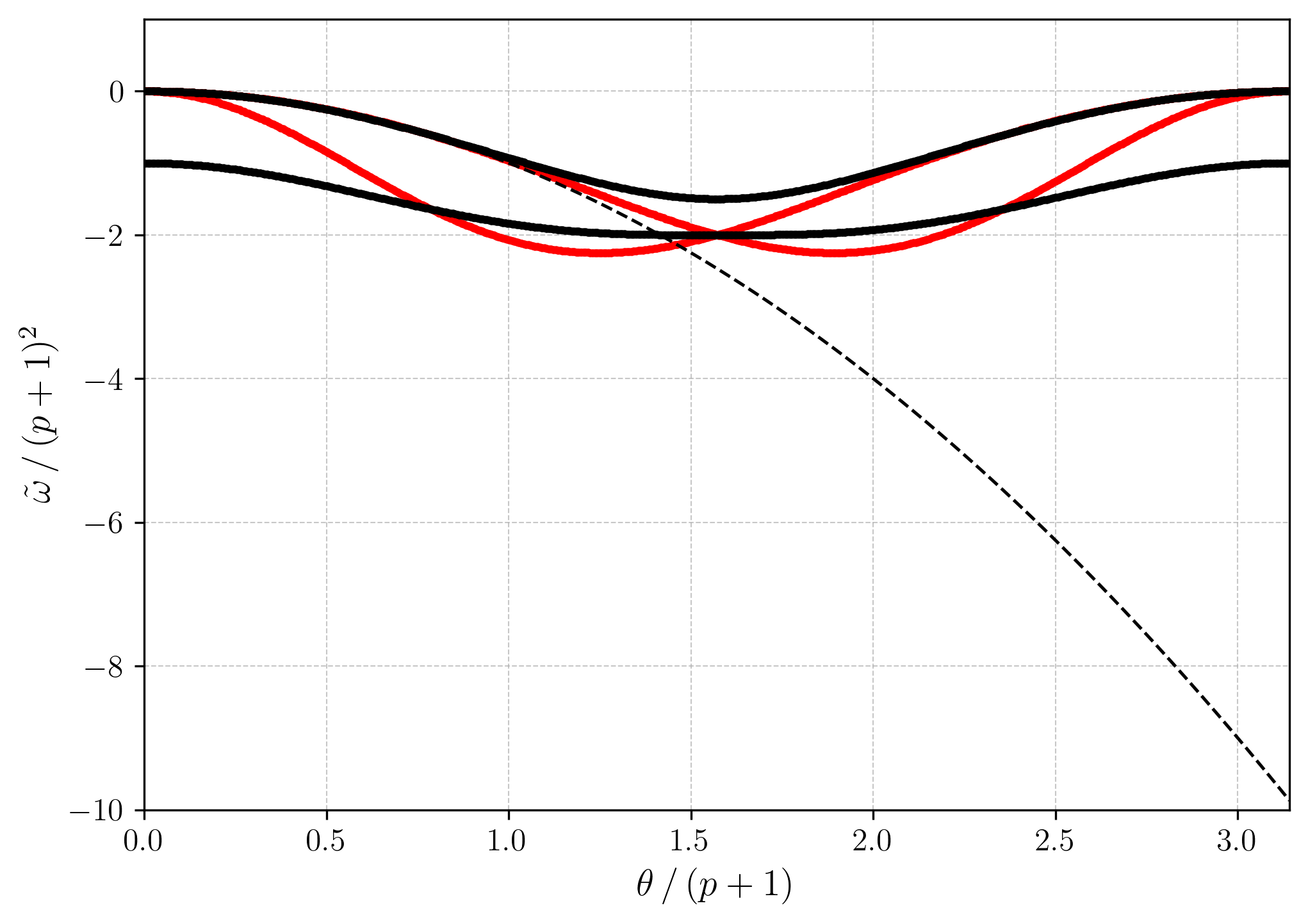}
 \end{subfigure}
 \begin{subfigure}{0.48\textwidth}
     \includegraphics[width=\textwidth]{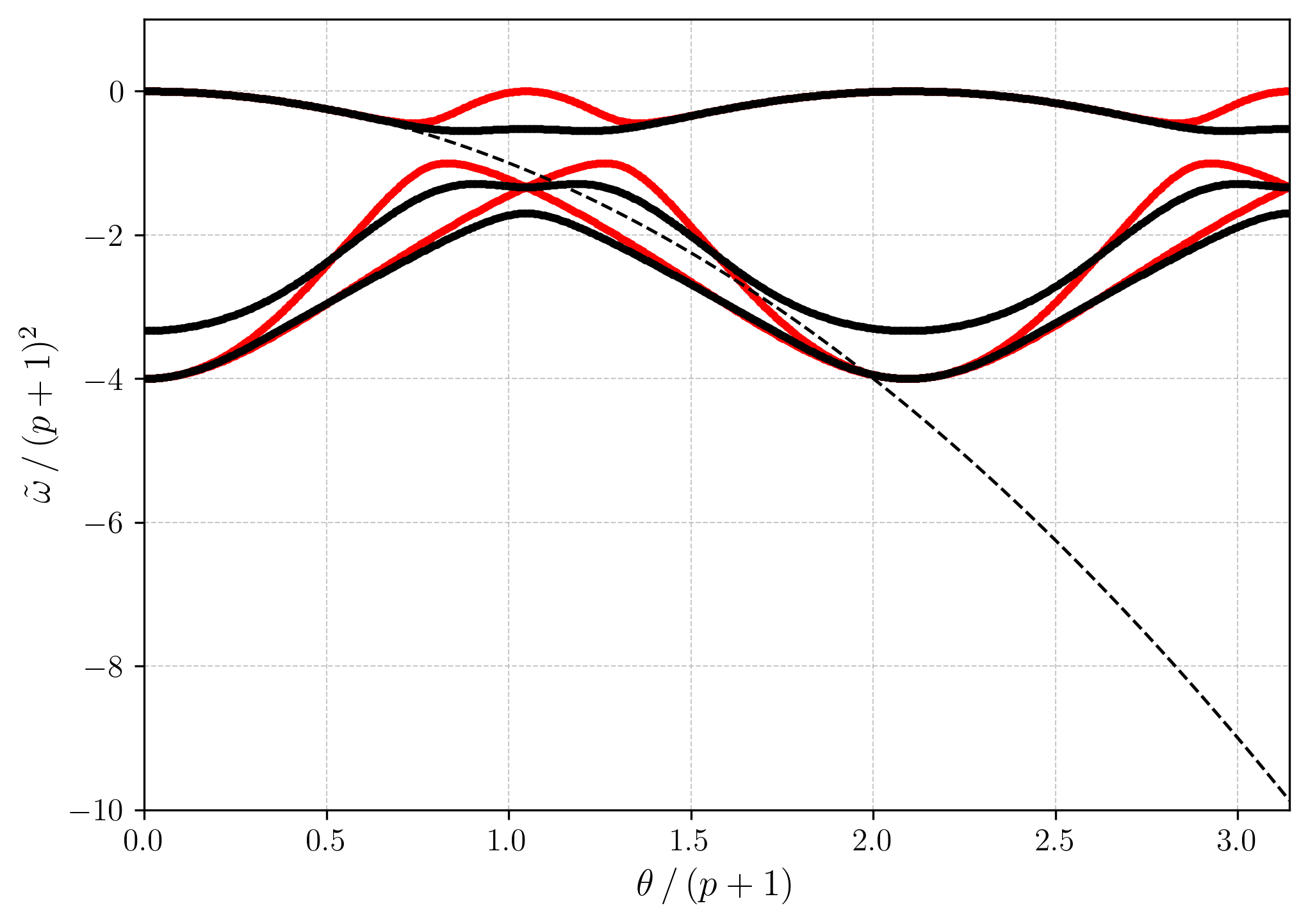}
 \end{subfigure}
 \begin{subfigure}{0.48\textwidth}
     \includegraphics[width=\textwidth]{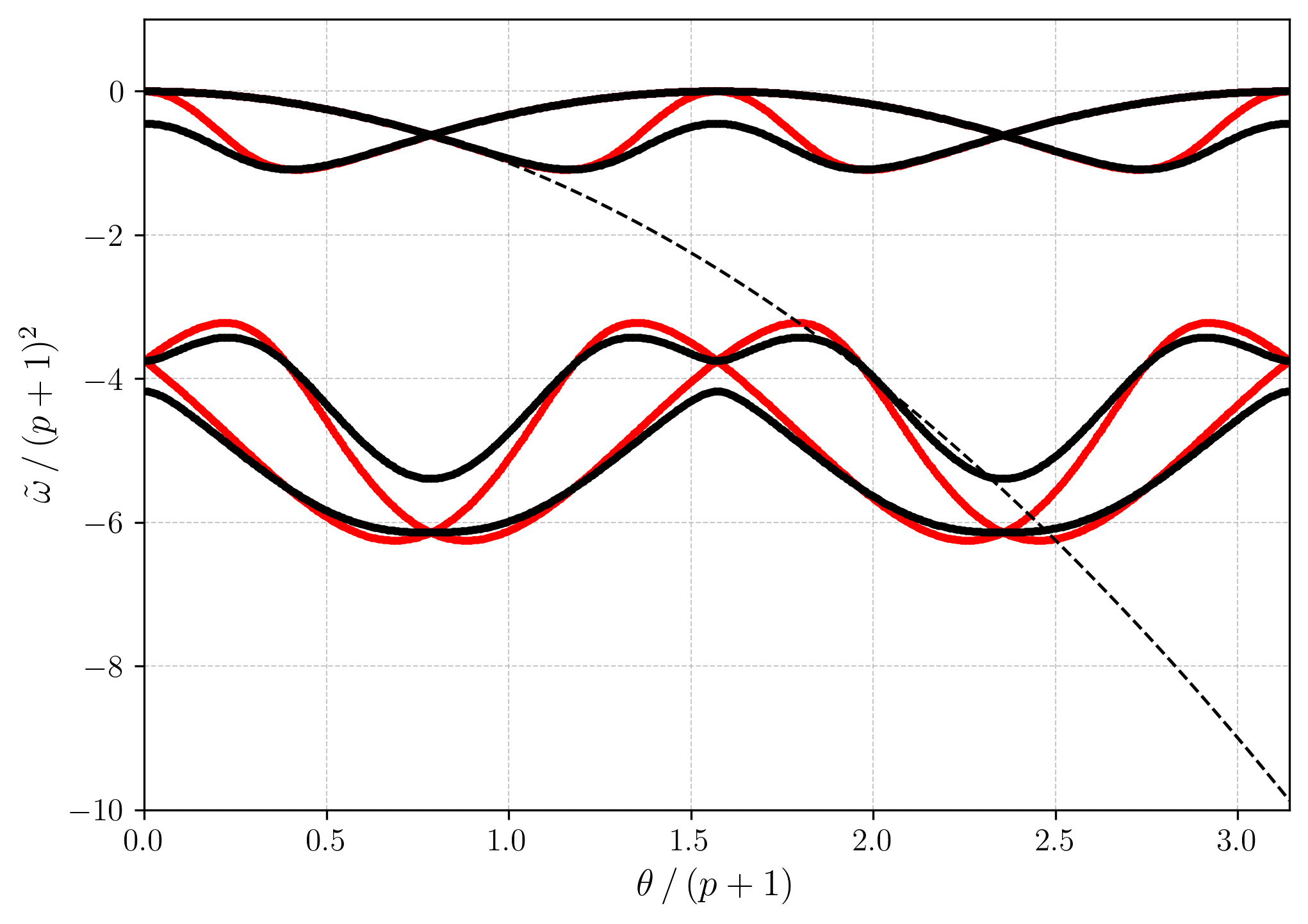}
 \end{subfigure}
 \begin{subfigure}{0.48\textwidth}
     \includegraphics[width=\textwidth]{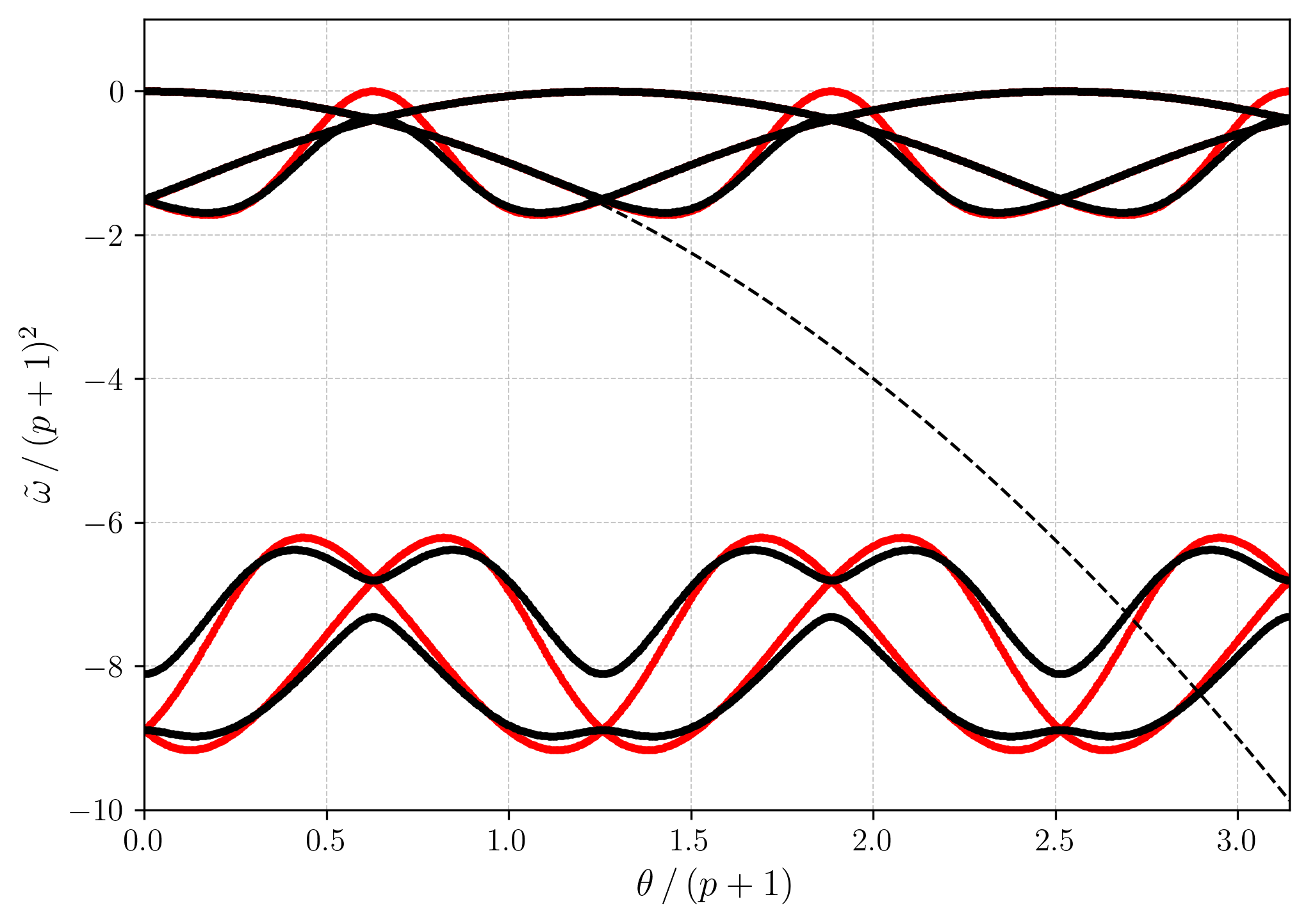}
 \end{subfigure}
\caption{Dissipation curves for the SD method for different order of approximation. From top to bottom, left to right: $p=1,2,3,4$. Black curves denotes the compact stencil formulation and red curves indicates the standard approach. Both schemes are fully centered.}
   \label{fig:eigen_diff}
\end{figure}

We now report results for the dissipation curves corresponding to different polynomial orders of approximation. At this stage, we restrict the analysis to dissipation curves only, since the centered schemes considered here (both extended and compact stencil variants) satisfy the condition $\mathrm{Re}(\widetilde{\omega}) = 0$. In figure~\ref{fig:eigen_diff} we show the dissipation curves for different orders of approximation, comparing both extended and compact stencils. 

We observe that, at least for low orders of approximation, the differences between the two schemes are significant: while the overall shapes remain similar, the standard extended stencil exhibits curves that more frequently and periodically reach the $\widetilde{\omega}=0$ line compared to the compact formulation. This behavior suggests that the compact stencil formulation may be more dissipative than the extended one. The curves presented here follow similar trends with respect to the ones presented by~\citet{huynh:09} under the terminology of centered-$g_{\rm Ga}$/SP-$g_{\rm Lump, Lo}$.

\begin{figure}[h!]
 \centering  
 \begin{subfigure}{0.48\textwidth}
     \includegraphics[width=\textwidth]{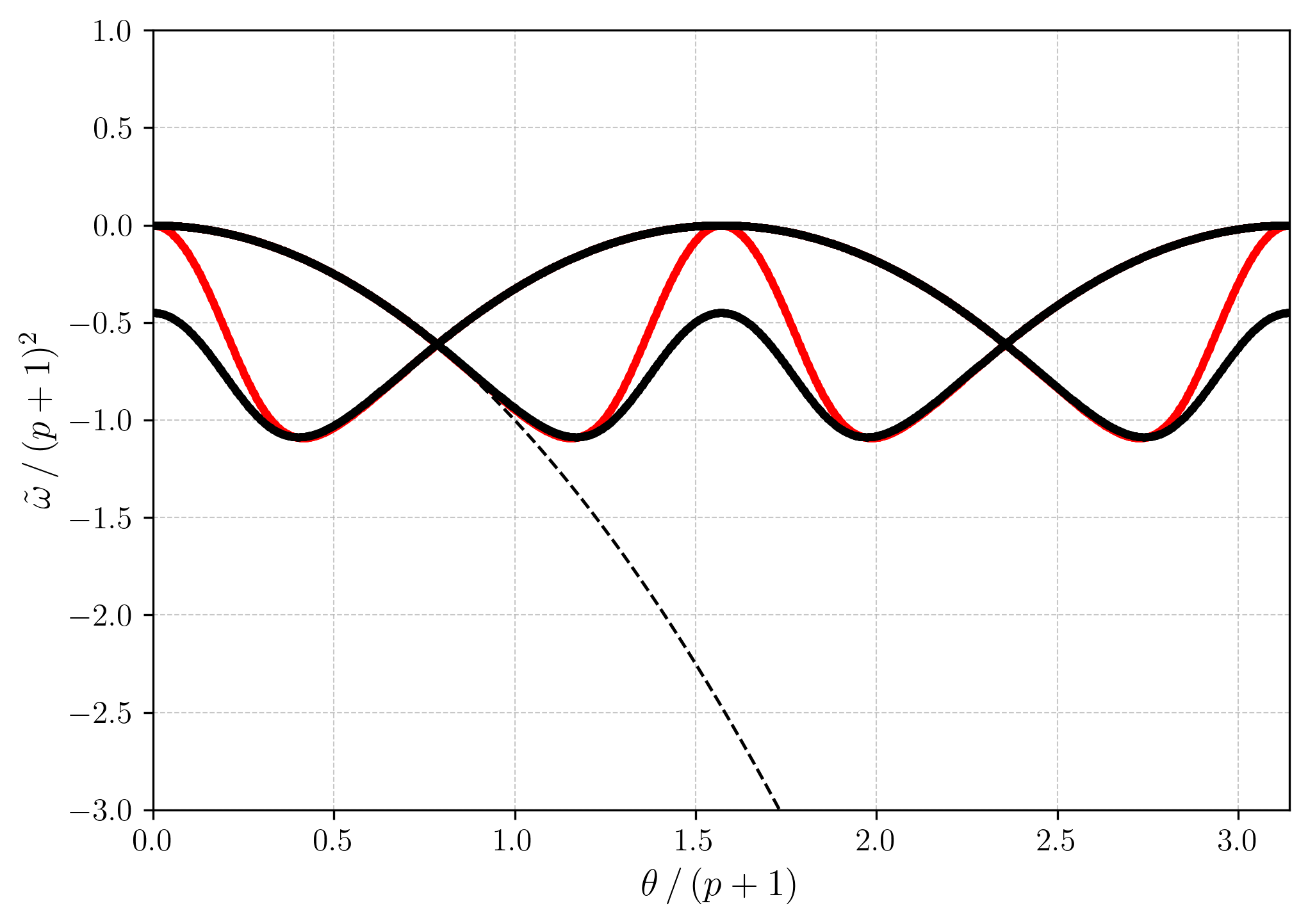}
 \end{subfigure}
 \begin{subfigure}{0.48\textwidth}
     \includegraphics[width=\textwidth]{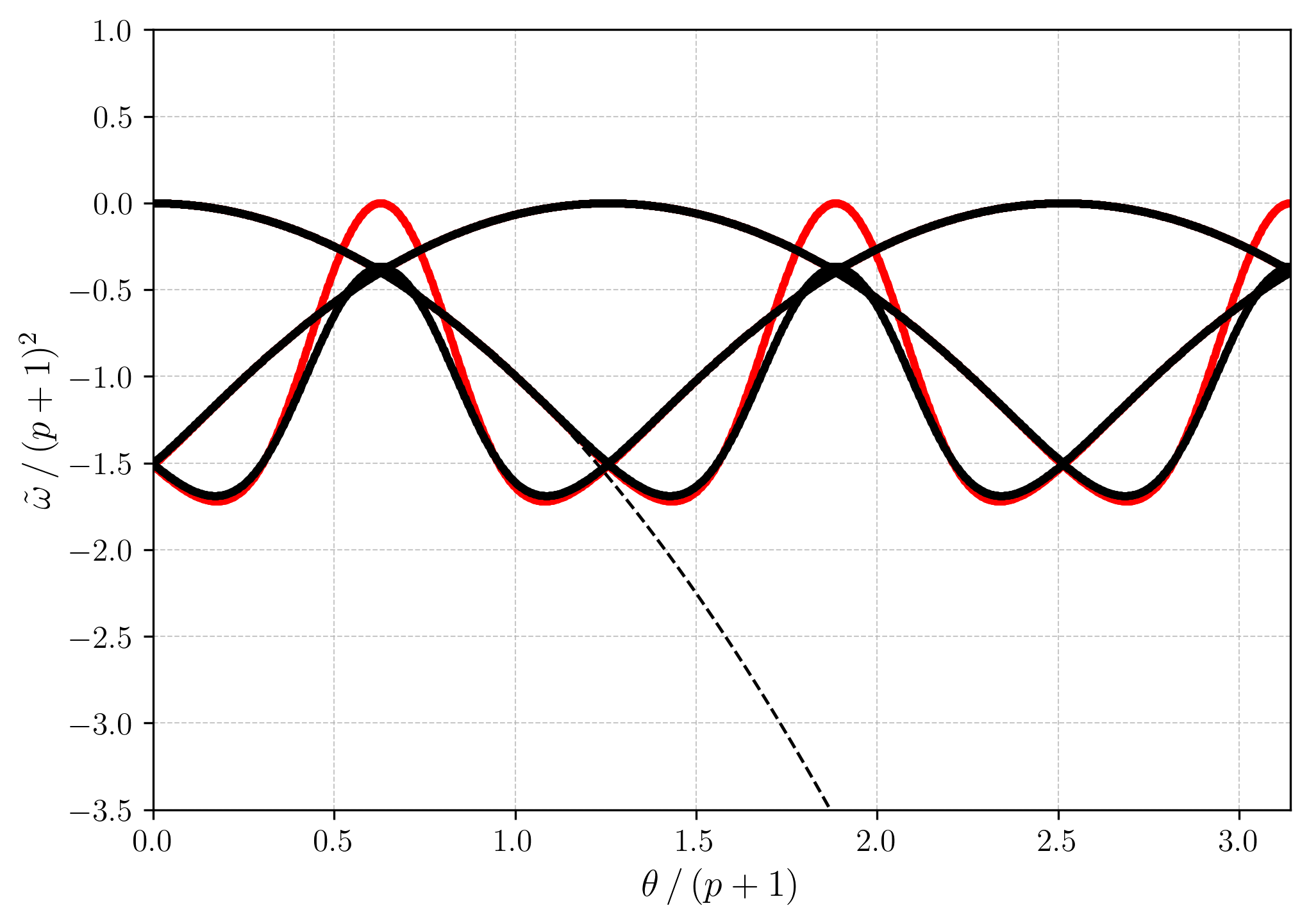}
 \end{subfigure}
 \begin{subfigure}{0.48\textwidth}
     \includegraphics[width=\textwidth]{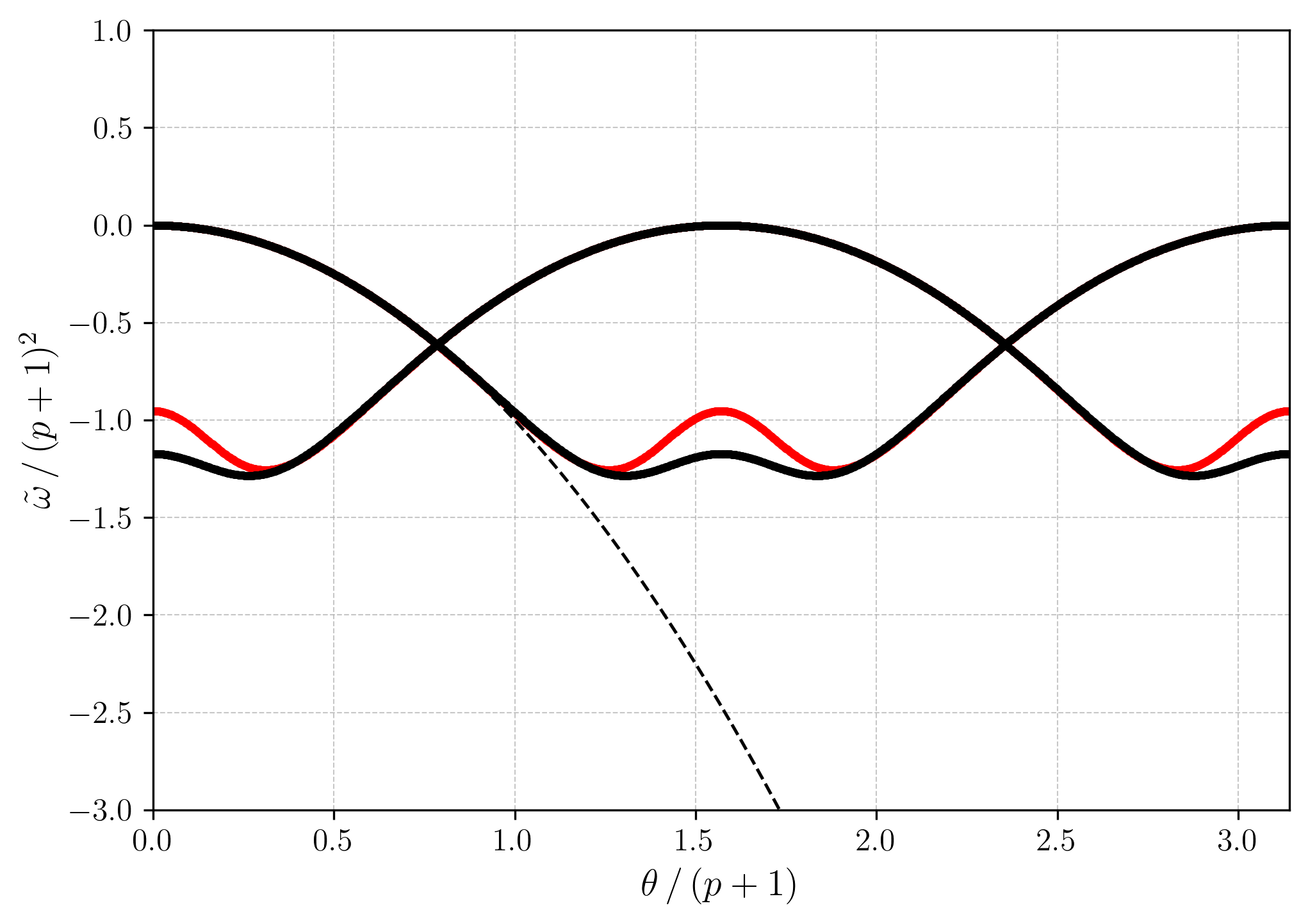}
 \end{subfigure}
 \begin{subfigure}{0.48\textwidth}
     \includegraphics[width=\textwidth]{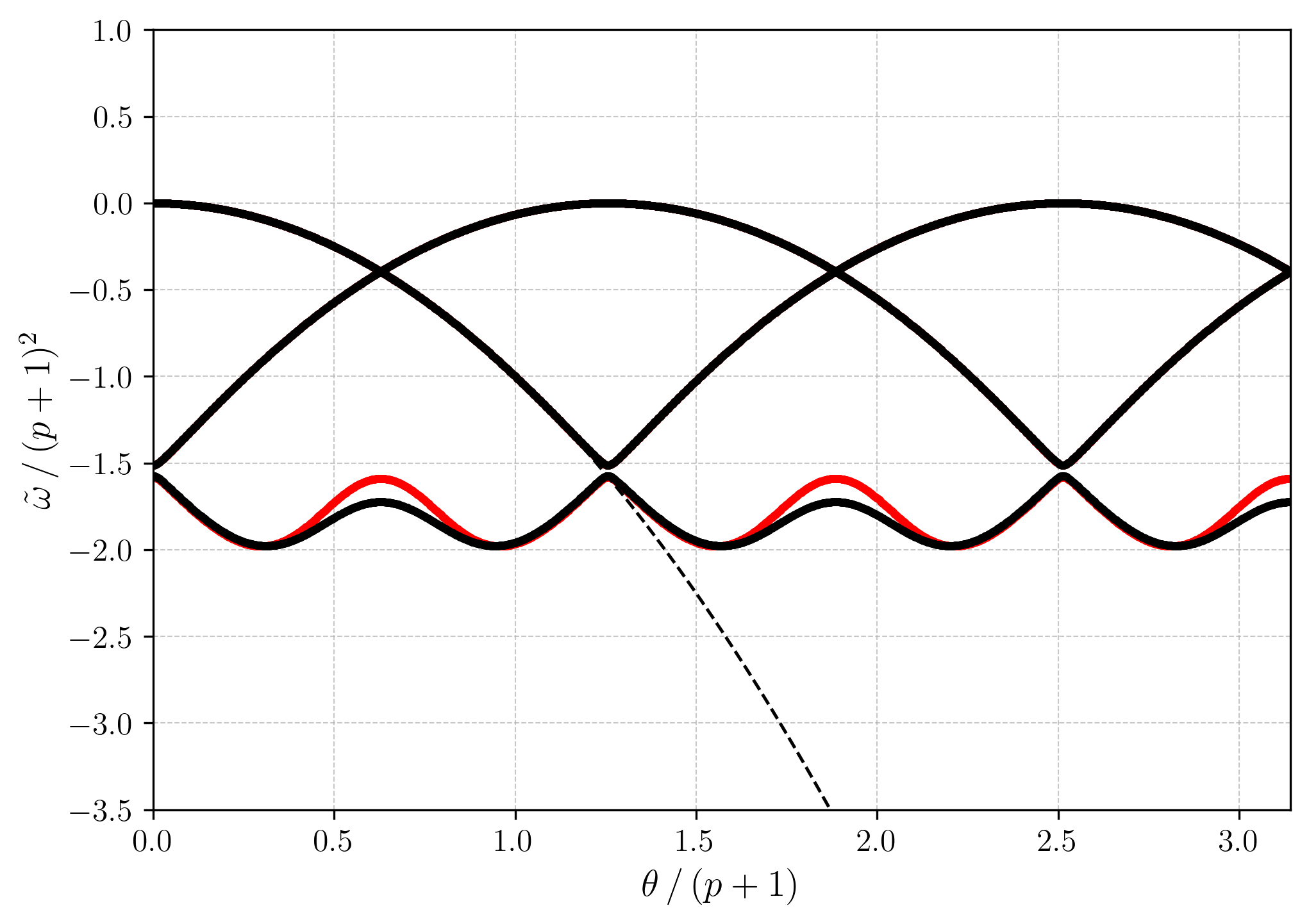}
 \end{subfigure}
  \begin{subfigure}{0.48\textwidth}
     \includegraphics[width=\textwidth]{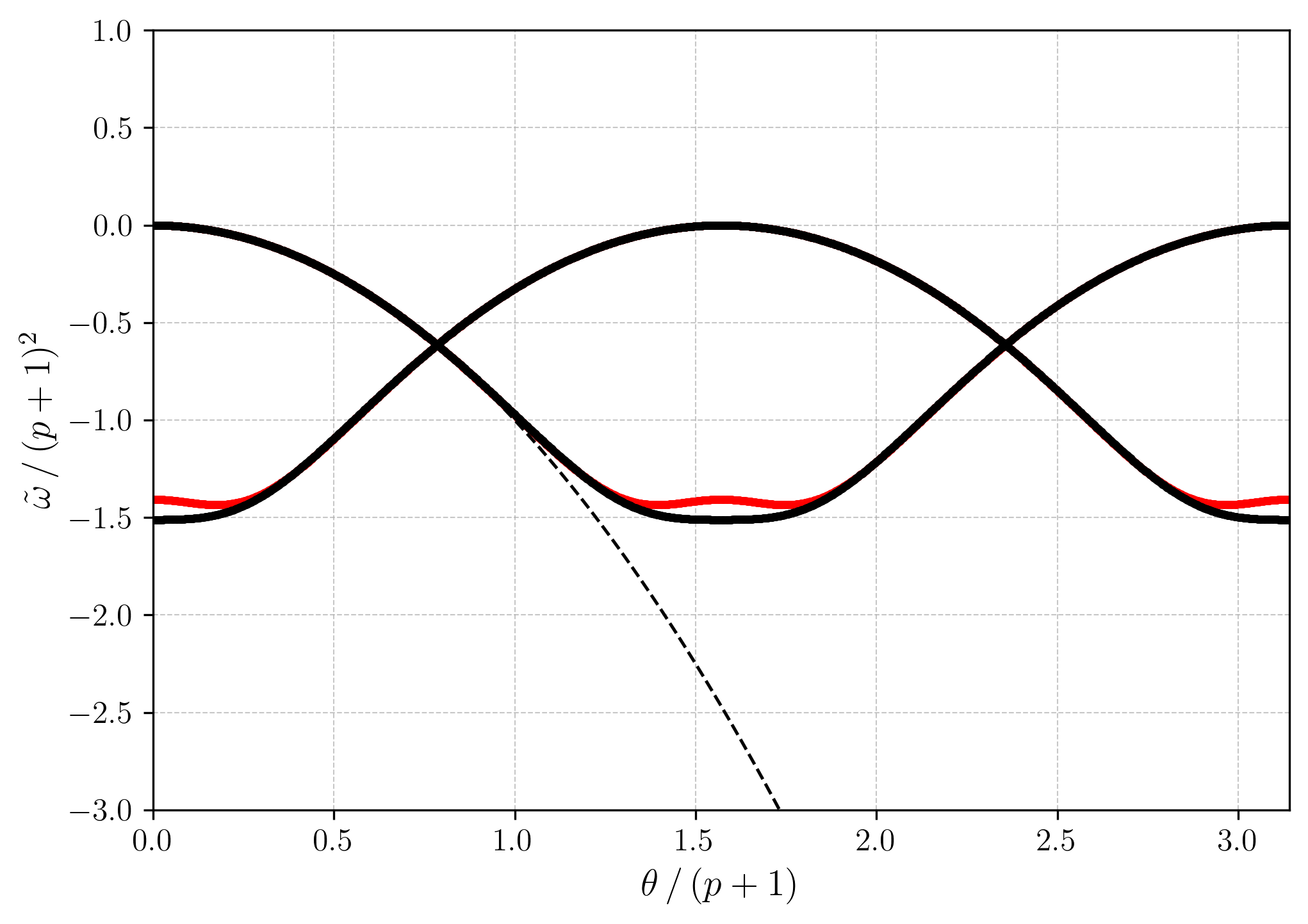}
 \end{subfigure}
 \begin{subfigure}{0.48\textwidth}
     \includegraphics[width=\textwidth]{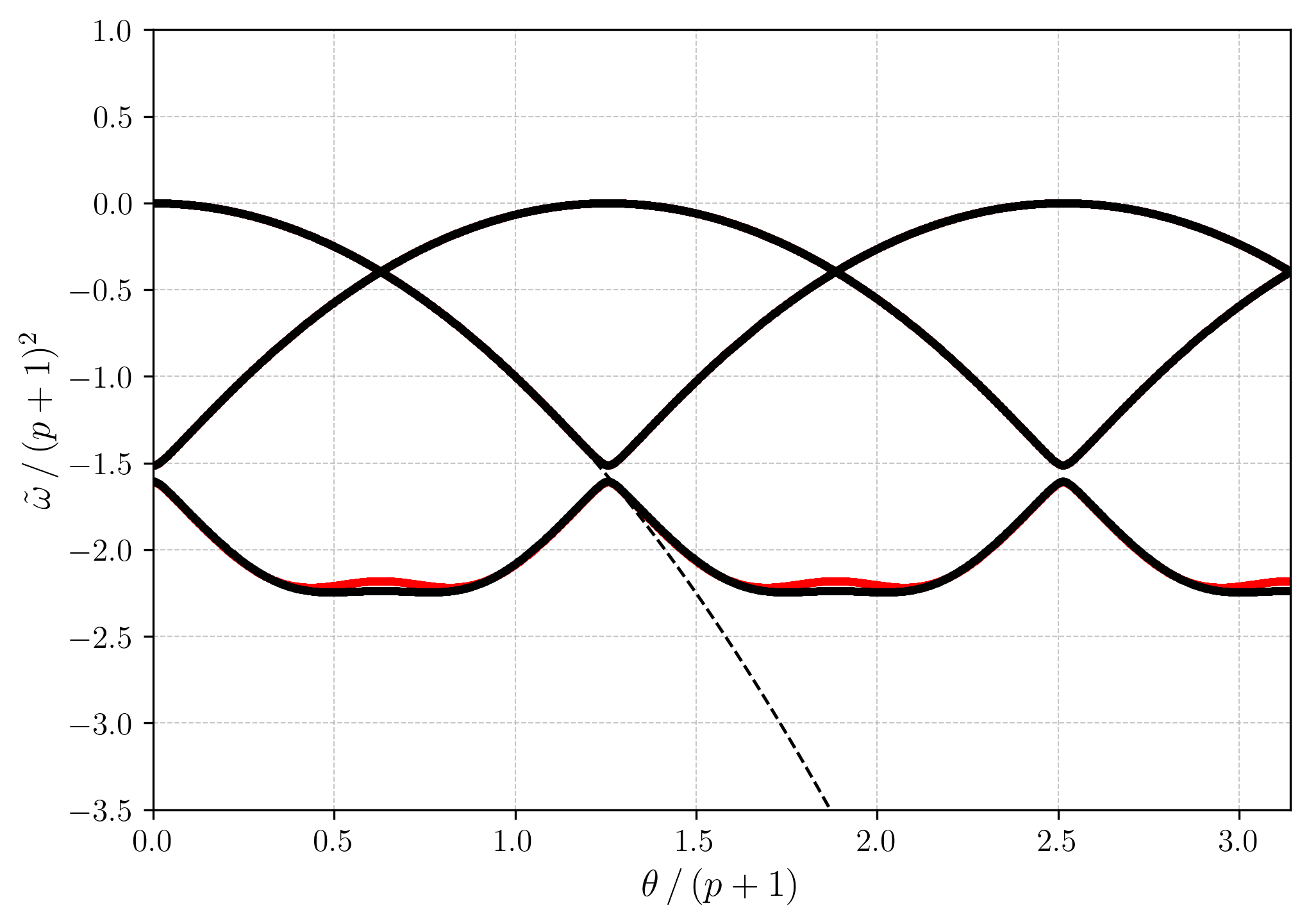}
 \end{subfigure}
  \begin{subfigure}{0.48\textwidth}
     \includegraphics[width=\textwidth]{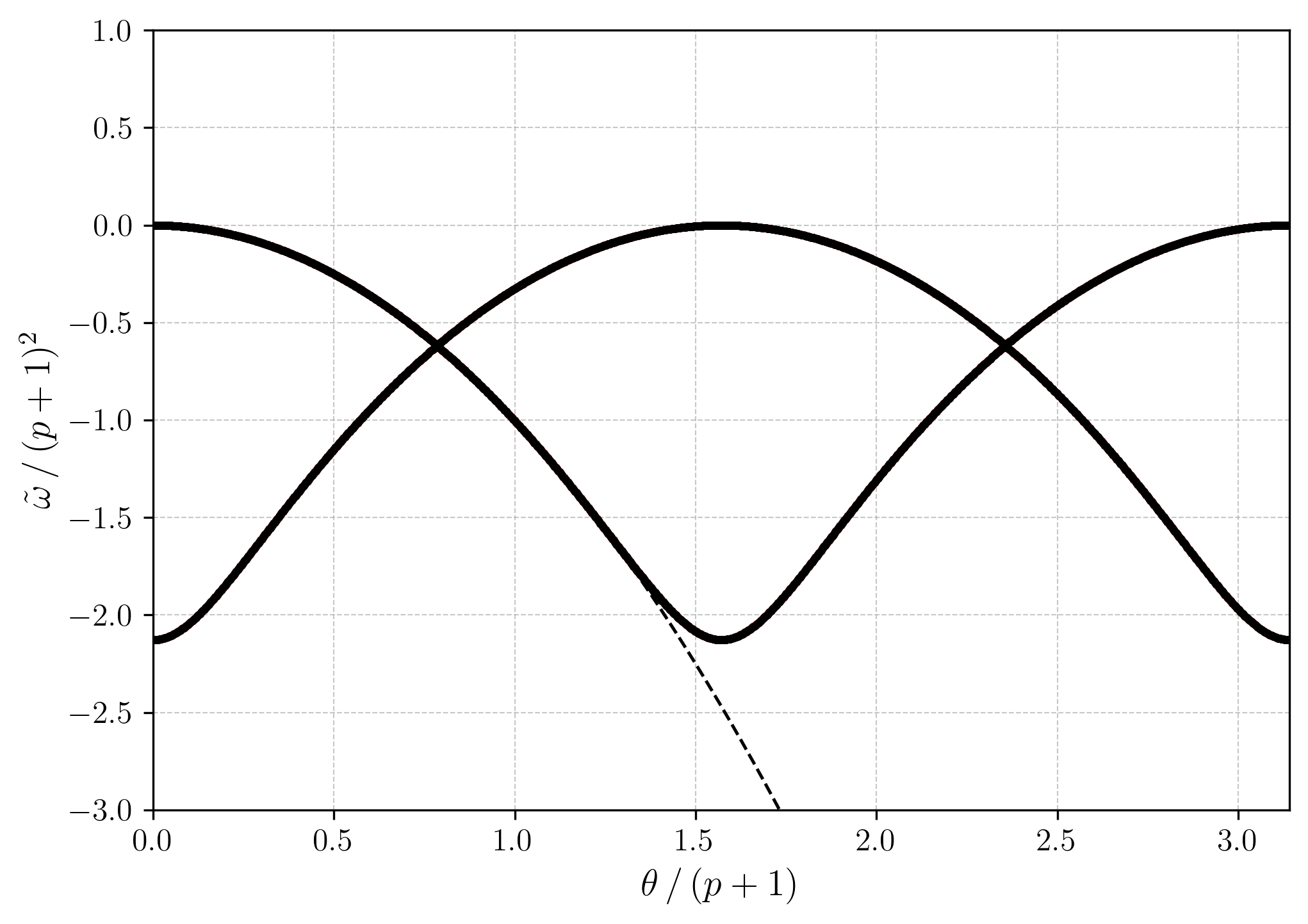}
 \end{subfigure}
 \begin{subfigure}{0.48\textwidth}
     \includegraphics[width=\textwidth]{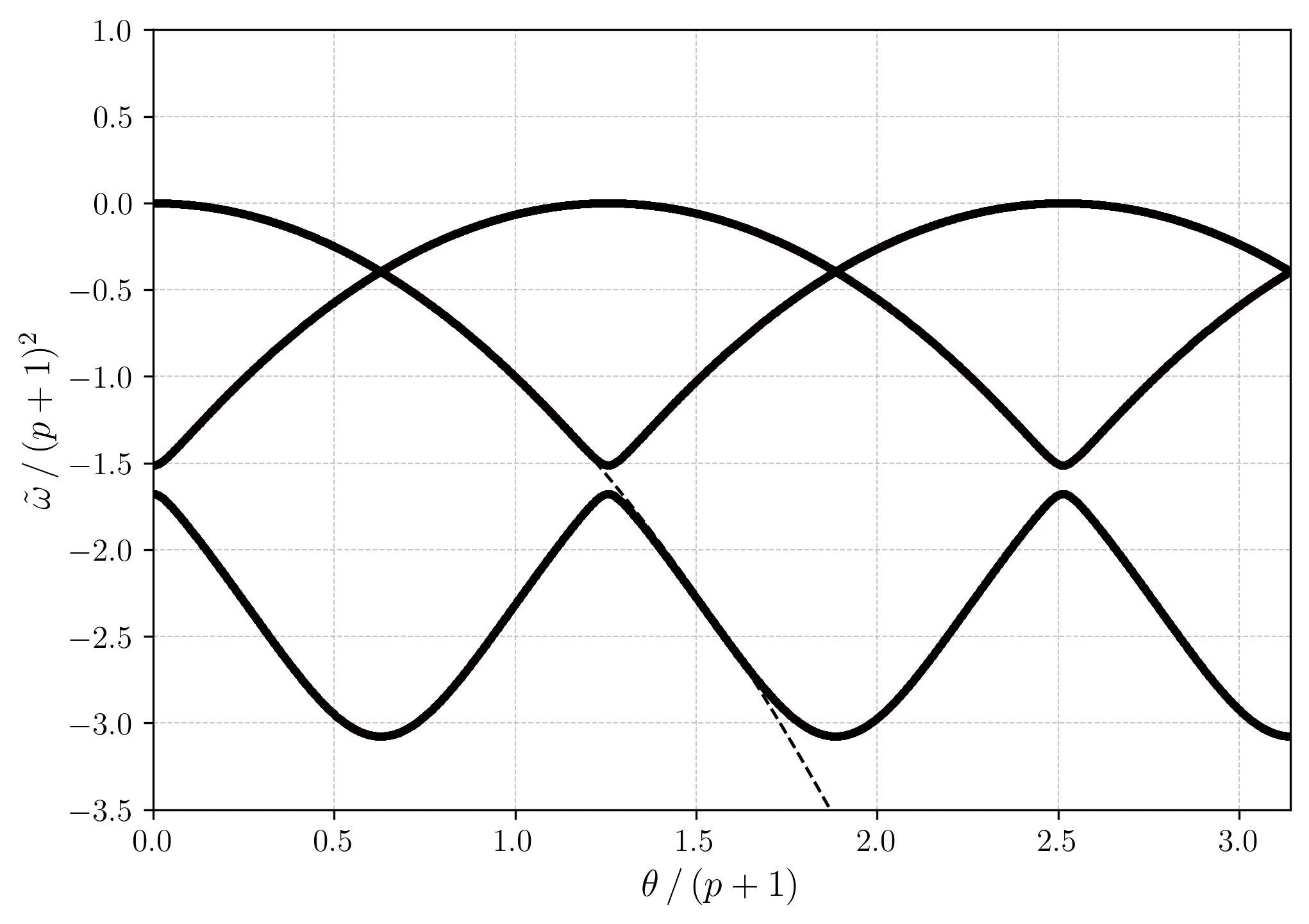}
 \end{subfigure}
\caption{Dissipation curves for the SD method for different $p=3$ (left) and $p=4$ (right) and different values of the interior penalty parameter $\eta_{\rm IP}$. From top to bottom $\eta_{\rm IP}=0.0,0.01,0.02,1.0$. Black curves denotes the compact stencil formulation and red curves indicates the standard approach.}
   \label{fig:eigen_diff_tau}
\end{figure}

In figure~\ref{fig:eigen_diff_tau}, instead, we investigate the influence of the interior penalty parameter on the newly proposed compact scheme for $p=3$ and $p=4$. For the sake of clarity, only the modes located close to the zero-dissipation line are shown. In fact, spurious eigenmodes, typically located far from the exact parabolic diffusion relation (dashed line), become increasingly damped as $\eta_{\rm IP}$ increases. From figure~\ref{fig:eigen_diff_tau}, we observe that increasing the penalty parameter $\eta_{\rm IP}$ progressively enhances numerical dissipation across all resolved modes for both the compact and standard approaches. As the penalty term increases, the differences between the two formulations become progressively less pronounced. This is expected, since for sufficiently large values of $\eta_{\rm IP}$ the penalty terms in equation~\eqref{eq:cmp:int:avg:2d} become dominant with respect to the particular choice of intermediate gradient state at the interface. For relatively small values of the interior penalty parameter, however, the differences between the two approaches remain broadly consistent with the zero-penalty case: the compact formulation provides a slightly larger amount of numerical dissipation at intermediate wavenumbers. Notice that, as commonly known in the literature, the interior penalty approach is certainly one viable option to increase numerical dissipation and consequently suppress numerical oscillations. However, this comes at the price of more restrictive CFL conditions. 
\section{Numerical results}\label{sec:3}
\subsection{Convergence study: diffusion of a well-resolved gaussian profile}
As a starting analysis, we consider the diffusion of a sufficiently well-resolved Gaussian to assess the convergence properties of the proposed formulation. 
In particular, we consider an initial Gaussian profile with a typical width $\Delta$ equal to 10\% of the domain length. In terms of a classical Gaussian filter, this relates to the time of pure diffusion of an initial Dirac, $t_i$, and to the corresponding initial variance, $\sigma^2$, via the well-known relation $D t_i = \sigma^2 / 2 = \Delta^2 / 24$, where $D$ is the diffusion coefficient. Hence, the initial profile is:
\begin{equation}
    u_{0}(x) = \frac{1}{\sqrt{4\pi D t_i}}\exp \left[ -\frac{(x-1/2)^2}{4 D t_i} \right] \quad \mathrm{with} \quad D t_i = \frac{\Delta^2}{24} = \frac{1}{2400}, \quad\text{and}\quad x\in[0,1].
\end{equation}

Notice that, to evaluate the error, the analytical solution after a total simulated time $T$ is easily obtained as
\begin{equation}
  u_{\rm a}(x) = \frac{1}{\sqrt{4\pi D (t_i + T)}}\exp \left[ -\frac{(x-1/2)^2}{4 D (t_i + T)} \right].
\end{equation}

The number of elements and the polynomial order are varied in order to perform a proper convergence study. In particular, four refinement levels are considered, based on the total number of degrees of freedom. Namely, we use $300$, $600$, $1200$, and $2400$ total degrees of freedom. For each refinement level, the combination of polynomial order and number of elements is selected such that the resulting discretization matches the prescribed total number of degrees of freedom. The specific resolutions considered for this study are reported in table~\ref{tab:resolutions}.
\begin{table}[!h]
\begin{center}
\begin{tabular}{|l | c | c | c | c | c|} 
 \hline
  \# DoF & $p=1$ & $p=2$ & $p=3$ & $p=4$ & $p=5$\\ 
 \hline
 300 &  $150$ & $100$ &  $75$ &  $60$ &  $50$ \\ 
 \hline
 600 &  $300$ & $200$ & $150$ & $120$ & $100$ \\ 
 \hline
1200 &  $600$ & $400$ & $300$ & $240$ & $200$ \\ 
 \hline
2400 & $1200$ & $800$ & $600$ & $480$ & $400$ \\ 
 \hline
\end{tabular}
\end{center}
\vspace{-0.4cm}
\caption{Detailed view of the total number of elements considered in the convergence study in order to match the total number of degrees of freedom (first column) for each polynomial order (first row).}
\label{tab:resolutions}
\end{table}

\begin{figure}
\centering
\includegraphics[width=0.75\textwidth]{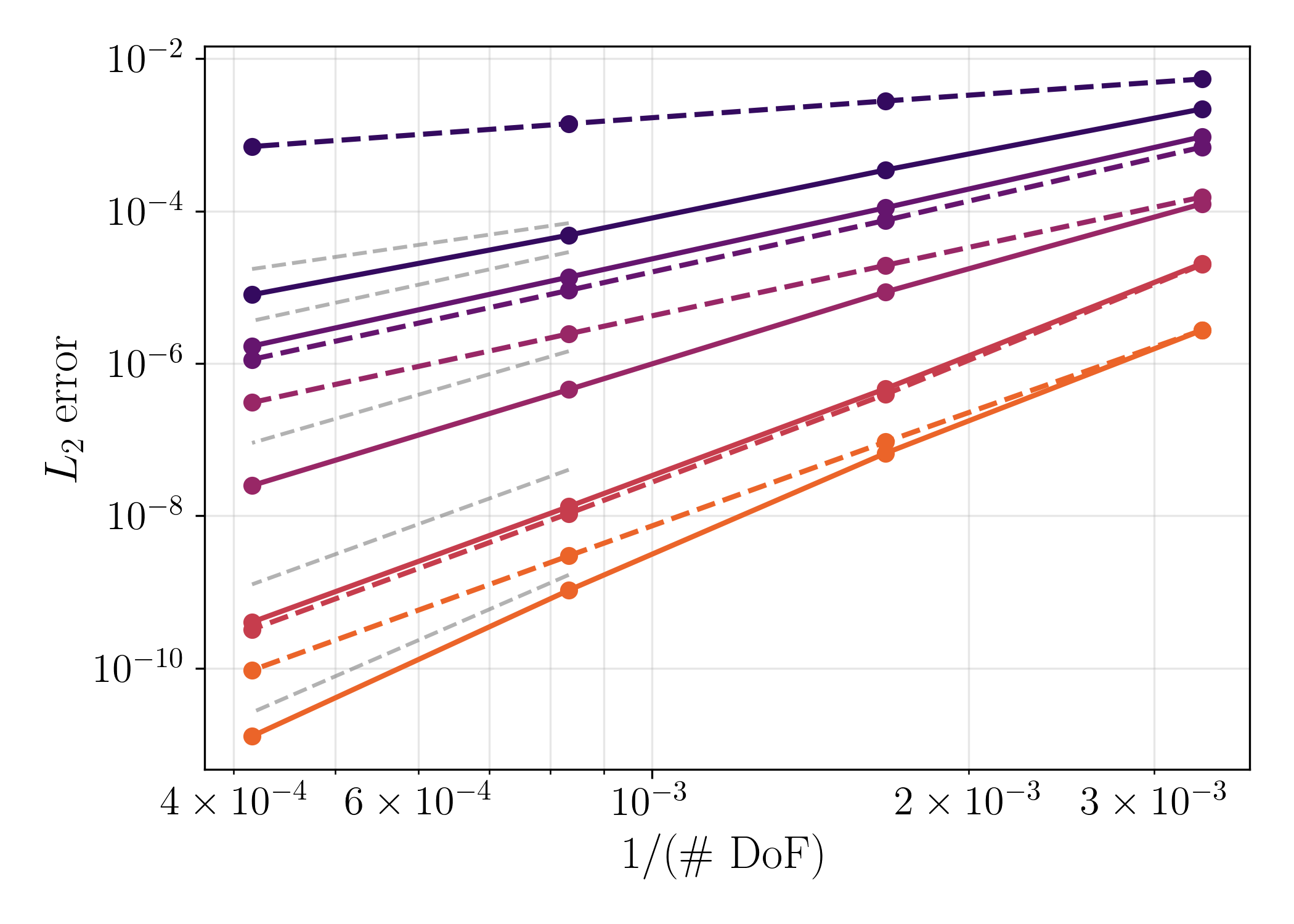}
\caption{$L_{2}$-error against the inverse of the total number of DoF for different orders of approximation. Color-gradient indicates different polynomial orders: from dark purple to orange, $p=1,2,3,4,5$. Solid lines, compact formulation; dashed lines, standard formulation. Dashed gray lines indicate the expected convergence rate.}
\label{fig:gauss_convergence}
\end{figure}

We evolve the initial condition for $T=2\E{-5}$ and compare the numerical solutions for different polynomial orders and levels of refinement. The $L_{2}$-errors for both compact and standard formulation are shown in figure~\ref{fig:gauss_convergence}. We can observe that, for the compact formulation, for every polynomial order, the expected convergence rate is recovered. For the standard approach, on the other hand, different behaviors are observed for even (\ie, $p=1,\,3,\,5$) and odd (\ie, $p=2,\,4$) orders.
In particular, if at odd orders the standard scheme recovers the expected convergence rate and performs similarly to the compact formulations, at even orders of approximation, it provides a much more degraded convergence rate, which seems to be at least one order slower than the expected convergence rate. This behavior has previously been observed in the literature within the framework of DG schemes~\cite{bassi1997high2,cockburn1998local}.
As a more quantitative assessment, table~\ref{tab:convergence} reports the convergence rates obtained for the two schemes at the different orders of approximation.
\begin{table}[!ht]
\begin{center}
\begin{tabular}{|l | c | c | c | c | c|} 
 \hline
   & $p=1$ & $p=2$ & $p=3$ & $p=4$ & $p=5$\\ 
 \hline
compact & 2.599 & 3.006 & 4.186  & 5.035 & 6.367\\ 
 \hline
standard & 0.998 & 3.001 & 2.996  & 5.050 & 4.990\\ 
 \hline
nominal & 2 & 3 & 4  & 5 & 6\\
 \hline
\end{tabular}
\end{center}
\vspace{-0.4cm}
\caption{Convergence rates for different orders of approximation for the compact and standard approaches.}
\label{tab:convergence}
\end{table}
The results confirm the behavior observed in figure~\ref{fig:gauss_convergence}. For even orders of approximation, the standard approach exhibits a reduced convergence rate, approximately one order lower than the theoretical expectation, while for odd orders the expected convergence rate is recovered. Conversely, the compact formulation consistently achieves the expected order of accuracy for all polynomial degrees, with convergence rates even slightly higher than the theoretical values in some cases.

Notice that previous studies on numerical fluxes for advection-diffusion equations using discontinuous spectral element methods~\citep{watkins2016numerical} reported that the so-called \emph{centered-centered} scheme provided the best performance for well-resolved profiles while under-performing for under-resolved cases. However, the analysis carried out by~\cite{watkins2016numerical} was limited to third-order approximations. The results reported in figure~\ref{fig:gauss_convergence} are consistent with these findings: the standard fully-centered approach slightly outperforms the compact formulation only for $p=2$, and the difference between the two methods is marginal. For even orders of approximation the performance gap becomes significantly larger in favor of the compact formulation. Also, the work by~\cite{ortleb2020comparative} highlighted the significant differences between odd and even polynomial orders in the context of a nodal DG discretizations. In particular, different formulations of the LDG flux exhibited different performance depending on the polynomial order of the underlying approximation.
\subsection{Diffusion of localized Dirac's delta}
While convergence studies based on smooth solutions are traditionally employed to assess the formal accuracy of a numerical scheme, its robustness must also be evaluated in strongly under-resolved regimes. These scenarios are particularly relevant in practical applications, where steep gradients or localized features may develop and become challenging to capture on coarse discretizations. 
\begin{figure}
\centering
\includegraphics[width=\textwidth]{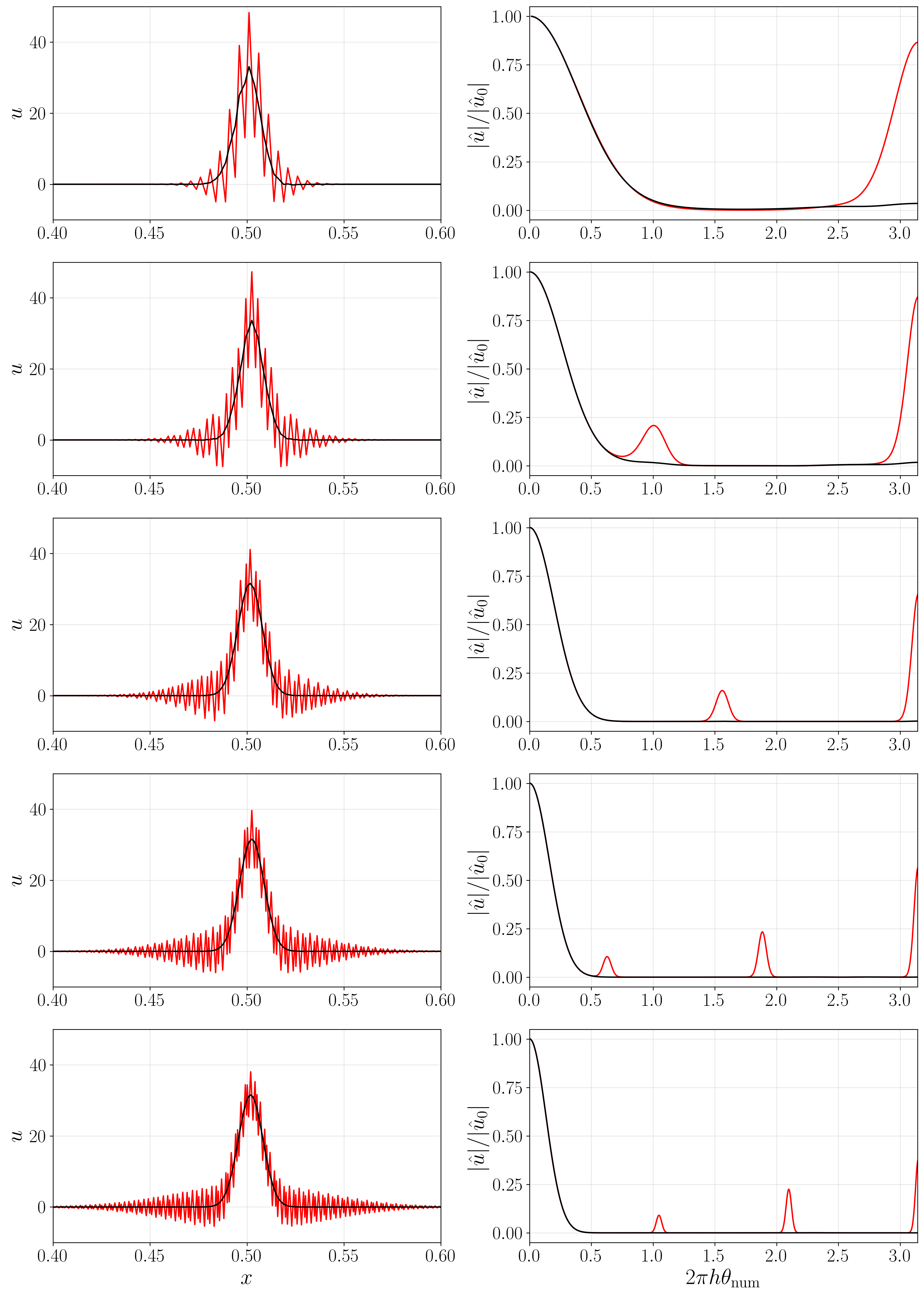}
\caption{Numerical solution in physical (left) and Fourier (right) spaces for a mesh composed of $200$ elements. From top to bottom $p=1,2,3,4,5$. Black curves denotes the compact stencil formulation and red curves indicates the standard approach.}
\label{fig:dirac_all}
\end{figure}

In order to reproduce such scenarios, we consider the diffusion of a localized Dirac's delta function. The integration time is the same considered in the previous test case. In particular, regarding the initial condition, we set the solution to zero at all the solution points except one, where a constant value is imposed as 
\begin{equation}
    u(x_{d}) = \frac{1}{w (x_{d})\Delta x}
\end{equation}
with $\Delta x$ the size of the one-dimensional element, $x_{d}$ the closest point to the middle of the domain and $w(x_{d})$ its corresponding quadrature weight. This formulation is so that, by varying the polynomial order or the resolution, the Gauss-quadrature  integral of the Dirac's delta is preserved and it is equal to unity.

In figure~\ref{fig:dirac_all} we compare the numerical solution for the diffusion of a localized Dirac's delta using the two fully centered formulations presented in the previous sections for different orders of approximation. First of all, by inspecting the solution in physical space only, we can notice more marked oscillations by employing the standard approach to compute viscous fluxes on the extended stencil. The compact stencil, instead, significantly reduces these oscillations.

To assess the robustness and the accuracy of the numerical method, we also consider the behavior of the numerical solution in the Fourier space. Given a numerical solution sampled at $N$ equally spaced grid points $x_j$, the discrete Fourier coefficients are computed as:
\begin{equation}
    \widehat{u}(k) =\frac{1}{2\pi}\sum_{j=1}^{N} u(x_j)\cos(2\pi k x_j),\label{eq:fft}
\end{equation}
where $k$ indicates the wavenumber.
The importance of the Fourier-space analysis becomes evident in the present case. In fact, considering the red curves indicating the standard approach, we can notice periodically excited frequencies which are instead normally suppressed by the compact formulation. By comparing figures~\ref{fig:eigen_diff} and~\ref{fig:dirac_all}, we can notice that the number of peaks and their location coincide exactly with those predicted by the eigenanalysis. Notice that similar trends in Fourier space (\ie, with localized excited frequencies) have been also observed in some recent work within the framework of DG schemes employing BR1 fluxes \cite{alhawwary2019study,alhawwary2020combined}.
\subsection{Non-linear diffusion of localized Dirac's delta}
In this section we present the one-dimensional solutions of the model equation:
\begin{equation}
\pd{u}{t} = \lapl{u^m} = \divg{(m\,u^{m-1} \grad{u})},\quad\text{for}\quad m \ge 0,
\label{eq:pme}
\end{equation}
also known as the porous medium equation (PME).
This non-linear model equation, which reduces to the heat equation when $m = 1$, can be used to describe mass diffusion in porous media (as the name implies) but can be also used in other areas of fluid mechanics, such as lubrication theory~\citep{vazquez:06,caballero:26}.

Although possibly more prone to introduce aliasing errors, the current implementation of the right-hand side of equation~\eqref{eq:pme} consists in simply computing the non-linear term, $m\,u^{m-1}$, from the interpolated solution at flux points.
Other, more elegant, implementations are possible. For instance, computing the non-linear term at solution points and then interpolate it at flux points such as to better constrain it into the correct polynomial space.
This last implementation was also tested and, although better behaved overall, did not suggest different conclusions on the comparison between the standard and the compact schemes. The relevant results, for the sake of brevity, are hence omitted. 

The solution of the PME to an initial Dirac located at the origin of the reference frame is known as the Zeldovich, Kompaneetsm and Barenblatt (ZKB) solution~\citep{zeldovich:50,barenblatt:52,barenblatt:96}, which takes the form:
\begin{equation}
u(\vect{x},t) = t^{-\alpha}\bigg[ C - \frac{\alpha (m-1)}{ 2 m d } \norm{\vect{x}}^2 t^{-2\beta} \bigg]_{+}^{\frac{1}{m-1}},
\label{eq:zkb}
\end{equation}
with $\alpha = d/[d (m-1) +2]$, $\beta = \alpha / d$, $d$ the dimensionality of the problem and the operator $[\;\cdot\;]_+ = \max \{\;\cdot\;, 0\}$.
In the case of applications to mass diffusion in porous media, the constant $C > 0$ is chosen to satisfy the mass conservation. 
In our one-dimensional (\ie, $d = 1$ and $r = \abs{x}$) tests, it is arbitrarily set equal to $0.1$. Concerning the main parameter $m$, numerical tests are conducted for $m = 5$.

The computation is initialized with the analytical solution at an initial time $t_0$, which is computed such as to have the profile within a single element. Accordingly, if $\Delta x$ is the element's width, the initial time, in the one-dimensional case, is obtained from equation~\eqref{eq:zkb} after imposing $u = 0$ and $\abs{x} = \Delta x /2$:
\begin{equation}
t_0 = \bigg[ C \frac{ 2 m }{\alpha_1 (m-1)} \bigg( \frac{2}{\Delta x} \bigg)^2 \bigg]^{-\frac{1}{2\alpha_1}},\quad\text{with}\quad \alpha_1 = \frac{1}{ (m-1) +2 }.
\end{equation}

\begin{figure}
\centering
\begin{tikzpicture}
\node (img) {
    \includegraphics[width=\textwidth]{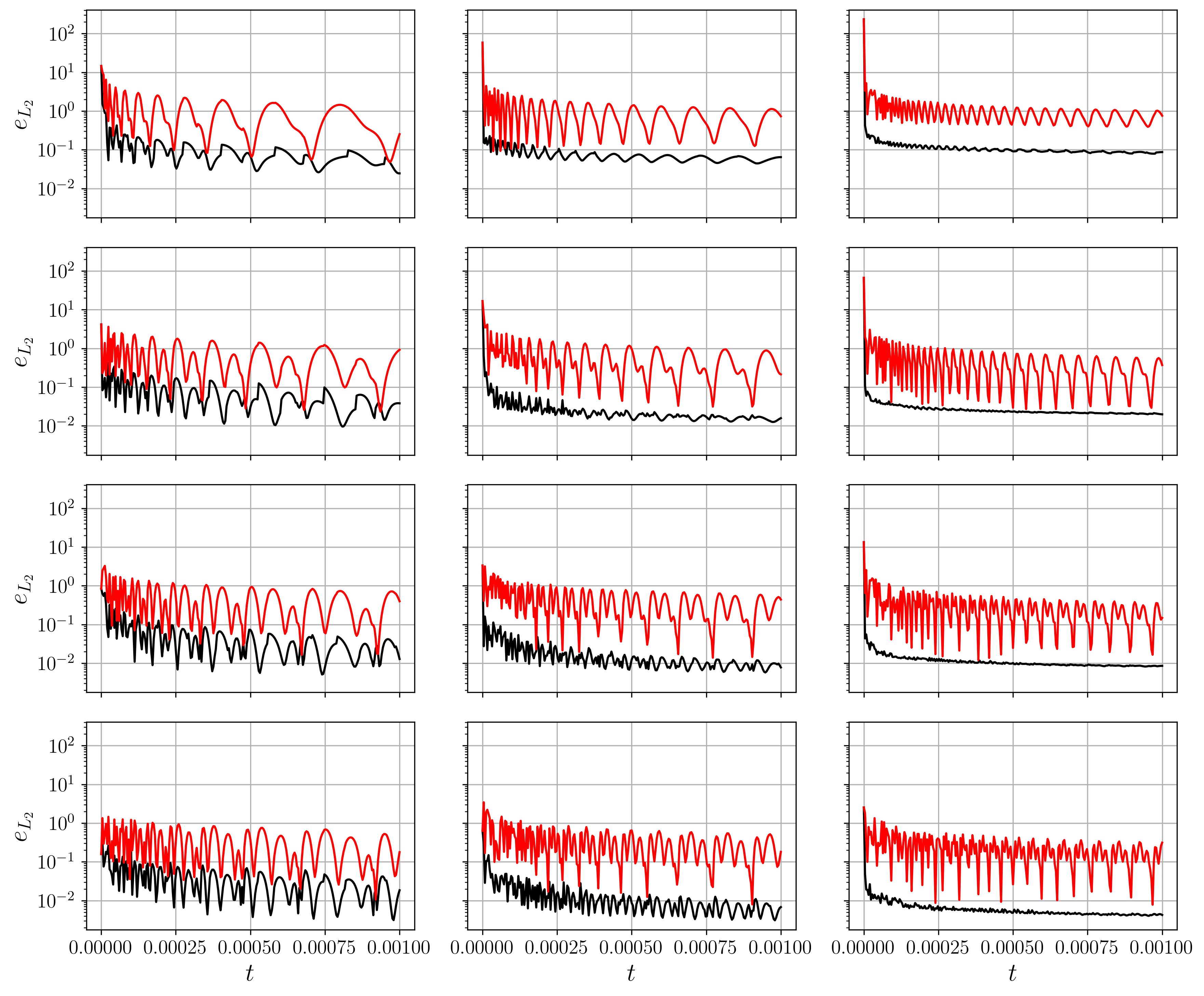}
};



\end{tikzpicture}
\caption{Evolution of the $L_{2}$ error as a function of time for different resolutions and different polynomial orders for the two schemes. From left to right the mesh is made of $n_{\rm el}=50,100,200$ equi-spaced elements. From top to bottom $p=1,2,3,4$. Black curves denotes the compact stencil formulation and red curves indicates the standard approach.}
\label{fig:nl_diff_err}
\end{figure}

In figure~\ref{fig:nl_diff_err} we show the time evolution of the $L_{2}$-error of the numerical simulations with respect to the exact solution using standard and compact formulations for different polynomial order and mesh resolutions.\footnote{
Notice that, as the one-dimensional solution has compact support over the interval 
$$\abs{x}^2 \le H(t),\quad\text{with}\quad H(t) = \frac{2 m C}{\alpha_1 (m-1)} t^{2\beta},$$
to prevent the $L_{2}$-error from being dominated by the error at the discontinuities in $x = \pm \sqrt H$, we restrict its evaluation within the interval $\abs{x}^2 < 0.9 H$, thus accounting for about 95\% of the bell-shaped solution at each time step.
}
From top to bottom we increase the polynomial order while keeping fixed the number of elements while, vice versa, from left to right, we increase the number of elements while keeping constant the polynomial order. First of all, we can notice that the error significantly oscillates in time for both the compact and standard approaches. However, we can also observe that the error using the compact formulation is always significantly smaller than the one employing the standard formulation, even more so for the most refined case. For all polynomial orders the compact formulation provides significantly smaller errors with respect to the standard one. 

\begin{figure}
\centering
\includegraphics[width=\textwidth]{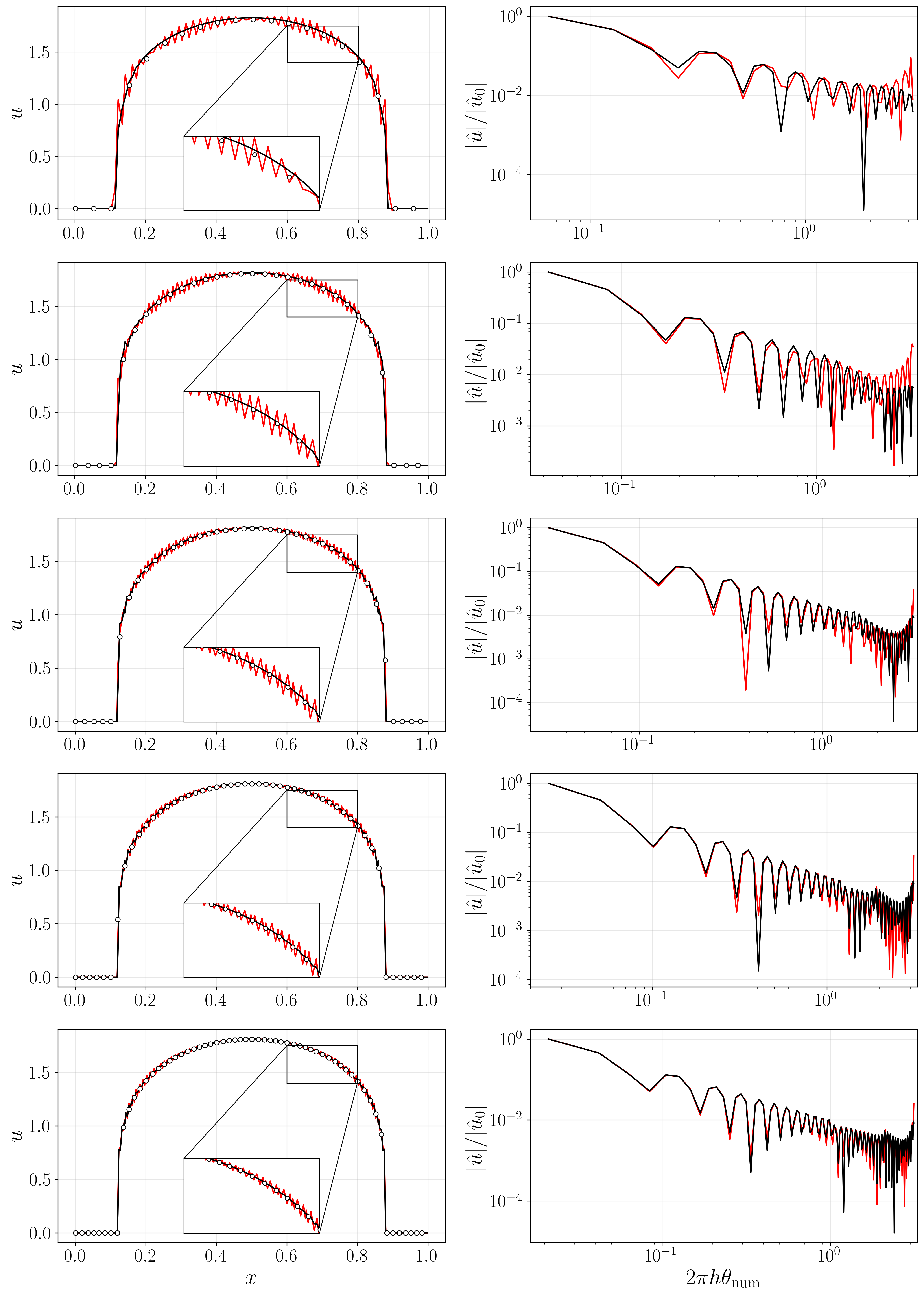}
\caption{Numerical solution in physical (left) and Fourier (right) spaces for a mesh composed of $50$ elements at $t= 9\E{-4}$. From top to bottom row $p=1,2,3,4,5$. Black curves denote the compact stencil formulation, red curves indicate the standard approach and white dots represent the exact solution. Both schemes are fully centered.}
\label{fig:nl_diff_sol_all}
\end{figure}

In figure~\ref{fig:nl_diff_sol_all} we show the numerical solution for $p=1,2,3,4,5$ at $t=9\E{-4}$. We can observe that at this time instance, the standard approach is characterized by visible oscillations around the exact solution. The compact approach, instead, suppresses in a significantly more effective way such numerical artifacts.
It is worthwhile mentioning that, due to the modulation observed in the $L_2$-error, the level of oscillations varies in time and the \emph{worst scenario} might take place at different times for different orders (\eg, the case $p = 5$ with the standard flux would be significantly more oscillatory at $t = 1\E{-3}$). For the sake of simplicity, we fix a common time of observation for all the tested orders. 

Focusing on the right panels of figure~\ref{fig:nl_diff_sol_all}, where the solution is plotted in Fourier space, we can notice a significant accumulation of energy at the high wavenumbers. Being the problem nonlinear, the interaction between each Fourier mode cannot evolve independently as shown in the eigenanalysis and in the pure diffusion equation~\cite{moura2015linear}. While the characteristic numerical dissipation of the two schemes affects different frequencies with different intensity, the nonlinear dynamics of the problem promotes energy transfer in Fourier space across different scales. In particular, a significant amount of energy is transferred at the high frequencies where it cannot be properly dissipated by either physical or numerical diffusion. This trend can be observed in both approaches but it is certainly more marked for the standard formulation.

\begin{figure}
\centering
\includegraphics[width=\textwidth]{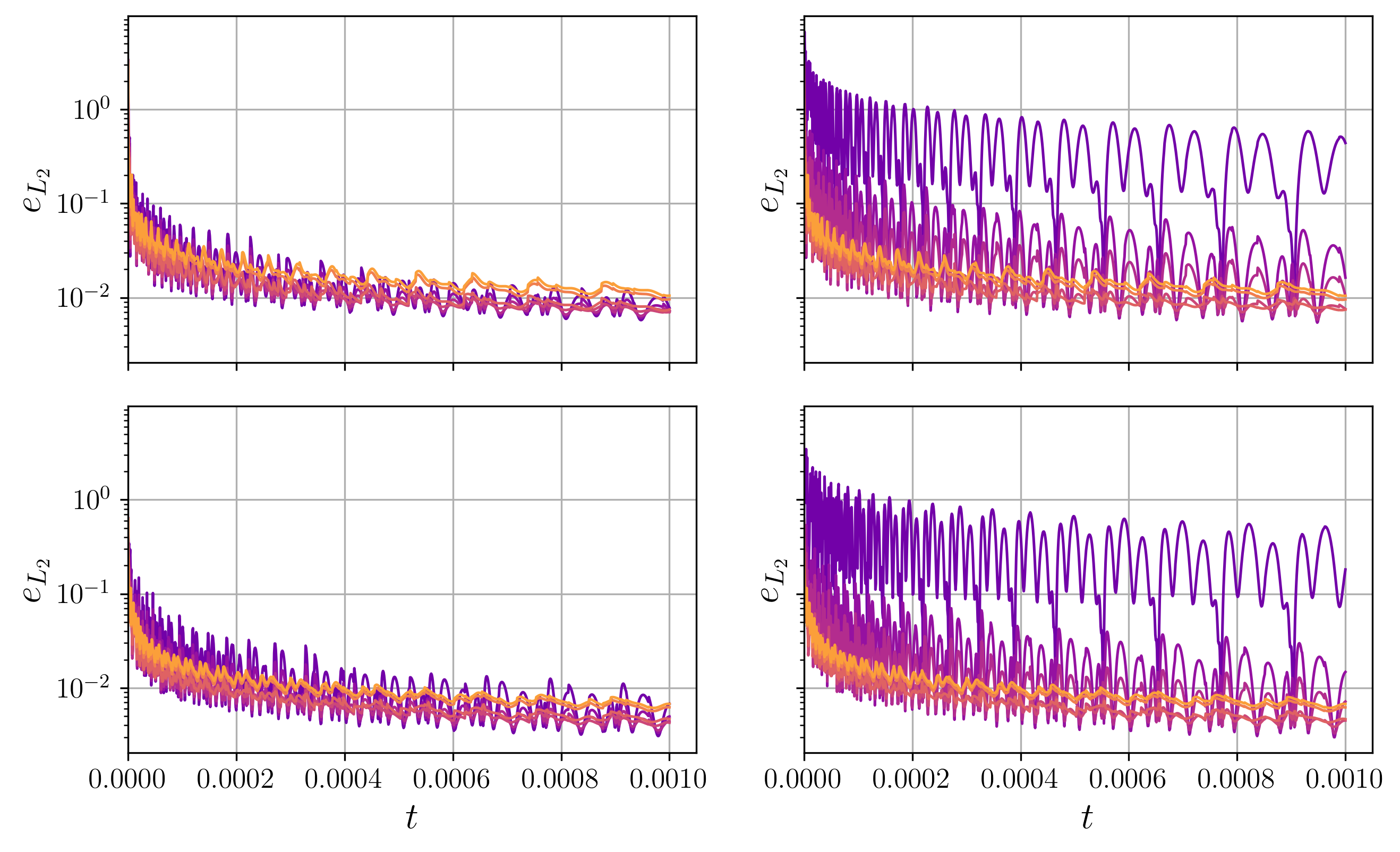}
\caption{Evolution of the $L_{2}$ error as a function of time for $p=3$ (top) and $p=4$ (bottom) for the compact (left) and standard (right) formulations on a grid of $100$ elements for different values of interior penalty parameter $\eta_{\rm IP}$. Color-gradient indicates different values of $\eta_{\rm IP}$: from purple to light orange, $\eta_{\rm IP}=0.0,0.01,0.02,0.05,0.1,0.5,1.0$.}
\label{fig:nl_diff_sol_ip_100}
\end{figure}

For this particular test case, we also consider the impact of an additional interior penalty within both formulations. In figure~\ref{fig:nl_diff_sol_ip_100} we show the $L_{2}$-error as a function of time when both formulations are equipped with  interior penalty terms. We consider $p=3$ and $p=4$ on a one-dimensional grid made of $100$ elements. 

First of all, by inspecting the right column of figure~\ref{fig:nl_diff_sol_ip_100} we can notice that without interior penalty the errors of the standard formulation are significantly higher. This, of course, reflects what is shown in figure~\ref{fig:nl_diff_err}, where the same curves are reported. The addition of an interior penalty term, regarding the standard approach, helps in mitigating numerical oscillations and initially decreases the error. This, instead, does not seem to be necessary for the compact formulation. The inclusion of an interior penalty term, in combination with the compact approach, does not significantly decrease the error as the compact formulation is already sufficiently robust to suppress numerical oscillations.

Finally, notice that after a certain threshold of the interior penalty parameter is reached, the error stops decreasing and it even slightly increases for relatively large values of $\eta_{\rm IP}$. This happens for both formulations. In fact, as shown in figure~\ref{fig:eigen_diff_tau}, for large values of $\eta_{\rm IP}$ the two schemes start to behave very similarly as the contribution of the interior penalty terms becomes dominant with respect to any particular choice of intermediate state for the viscous fluxes.  

\begin{figure}
\centering
\includegraphics[width=\textwidth]{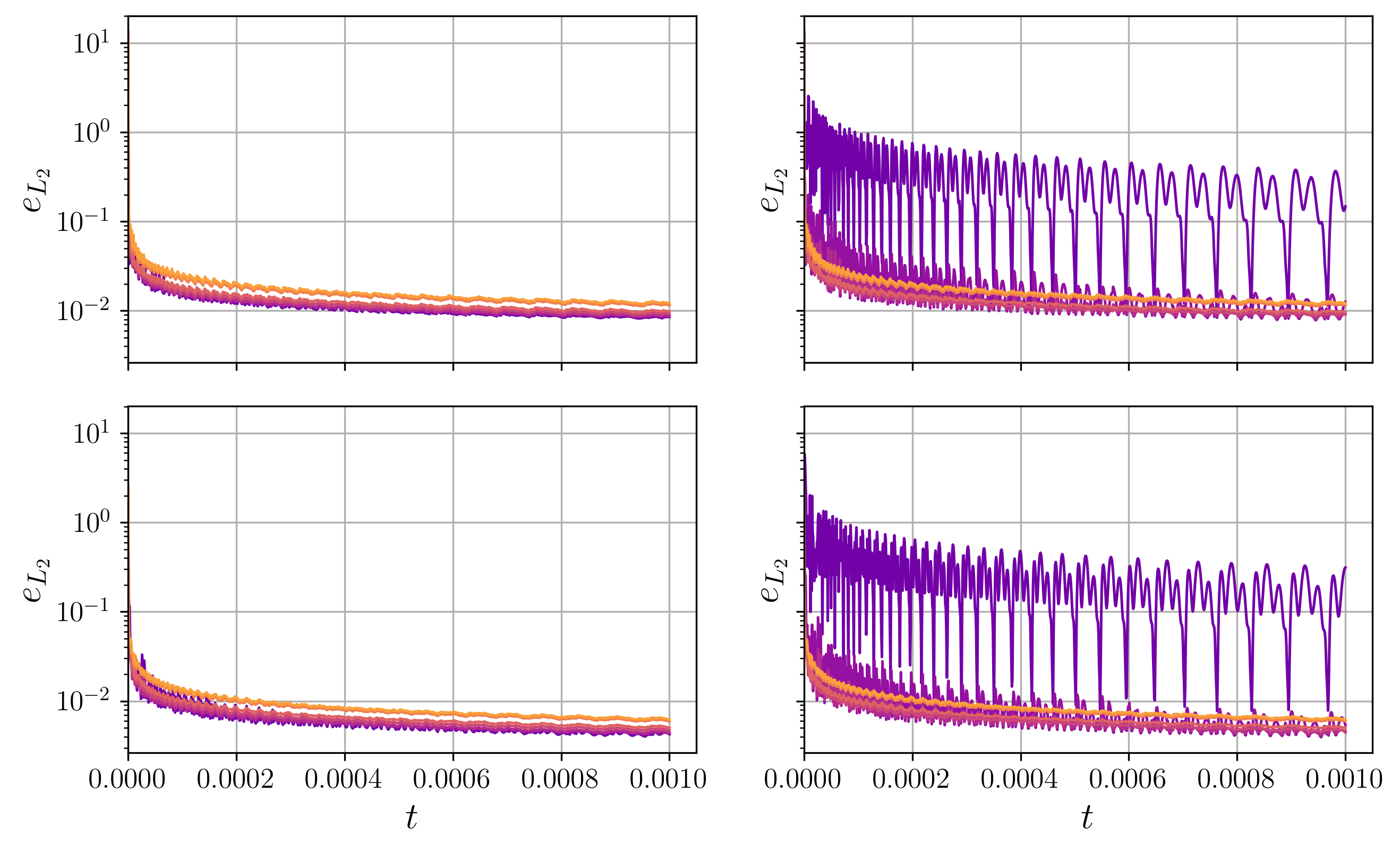}
\caption{Evolution of the $L_{2}$ error as a function of time for $p=3$ (top) and $p=4$ (bottom) for the compact (left) and standard (right) formulations on a grid of $200$ elements for different values of interior penalty parameter $\eta_{\rm IP}$. Color-gradient indicates different values of $\eta_{\rm IP}$: from purple to light orange, $\eta_{\rm IP}=0.0,0.01,0.02,0.05,0.1,0.5,1.0$.}
\label{fig:nl_diff_sol_ip_200}
\end{figure}

When we further increase the resolution to $200$ elements (see figure~\ref{fig:nl_diff_sol_ip_200}) for the same polynomial orders, the addition of interior penalty only deteriorates the compact formulation, leading to an increase in numerical error. The standard formulation, instead, greatly benefits from the use of interior penalty in order to suppress numerical oscillations.
%
\subsection{Taylor-Green vortex flow}
To explore the performance of the compact formulation for the diffusive flux in more realistic three-dimensional settings, we consider a series of ILES of the Taylor-Green vortex (TGV) flow~\citep{brachet:83} at Reynolds numbers of $1\,600$ and $5\,000$.

For this test, the full set of three-dimensional Navier-Stokes equations are integrated:
\begin{equation}
\pd{\vect{U}}{t} + \divg{ [\tens{F}_{\rm c}(\vect{U}) - \tens{F}_{\rm v}(\vect{U},\grad{\vect{U}}) ]} = \vect{0},
\end{equation}
where $\vect{U} = (\rho,\; \rho\vect{u},\; \rho E)^\top$ is the vector of conservative variables and $\tens{F}_{\rm c}$, $\tens{F}_{\rm v} \in \mathbb{R}^{5 \times 3}$ are the inviscid and viscous fluxes, respectively. These read:
\begin{equation}
\tens{F}_{\rm c} =
\begin{pmatrix}
\rho \vect{u}^\top\\
\rho \vect{u} \otimes \vect{u} + p \tens{I}\\
(\rho E + P)\rho \vect{u}^\top
\end{pmatrix},\quad\text{and}\quad
\tens{F}_{\rm v} =
\begin{pmatrix}
\vect{0}^\top\\
\tens{\tau}\\
\vect{u}^\top \cdot \tens{\tau} - \kappa \grad{T}^\top
\end{pmatrix},
\end{equation}
where $\rho$ is the density, $\vect{u}$ is the velocity vector, $P$ is the pressure, $E$ is the total energy, $\kappa$ is the thermal conductivity and $\tens{I} \in \mathbb{R}^{3 \times 3}$ is the identity matrix.
Considering a Newtonian fluid satisfying the ideal gas law, the above equations are closed by the following constitutive relations:
\begin{equation}
\tens{\tau} = 2\mu \tens{S} + \lambda (\divg{\vect{u}})\tens{I},\quad\text{and}\quad
P = (\gamma - 1)( \rho E - \tfrac{1}{2}\rho \vect{u} \cdot \vect{u} ),
\end{equation}
with $\tens{S} = \frac{1}{2}(\grad{\vect{u}} + \grad{\vect{u}}^\top)$ the deformation tensor, $\mu$ the dynamic viscosity, $\lambda$ the second coefficient of viscosity---which we set equal to zero~\citep{buresti:15,rajagopal:13}---and $\gamma$ the heat capacity ratio, set equal to 1.4 as for a diatomic gas.
We consider air, for which we set $\kappa = \mu C_P / \Pr$ with $C_P$ the relevant heat coefficient at constant pressure and $\Pr = 0.72$, and we assume the validity of the Sutherland law to compute $\mu$ as a function of the temperature.

The inviscid fluxes are computed at the elements' interfaces via the Roe flux with entropy fix~\citep{roe:81,harten:83}, whereas time integration is performed using the fourth-order, five-stages, strong stability preserving (RK45-SSP) scheme~\citep{spiteri:02}.
Concerning the viscous fluxes, the proposed compact formulation is applied in the calculation of the gradients of conservative variables, from which primitive variables' gradients are then obtained as needed. For the present tests, in particular, to better asses the different performances of the standard and compact formulations, a centered viscous flux is adopted, computed on the average solution and the average gradient at interfaces and the IP term is not applied unless stated otherwise.

The computational domain $\Omega$ is a cube of side length $H = 2 \pi L$, with $L$ a reference size.
The flow is initialized using the following relations:
\begin{align}
\rho(\vect{x}^*,0) &= \rho_0,\\
\vect{u}(\vect{x}^*,0) &= u_0
\begin{pmatrix}
\sin (x^*) \cos (y^*) \cos (z^*)\\
-\cos (x^*) \sin (y^*) \cos (z^*)\\
0
\end{pmatrix},\\
p(\vect{x}^*,0) &= p_0 + \frac{\rho_0 u_0^2}{16} [\cos (2 x^*) \cos (2 y^*)][\cos (2 z^*) + 2],
\end{align}
where $\vect{x}^* = (\vect{x} - \vect{x}_{\rm c})/L$, $\vect{x}_{\rm c}$ being the location of the center of the domain, and $\rho_0$, $u_0$, $p_0$ are the reference density, velocity and pressure, respectively. These are set, together with the reference viscosity $\mu_0$, to obtain the required Reynolds and Mach numbers according to the following relations:
\begin{equation}
\mu_0 = \frac{\rho_0 u_0 L}{\Rey},\quad\text{and}\quad p_0 = \frac{\rho_0}{\gamma}\bigg( \frac{u_0}{\Ma} \bigg)^2,
\end{equation}
where $\Ma$ was set below 0.3 to remain in the incompressible regime and, as it is common practice, the flow was let evolve for a total normalized time of $t^* = t u_0 / L = 20$.
In the following analysis, we adopt as flow diagnostic quantities the kinetic energy and the enstrophy, integrated over the whole volume. These are defined, respectively, as:
\begin{equation}
E_k(t) = \frac{1}{2}\int_{\Omega} \vect{u}\cdot \vect{u} \,\ud \Omega,\quad\text{and}\quad 
\omega^2(t) = \int_{\Omega} \vect{\omega}\cdot \vect{\omega} \,\ud \Omega,\quad\text{with}\quad 
\vect{\omega} = \curl{\vect{u}}.
\end{equation}

\begin{table}[!ht]
\centering

\begin{tabular}{|l|c|c|c|c|c|}
\hline
\# DoF & flux & $p=5$ & $p=6$ & $p=7$ & $p=8$\\
\hline

\multirow{2}{*}{$128^3$} 
& standard 
& \cellcolor{stablegreen}stable
& \cellcolor{stablegreen}stable
& \cellcolor{unstablered}\bf unstable
& \cellcolor{stablegreen}stable
\\

& compact
& \cellcolor{stablegreen}stable
& \cellcolor{stablegreen}stable
& \cellcolor{stablegreen} \bf stable
& \cellcolor{stablegreen}stable
\\
\hline

\multirow{2}{*}{$64^3$} 
& standard
& \cellcolor{stablegreen}stable
& \cellcolor{stablegreen}stable
& \cellcolor{unstablered}unstable
& \cellcolor{unstablered}unstable
\\

& compact
& \cellcolor{stablegreen}stable
& \cellcolor{stablegreen}stable
& \cellcolor{stablegreen}stable
& \cellcolor{unstablered}unstable
\\

\hline
\end{tabular}

\caption{Stability of the TGV flow simulation at $\Rey=1\,600$. The simulations discussed in details in this section are highlighted in bold.}
\label{tab:tgv:1600}

\end{table}

\begin{table}[!ht]
\centering


\begin{tabular}{|l | c | c | c | c | c|}
\hline

\# DoF & flux & $p=4$ & $p=5$ & $p=6$ & $p=7$\\

\hline

\multirow{2}{*}{$96^3$} 
& standard 
& \cellcolor{stablegreen}stable 
& \cellcolor{stablegreen}\bf stable 
& \cellcolor{stablegreen}stable 
& \cellcolor{unstablered}unstable
\\

& compact 
& \cellcolor{stablegreen}stable 
& \cellcolor{stablegreen}\bf stable 
& \cellcolor{stablegreen}stable 
& \cellcolor{unstablered}unstable
\\

\hline

\multirow{2}{*}{$84^3$} 
& standard 
& \cellcolor{stablegreen}stable 
& \cellcolor{stablegreen}stable 
& \cellcolor{unstablered} \bf unstable 
& \cellcolor{unstablered}unstable
\\

& compact 
& \cellcolor{stablegreen}stable 
& \cellcolor{stablegreen}stable 
& \cellcolor{stablegreen}\bf stable 
& \cellcolor{unstablered}unstable
\\

\hline

\multirow{2}{*}{$72^3$} 
& standard 
& \cellcolor{stablegreen}stable 
& \cellcolor{stablegreen}stable 
& \cellcolor{unstablered}unstable 
& \cellcolor{unstablered}unstable
\\

& compact 
& \cellcolor{stablegreen}stable 
& \cellcolor{stablegreen}stable 
& \cellcolor{unstablered}unstable 
& \cellcolor{unstablered}unstable
\\

\hline

\multirow{2}{*}{$64^3$} 
& standard 
& \cellcolor{stablegreen}stable 
& \cellcolor{stablegreen}stable 
& \cellcolor{unstablered}unstable 
& \cellcolor{unstablered}unstable
\\

& compact 
& \cellcolor{stablegreen}stable 
& \cellcolor{stablegreen}stable 
& \cellcolor{stablegreen}stable 
& \cellcolor{unstablered}unstable
\\

\hline

\end{tabular}

\caption{Stability of the TGV flow simulation at $\Rey=5\,000$. The simulations discussed in details in this section are highlighted in bold.}

\label{tab:tgv:5000}

\end{table}

Several computations  were performed at the two selected Reynolds numbers, namely $1\,600$ and $5\,000$, to highlight differences between the standard BR1 flux and the proposed compact flux; see tables~\ref{tab:tgv:1600} and~\ref{tab:tgv:5000}. For reference, similar simulations, using the same code, are reported by~\citet{chapelier2016spectral,chapelier:16b}.
In particular, the selected resolutions and Reynolds were chosen to test the behavior of the compact flux in the case of well-resolved ($\Rey = 1\,600$) and under-resolved ($\Rey = 5\,000$) ILES, respectively.
It is worth stressing that, given the relatively high Reynolds numbers adopted, the impact of the convective scheme cannot be considered negligible and the presented results shall be regarded as representative of a \emph{typical} numerical setup for (I)LES.
A more detailed study of the combined performances of viscous and inviscid fluxes, for different choices of these last, is beyond the scope of the present paper. 

Looking at what is reported in tables~\ref{tab:tgv:1600} and~\ref{tab:tgv:5000}, a first evident observation is that, depending on the resolving power of the adopted setups, some computations failed. This was especially the case for higher-order simulations at low resolutions, as the last column of table~\ref{tab:tgv:5000} testifies.
In general, as expected, the instability was observed around the time of the peak of enstrophy, when the most intense small-scale activity occurs. Indeed, in such scenario, the scheme's ability to damp high, unresolved wavenumbers is a crucial aspect.
Overall, in agreement with previous analyses (\cf~section~\ref{sec:2}), the compact flux appears to have an advantage, especially at low resolutions. Among all the computations in tables~\ref{tab:tgv:1600} and~\ref{tab:tgv:5000}, we selected a few representative runs which are detailed below. 

\begin{figure}[t]
 \centering  
\includegraphics[width=\textwidth]{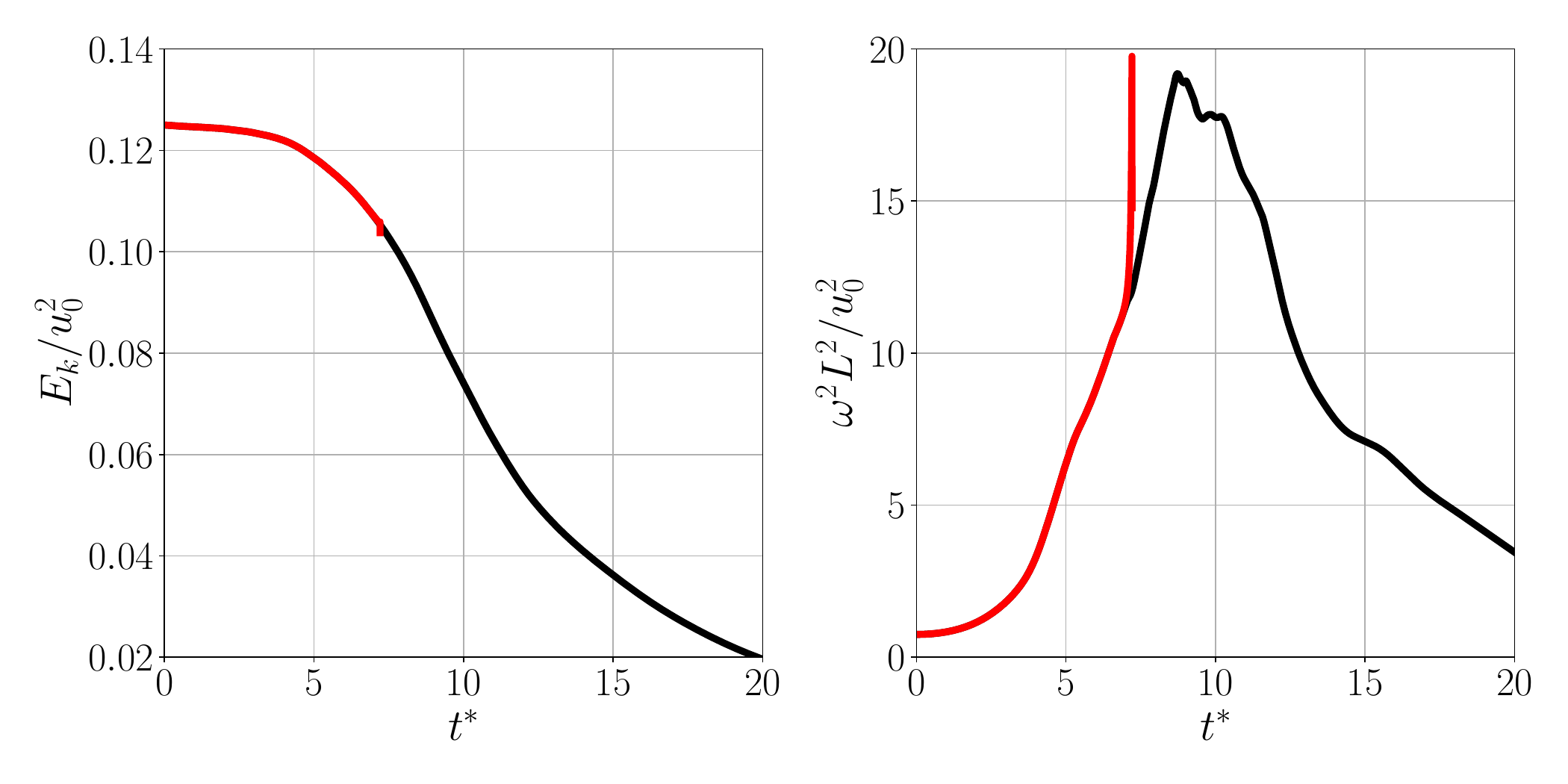}
\caption{Time evolution of normalized kinetic energy (left) and enstrophy (right) for the TGV $16^3p7$, $\Rey = 1\,600$ test case: black lines, compact flux; red lines, standard flux.}
   \label{fig:tgv:16p7:1600}
\end{figure}

\begin{figure}[t]
 \centering  
 \begin{subfigure}{0.48\textwidth}
 \includegraphics[width=\textwidth]{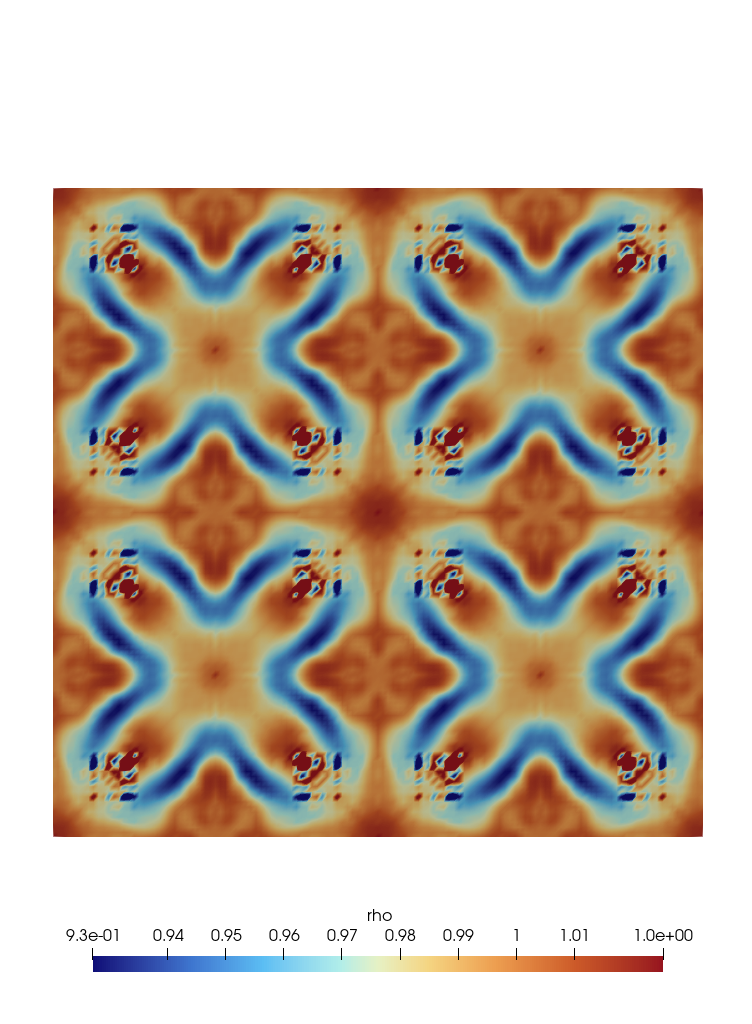}
 \end{subfigure}
 \begin{subfigure}{0.48\textwidth}
 \includegraphics[width=\textwidth]{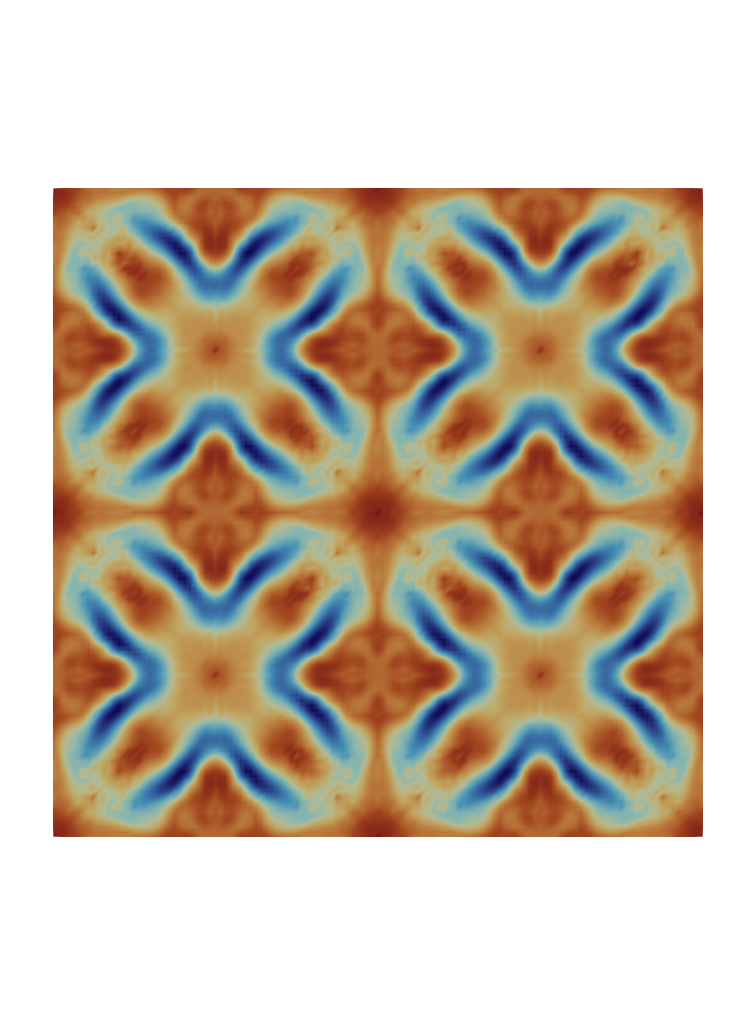}
 \end{subfigure}
\caption{Contours of density at $t^* = 7.2$ for the TGV $16^3p7$, $\Rey = 1\,600$ test case: left, standard flux; right, compact flux. Plane location shown by red outline in figure~\ref{fig:tgv:16p7:1600:Q}.}
   \label{fig:tgv:16p7:1600:rho}
\end{figure}

\begin{figure}[t]
 \centering  
 \begin{subfigure}{0.48\textwidth}
 \includegraphics[width=\textwidth]{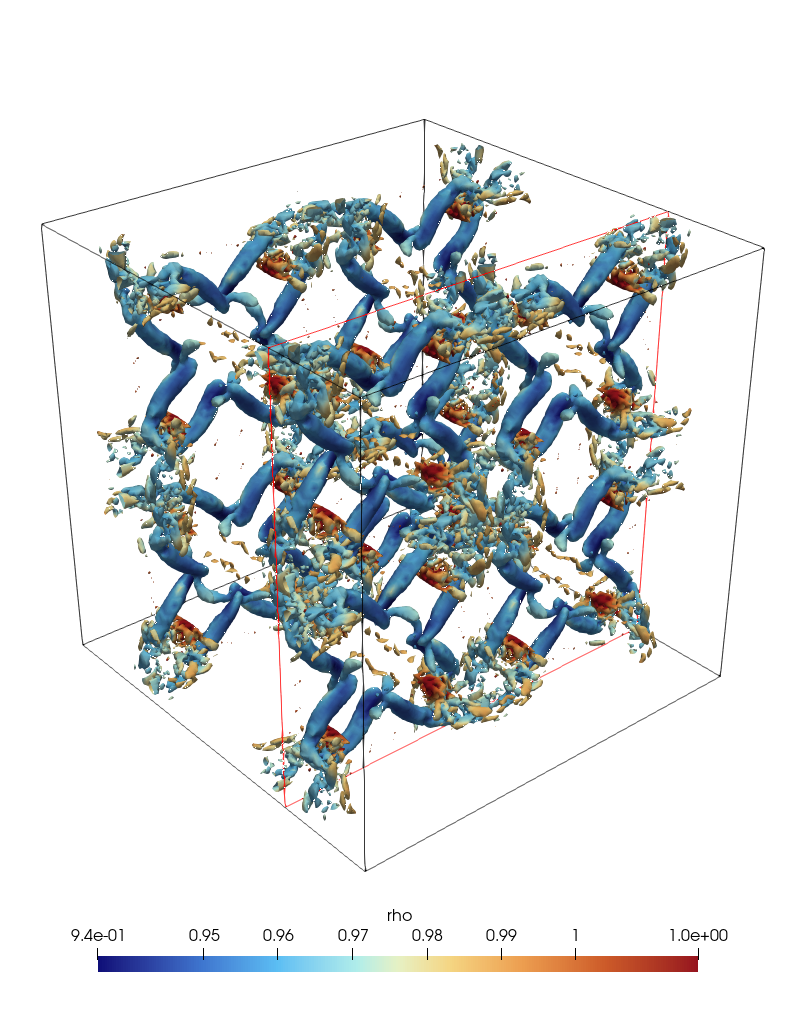}
 \end{subfigure}
 \begin{subfigure}{0.48\textwidth}
 \includegraphics[width=\textwidth]{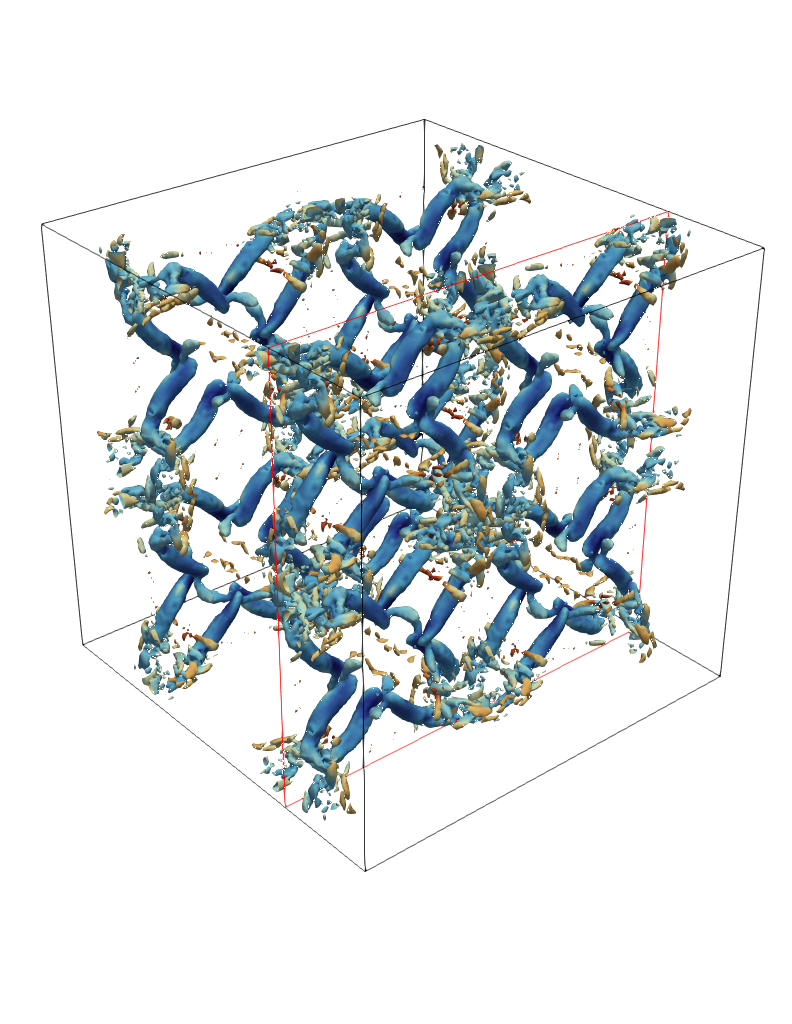}
 \end{subfigure}
\caption{Iso-surfaces of $Q$-criterion  at $t^* = 7.2$ for the TGV $16^3p7$, $\Rey = 1\,600$ test case: left, standard flux; right, compact flux.}
   \label{fig:tgv:16p7:1600:Q}
\end{figure}

Starting with the relatively well-resolved test at $\Rey = 1\,600$, performed with $16^3$ 7-th order elements (\ie, $128^3$ DoF), the time history of the  kinetic energy and enstrophy, \cf~figure~\ref{fig:tgv:16p7:1600}, shows the failure of the run with the standard flux at a normalized time $t^* = 7.2$, well before the peak of enstrophy.
As shown in figures~\ref{fig:tgv:16p7:1600:rho} and~\ref{fig:tgv:16p7:1600:Q}, which depict, respectively, the density contours and the iso-surfaces of $Q$-criterion at the moment of the numerical instability, the computation with the standard flux shows nonphysical peaks of density in regions where small-scale vortices are developing (see red spots in figure~\ref{fig:tgv:16p7:1600:Q}).
Similar behaviors are observed at lower resolutions, $64^3$ DoF, except that the highest order computation---characterized by a very low level of numerical dissipation overall---resulted unstable for both fluxes.
Note that, even the addition of the IP term with the recommended value of $\eta_{\rm IP} = 1.0$ could not stabilize the simulations. At that combination of order and resolution, the use of an explicit sub-grid scale model is quite possibly required.
%
%

Switching to the higher Reynolds number, a first simulation was performed on a grid with $12^3$ 7-th order elements, for a total of $96^3$ DoF. This resolution was selected to compare with an identical simulation that was reported unstable with the standard flux without some form of regularization or sub-grid scale modeling~\citep{chapelier2016spectral}.
Indeed, the same outcome is obtained here: even the improved damping at high wavenumbers brought by the compact formulation is not able to counteract the energy pile-up at small scales and the computation goes unstable at about the moment the peak of enstrophy occurs (more precisely, $t^{*} = 9.1$ for the compact flux and $t^{*} = 9.0$ for the standard flux). 
Checking the time history of the  kinetic energy and enstrophy, not shown, we observe marginal differences in the former and  higher values in the latter for $t^{*} \ge 8.0$ when the standard flux is adopted.

\begin{figure}[t]
 \centering  
\includegraphics[width=\textwidth]{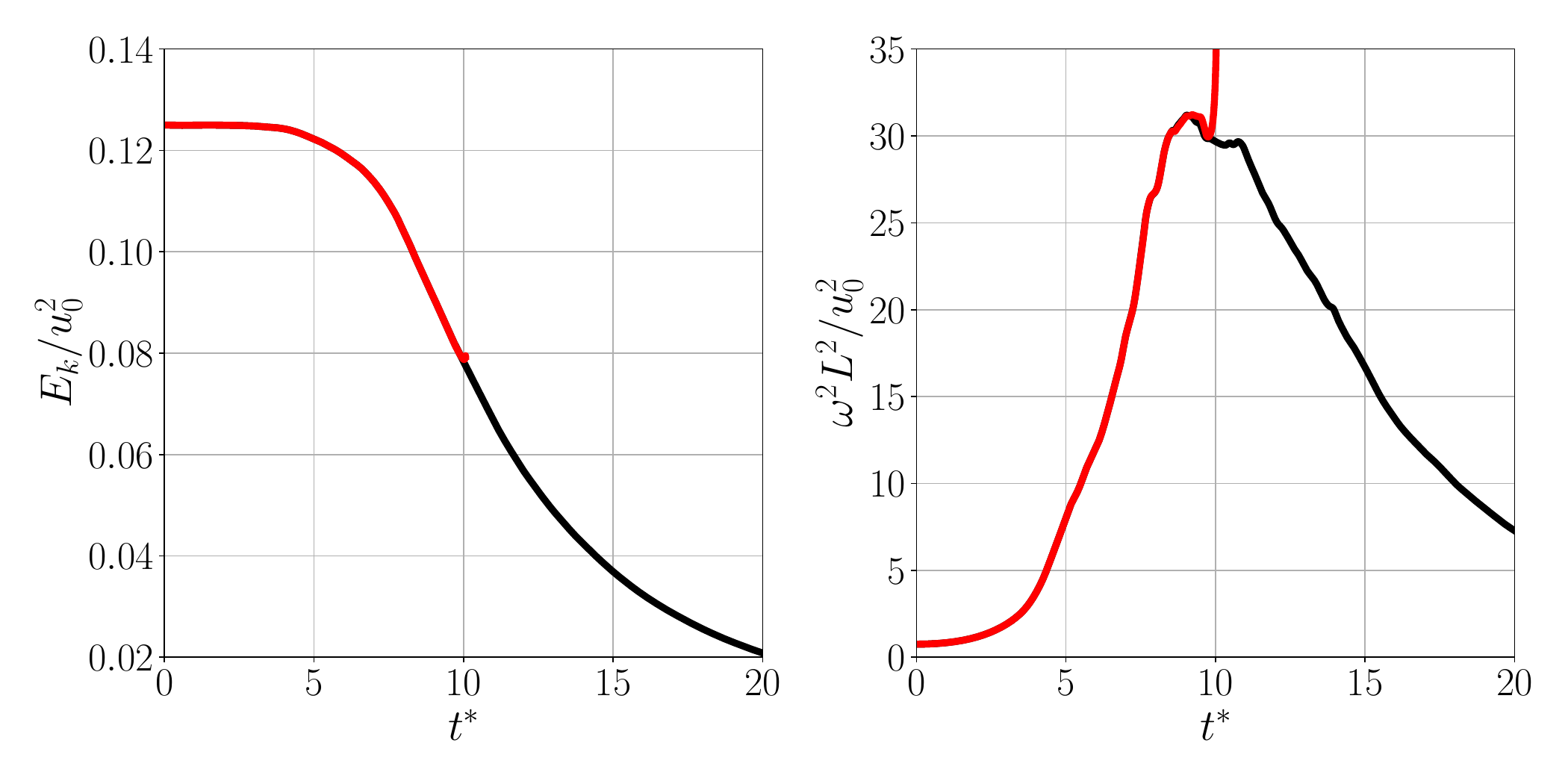}
\caption{Time evolution of normalized kinetic energy (left) and enstrophy (right) for the TGV $12^3p6$, $\Rey = 5\,000$ test case: black lines, compact flux; red lines, standard flux.}
   \label{fig:tgv:12p6}
\end{figure}

\begin{figure}[t]
 \centering  
 \begin{subfigure}{0.48\textwidth}
 \includegraphics[width=\textwidth]{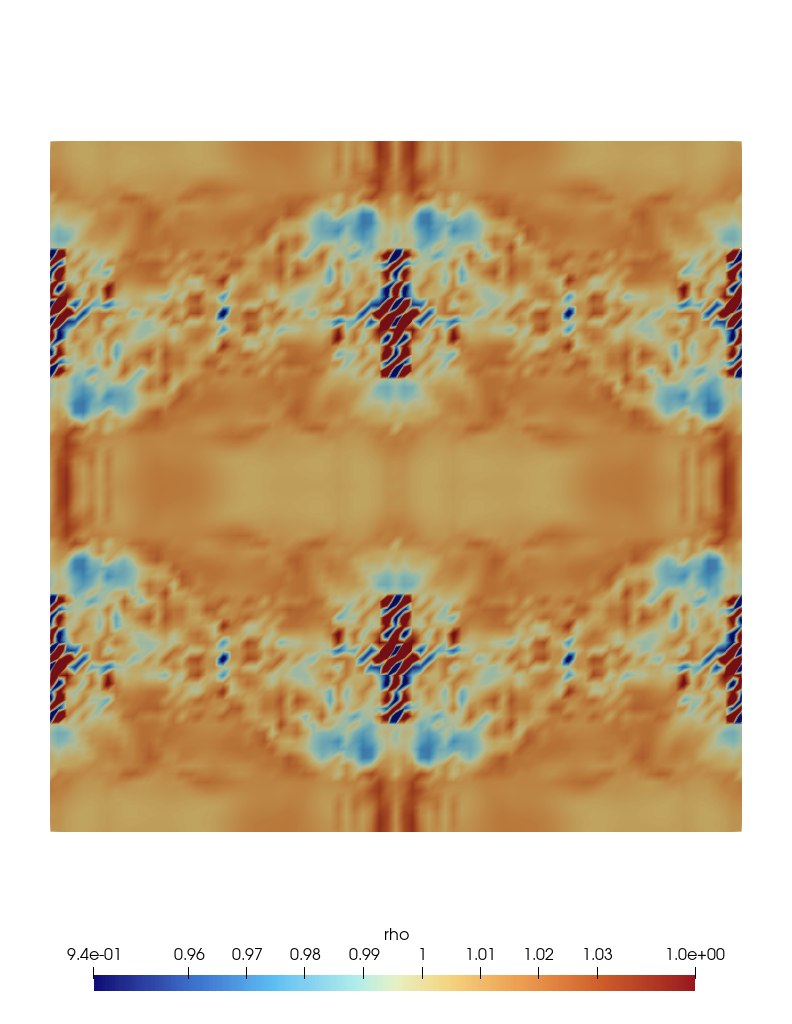}
 \end{subfigure}
 \begin{subfigure}{0.48\textwidth}
 \includegraphics[width=\textwidth]{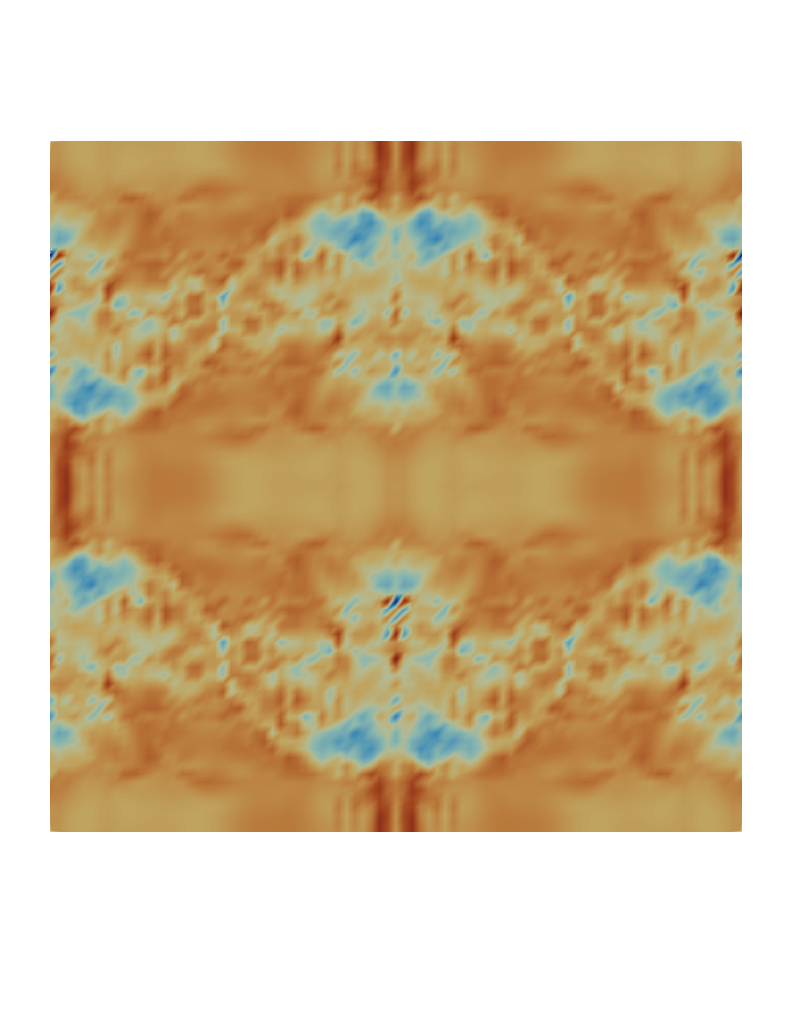}
 \end{subfigure}
\caption{Contours of density at at $t^* = 10.06$ for the TGV $12^3p6$, $\Rey = 5\,000$ test case: left, standard flux; right, compact flux. Plane location shown by red outline in figure~\ref{fig:tgv:12p6:5000:Q}.}
   \label{fig:tgv:12p6:5000:rho}
\end{figure}

\begin{figure}[t]
 \centering  
 \begin{subfigure}{0.48\textwidth}
 \includegraphics[width=\textwidth]{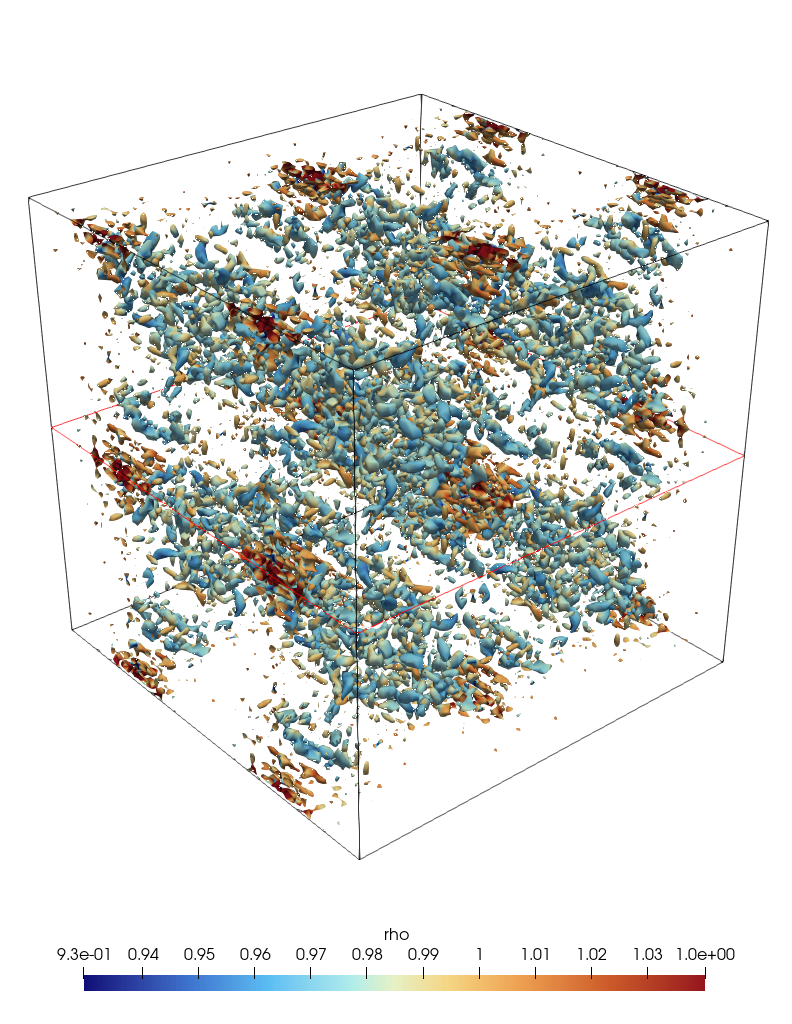}
 \end{subfigure}
 \begin{subfigure}{0.48\textwidth}
 \includegraphics[width=\textwidth]{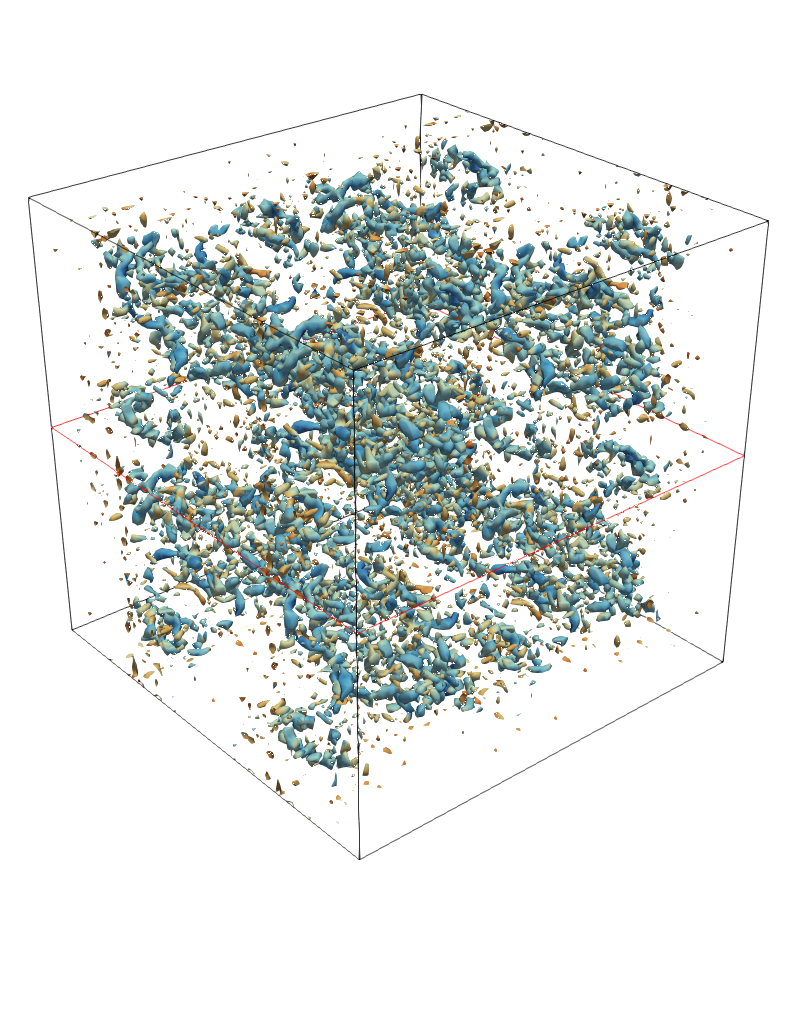}
 \end{subfigure}
\caption{Iso-surfaces of $Q$-criterion  at $t^* = 10.06$ for the TGV $12^3p6$, $\Rey = 5\,000$ test case: left, standard flux; right, compact flux.}
   \label{fig:tgv:12p6:5000:Q}
\end{figure}

To gauge the performances of the two fluxes close to the limit of numerical stability, despite a marginal loss of resolution, $84^3$ DoF, a 6-th polynomial order was then selected to slightly increase numerical dissipation.  
The time history of  kinetic energy and enstrophy are depicted in figure~\ref{fig:tgv:12p6}.
At this resolution and with this relatively high polynomial order, the spectral difference scheme is not expected to promote high levels of numerical dissipation~\citep{jameson:14,chapelier:16b,chapelier:17b} and the number of selected degrees of freedom is definitely too low to resolve the fine structures which are generated when the peak of enstrophy occurs.
Indeed, the simulation with the standard flux still fails, for $t^* \sim 10$, with a sudden increase in enstrophy. 
Contours of density and iso-surfaces of $Q$-criterion right before the instability are depicted in figures~\ref{fig:tgv:12p6:5000:rho} and~\ref{fig:tgv:12p6:5000:Q}, respectively, where numerical artifacts can be distinctly seen in regions of strong small-scales activity.\footnote{Notice that some mild artifacts are also visible in the simulation with the compact flux (see the region at a distance of about $H/3$ from the bottom in figure~\ref{fig:tgv:12p6:5000:rho}), but these were not severe enough to destabilize the simulation.}
Clearly, too much energy is passed down at small scales which are not sufficiently damped by numerical dissipation.
The use of the compact scheme, on the other hand, maintain the simulation stable.
In view of previous results by the eigenanalysis, this can be explained by the improved damping of the spurious modes provided by the compact formulation.
%

\begin{figure}[t]
 \centering  
\includegraphics[width=\textwidth]{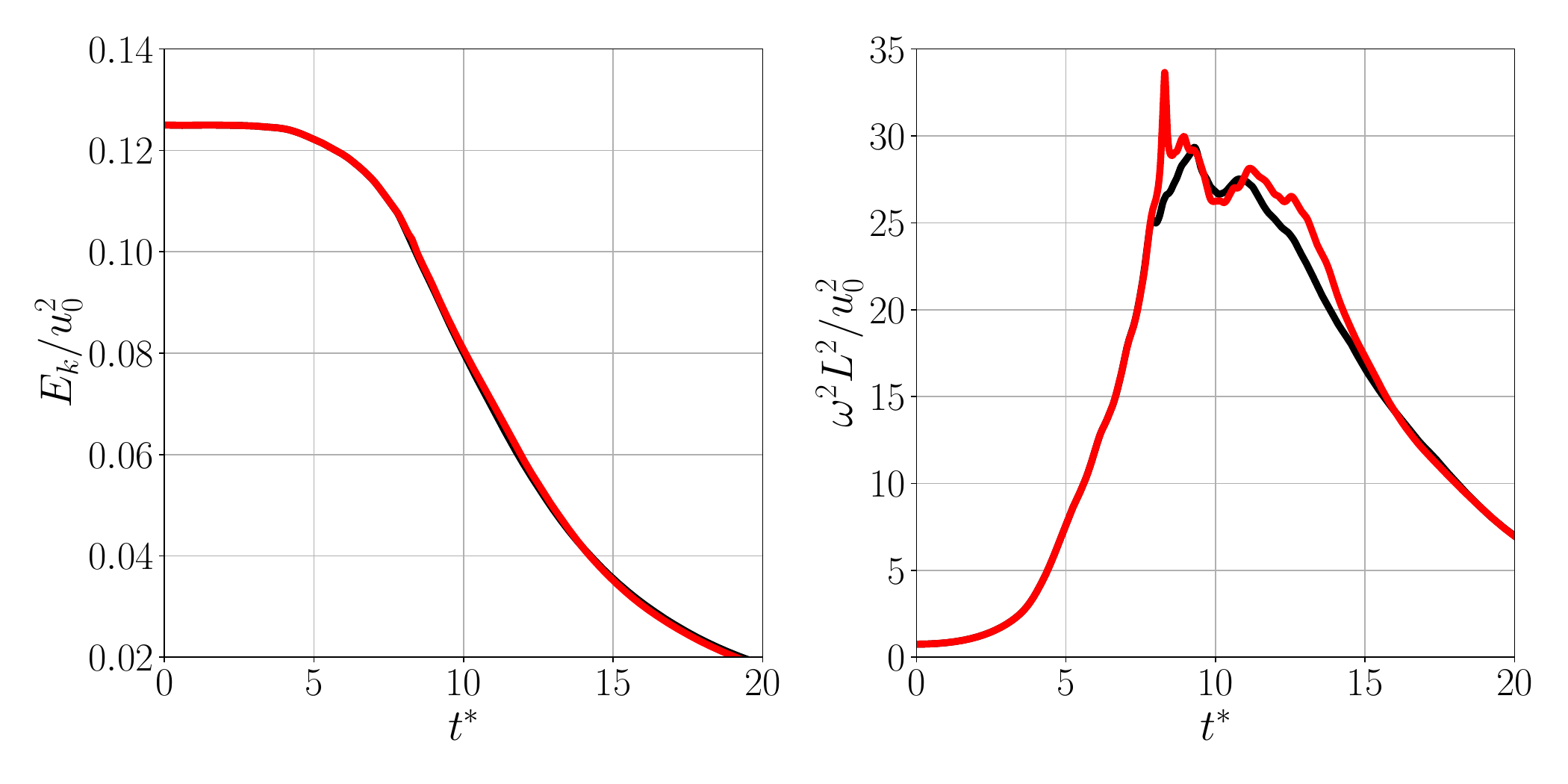}
\caption{Time evolution of normalized kinetic energy (left) and enstrophy (right) for the TGV $16^3p5$, $\Rey = 5\,000$ test case: black lines, compact flux; red lines, standard flux.}
   \label{fig:tgv:16p5:5000}
\end{figure}

When the order of the scheme is reduced, the increased level of numerical dissipation injected allows stable computations with both formulations of the viscous fluxes. 
This is the case, for instance, for the $\Rey = 5\,000$, $96^3$ DoF test case performed on $16^3$ 5-th order elements, in which both the compact and the standard scheme were stable (\cf~table~\ref{tab:tgv:5000}).
The relevant results, in terms of kinetic energy and enstrophy, are depicted in figure~\ref{fig:tgv:16p5:5000}.
Despite the stability of both fluxes, some marked differences can be observed when looking at the enstrophy: the standard flux appears to promote a significantly stronger small-scales activity compared to the compact flux.
In particular, a peak in normalized enstrophy of almost 35 is observed in the computation performed with the standard flux, a value which is well beyond the maximum value registered---with both types of fluxes---in the higher-order, more resolved test (\cf~figure~\ref{fig:tgv:12p6}).
Such a high value of enstrophy on a less resolved computation is, of course, completely unexpected and quite possibly due to undamped spurious high-wavenumber modes, as highlighted in previous analyses.

Overall, the tests performed on the TGV, which, as already mentioned, represents a rather \emph{atypical} choice to  ascertain the performances of the viscous part alone, seem to confirm the slightly more dissipative behavior of the compact flux with respect to spurious  high-wavenumber modes that would be otherwise undamped by the BR1 flux. This is especially true for under-resolved simulations and, supposedly, for more challenging setups involving viscous-dominated problems and/or the presence of stiff source terms (\eg, RANS modeling, multiphase or reactive flows ).

\section{Conclusions}\label{sec:4}
In this work, a novel compact formulation for the discretization of second-order fluxes within the Spectral Difference method has been presented. Building on the ideas originally proposed by Huynh for the Flux Reconstruction approach, the new scheme modifies the way the continuous flux used to build the auxiliary gradient variable is constructed at element interfaces. By introducing one-sided, interface-dependent continuous fluxes and averaging their contributions only at the corresponding interface, the resulting formulation reduces the stencil required for the viscous discretization from five to three elements (in the one-dimensional case), while retaining the centered, parameter-free character of the classical BR1 approach. The methodology has been derived in detail for the one-dimensional case and subsequently extended to multiple dimensions, exploiting the tensor-product structure of the SD elements.

A temporal eigenanalysis of the linear diffusion equation has been used to characterize the dissipation properties of the compact scheme in comparison with the standard, extended-stencil BR1 formulation, both with and without the addition of interior penalty terms. The results indicate that the compact formulation is generally more dissipative than the standard approach for high wavenumbers, particularly at low polynomial orders, and that the differences between the two formulations progressively vanish as the interior penalty parameter increases and its contribution becomes dominant. These findings have been further corroborated using a combined-mode Fourier analysis, which confirmed that the compact scheme provides enhanced damping specifically in the medium- and high-wavenumber range, \ie, in the regions most relevant to the robustness and stability of the numerical scheme.

A set of numerical experiments of increasing complexity has then been used to assess the performance and robustness of the proposed approach. A convergence study for the diffusion of a well-resolved Gaussian profile has shown that the compact formulation consistently recovers the expected order of accuracy for all tested polynomial degrees, whereas the standard BR1 scheme exhibits a degraded, sub-optimal convergence rate at even orders of approximation, in agreement with previous findings reported in the DG literature. The diffusion of an under-resolved, localized Dirac's delta has further highlighted the improved robustness of the compact scheme, which significantly reduces the spurious oscillations and the periodically excited high-wavenumber modes observed with the standard formulation, consistently with the predictions of the temporal eigenanalysis. Similar conclusions have been drawn from the nonlinear porous medium equation, where the compact scheme has been shown to consistently yield smaller errors and reduced oscillatory behavior compared to the standard approach, with the addition of interior penalty terms proving beneficial mainly for the standard formulation and, at higher resolutions, even detrimental for the compact one.

Finally, the compact formulation has been assessed in a fully three-dimensional, nonlinear setting by means of implicit large-eddy simulations of the Taylor-Green vortex at $\Rey = 1\,600$ and $5\,000$. In these under-resolved turbulent computations, the compact scheme has been found to remain numerically stable in several configurations for which the standard BR1 formulation fails, owing to its improved ability to damp spurious high-wavenumber content generated by the nonlinear energy cascade. In cases where both formulations remain stable, the compact scheme has also been shown to produce lower levels of small-scale enstrophy activity, consistent with its enhanced dissipative characteristics at high wavenumbers.

Overall, the results presented in this work demonstrate that the proposed compact formulation constitutes an effective and computationally attractive alternative to the classical BR1 scheme for the discretization of second-order fluxes within the SD method. Its reduced stencil makes it particularly well suited for large-scale while its improved accuracy and robustness properties are especially relevant for implicit LES applications in under-resolved flows.

\section*{Dedication}
This paper is dedicated to the memory of Antony Jameson, an extraordinary mind, a mentor and a dear friend, who pioneered the science of computational fluid dynamics. His seminal contributions paved the way for the development of some of the most elegant techniques currently used in scientific computing. GL

\section*{Acknowledgments}
This study was funded by the European Union - NextGenerationEU, in the framework of the iNEST - Interconnected Nord-Est Innovation Ecosystem (iNEST ECS00000043 – CUP G93C22000610007). The views and opinions expressed are solely those of the authors and do not necessarily reflect those of the European Union, nor can the European Union be held responsible for them. NT acknowledge the support by INdAM-GNCS: Istituto Nazionale di Alta Matematica –– Gruppo Nazionale di Calcolo Scientifico. Fruitful discussions with Nicola Clinco about the manuscript and the results therein and the use of the SD solver originally developed by Antony Jameson’s group at Stanford University are gratefully acknowledged.

\appendix
\section{Combined-mode Fourier analysis}\label{sec:5}
The combined-mode Fourier analysis has been recently proposed by \cite{alhawwary2018fourier,alhawwary2019study,alhawwary2020combined} to avoid the, sometimes misleading, interpretation of the $p+1$ eigen-curves obtained from the classical temporal eigenanalysis applied to DSEMs.

The mathematical formulation builds upon the standard temporal eigenanalysis. The eigenvalue problem 
\begin{equation}
-\iota \widetilde{\omega} \vect{u}_e^{\rm S} = \tens{\mathcal{B}}(\theta)\cdot\vect{u}_e^{\rm S},
\end{equation}
is derived as a representation of the numerical discretization of the diffusion equation using the SD methods equipped with the two numerical fluxes herein considered.

The semi-discrete operator admits $p+1$ eigenpairs, $(\widetilde{\omega}_j,\vect{\mu}_j)$, $j=0,\ldots,p$, which form a basis
for the solution space. Therefore, the initial Fourier mode can be expanded as
\begin{equation}
\vect{u}_e^{\rm S}(t_{0}) = \sum_{j=0}^{p}\gamma_j\vect{\mu}_j = \tens{\mathcal{N}} \cdot \vect{\gamma},
\end{equation}
where $\tens{\mathcal{N}}=[\vect{\mu}_0,\ldots,\vect{\mu}_p]$ is the eigenvector matrix and $\vect{\gamma}$ contains the expansion coefficients.
The temporal evolution of the numerical solution is then obtained as
\begin{equation}
\vect{u}_e^{\rm S}(t)= \sum_{j=0}^{p} \gamma_j\vect{\mu}_j \exp(-\iota\widetilde{\omega}_j t).
\end{equation}

Unlike the classical eigenmode analysis, where the numerical properties are evaluated by selecting a single eigenmode, the combined-mode approach retains the contribution of all eigenvectors involved in the expansion of the original Fourier mode. The effective amplification factor associated with the complete Fourier mode is defined as
\begin{equation}
G(\theta,t) = \frac{ [ \vect{u}_e^{\rm S}(t_{0})]^{*} \cdot \vect{u}_e^{\rm S}(t)}{[\vect{u}_e^{\rm S}(t_{0})]^{*} \cdot \vect{u}_e^{\rm S}(t_{0})},
\end{equation}
with $[\;\cdot\;]^{*}$ denoting the complex conjugate operator.

The corresponding combined numerical frequency is obtained from
\begin{equation}
\widetilde{\omega} = -\frac{\iota}{t}\log\left[G(\theta,t)\right].
\end{equation}

This procedure accounts for the interaction among all numerical modes and provides the dissipative and dispersive properties of the complete Fourier component, without requiring the identification of a single physical eigenmode.
\begin{figure}[h!]
 \centering  
 \begin{subfigure}{0.48\textwidth}
     \includegraphics[width=\textwidth]{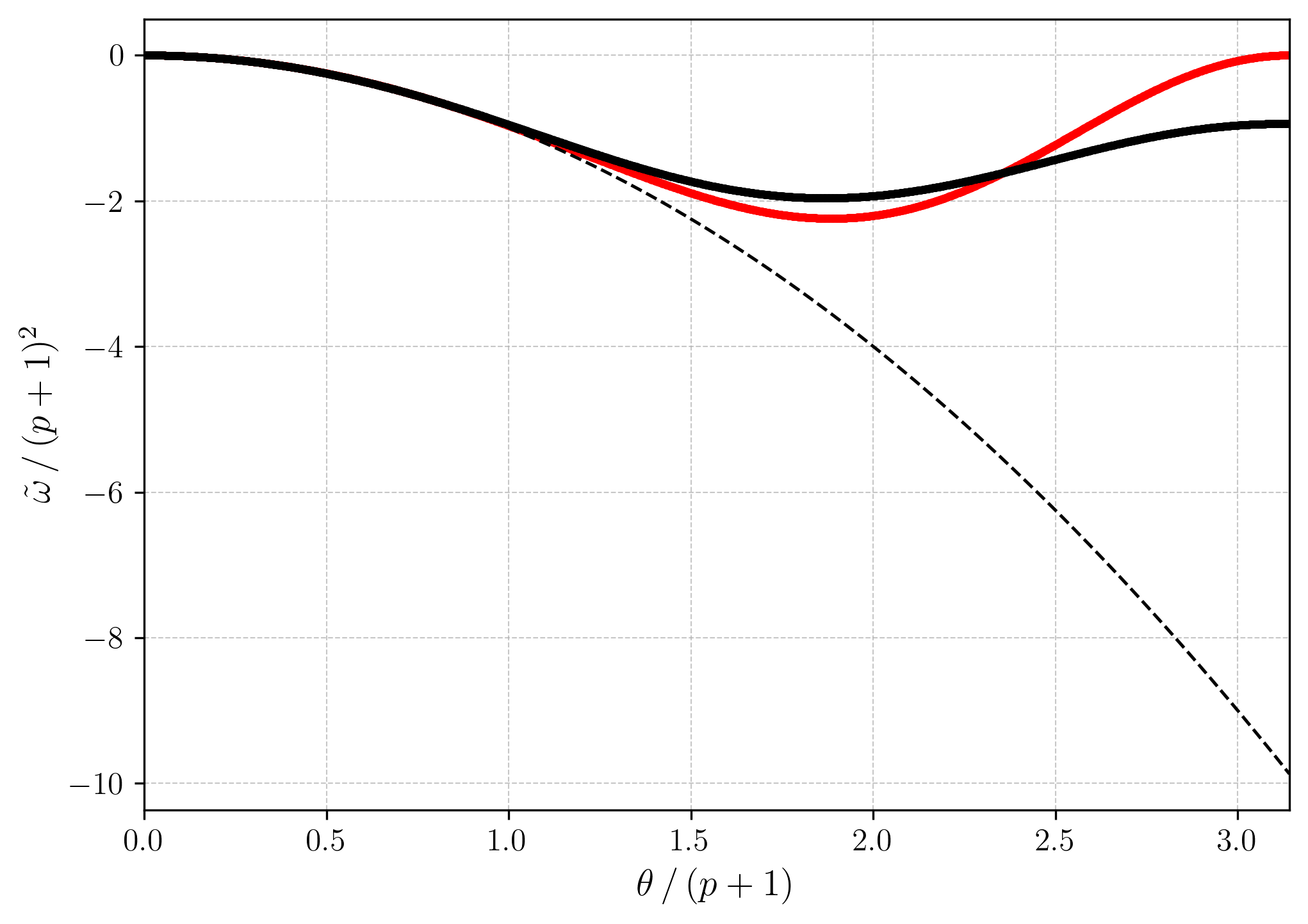}
 \end{subfigure}
 \begin{subfigure}{0.48\textwidth}
     \includegraphics[width=\textwidth]{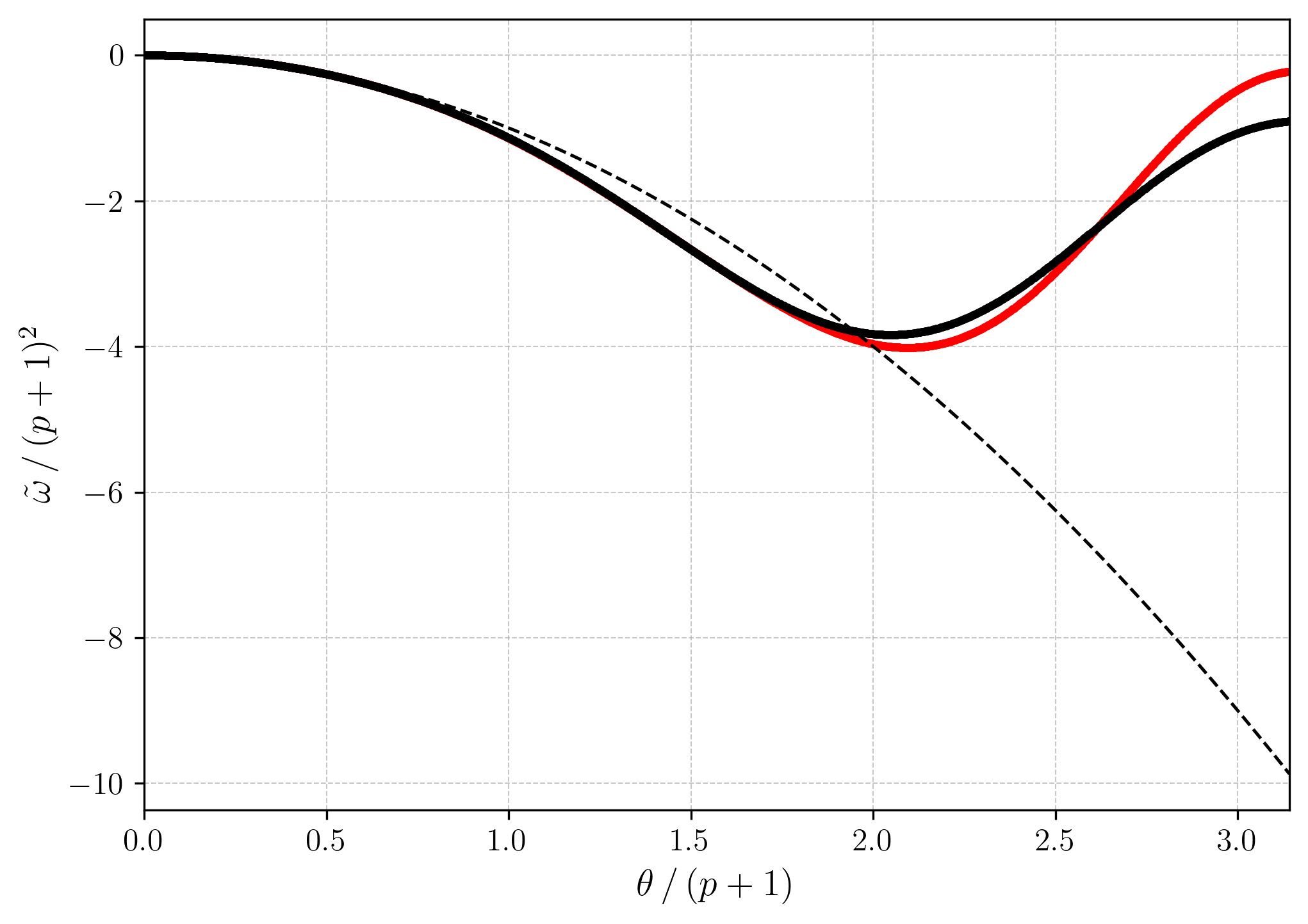}
 \end{subfigure}
 \begin{subfigure}{0.48\textwidth}
     \includegraphics[width=\textwidth]{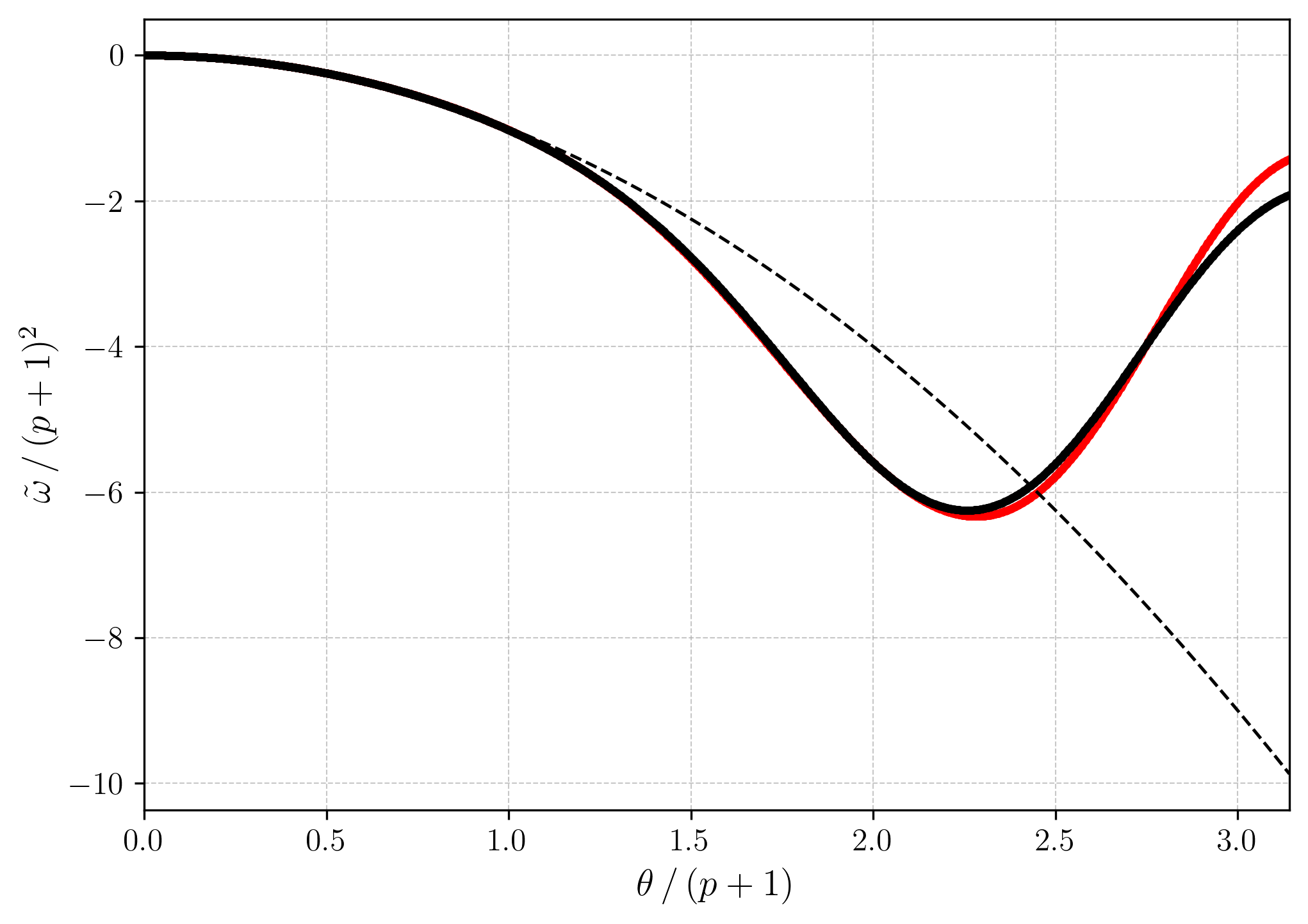}
 \end{subfigure}
 \begin{subfigure}{0.48\textwidth}
     \includegraphics[width=\textwidth]{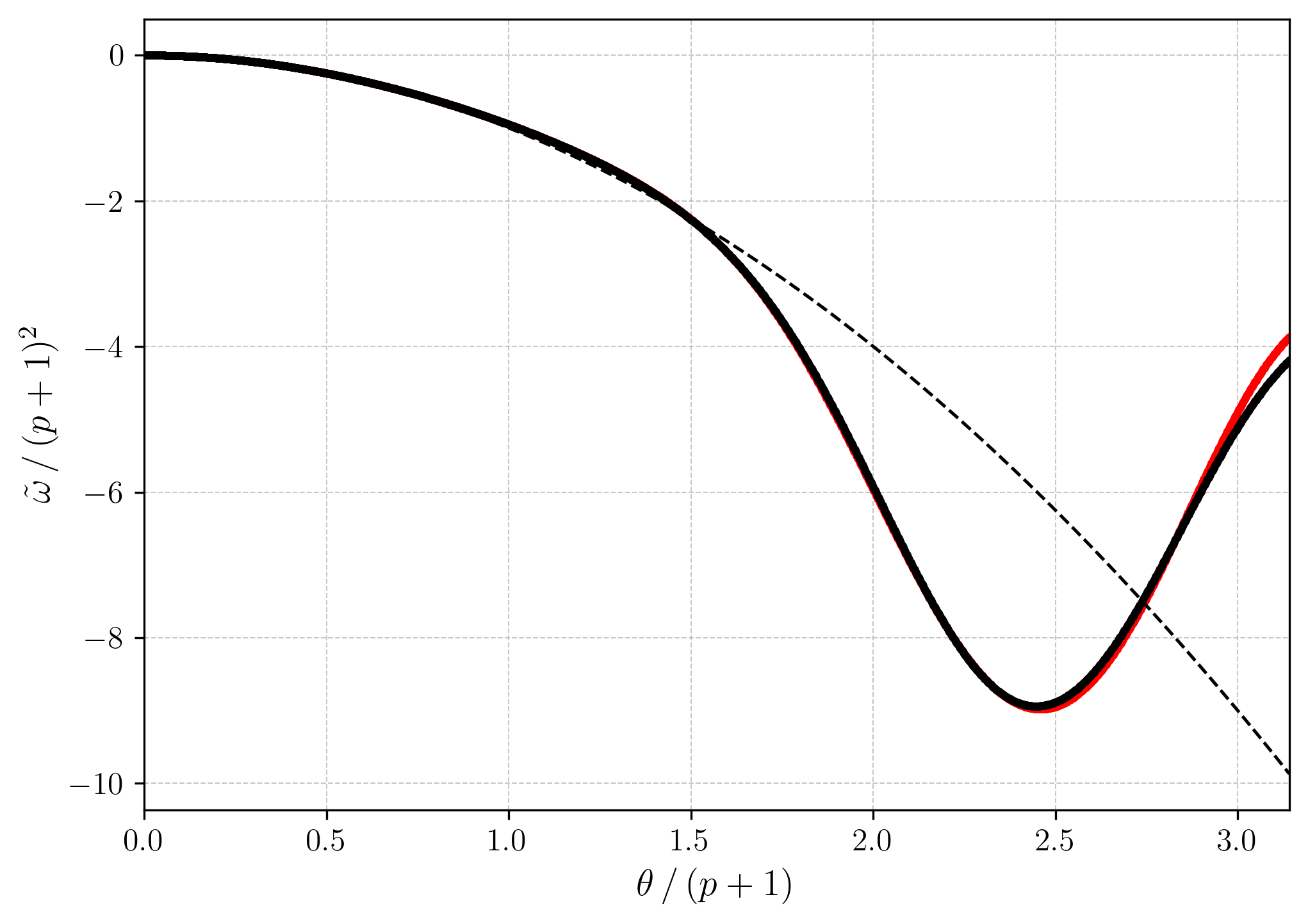}
 \end{subfigure}
\caption{Combined-mode dissipation curves for the SD method for different order of approximation. From top to bottom, left to right: $p=1,2,3,4$. Black curves denotes the compact stencil formulation and red curves indicates the standard approach. Both schemes are fully centered.}
   \label{fig:eigen_diff_single}
\end{figure}
For completeness, in figure~\ref{fig:eigen_diff_single} we report the combined numerical frequency $\widetilde{\omega}$ obtained for the SD scheme employing both the standard and the compact approaches for different polynomial orders. We can notice similar results with respect to the standard eigenanalysis. The standard formulation is constantly characterized by a reduced amount of numerical dissipation, in particular for low-orders of approximation. The compact formulation, instead, is quite similar to the standard approach for extended ranges of wavenumber and it only deviates from it inside medium and high-wavenumber regions. This is particularly interesting as those frequencies are the most sensitive ones in terms of robustness and stability of the scheme. 

\begin{figure}[h!]
 \centering  
 \begin{subfigure}{0.48\textwidth}
     \includegraphics[width=\textwidth]{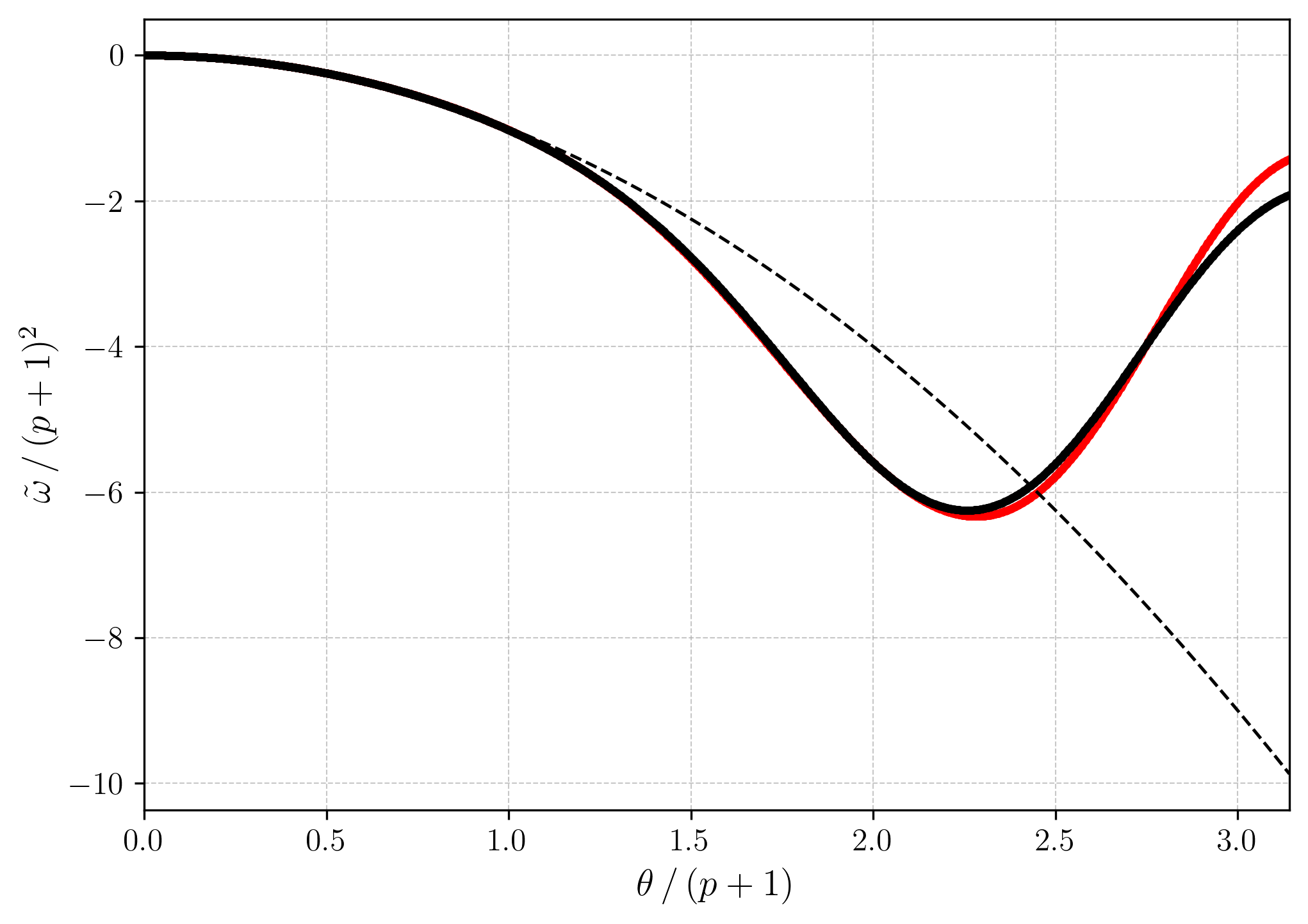}
 \end{subfigure}
 \begin{subfigure}{0.48\textwidth}
     \includegraphics[width=\textwidth]{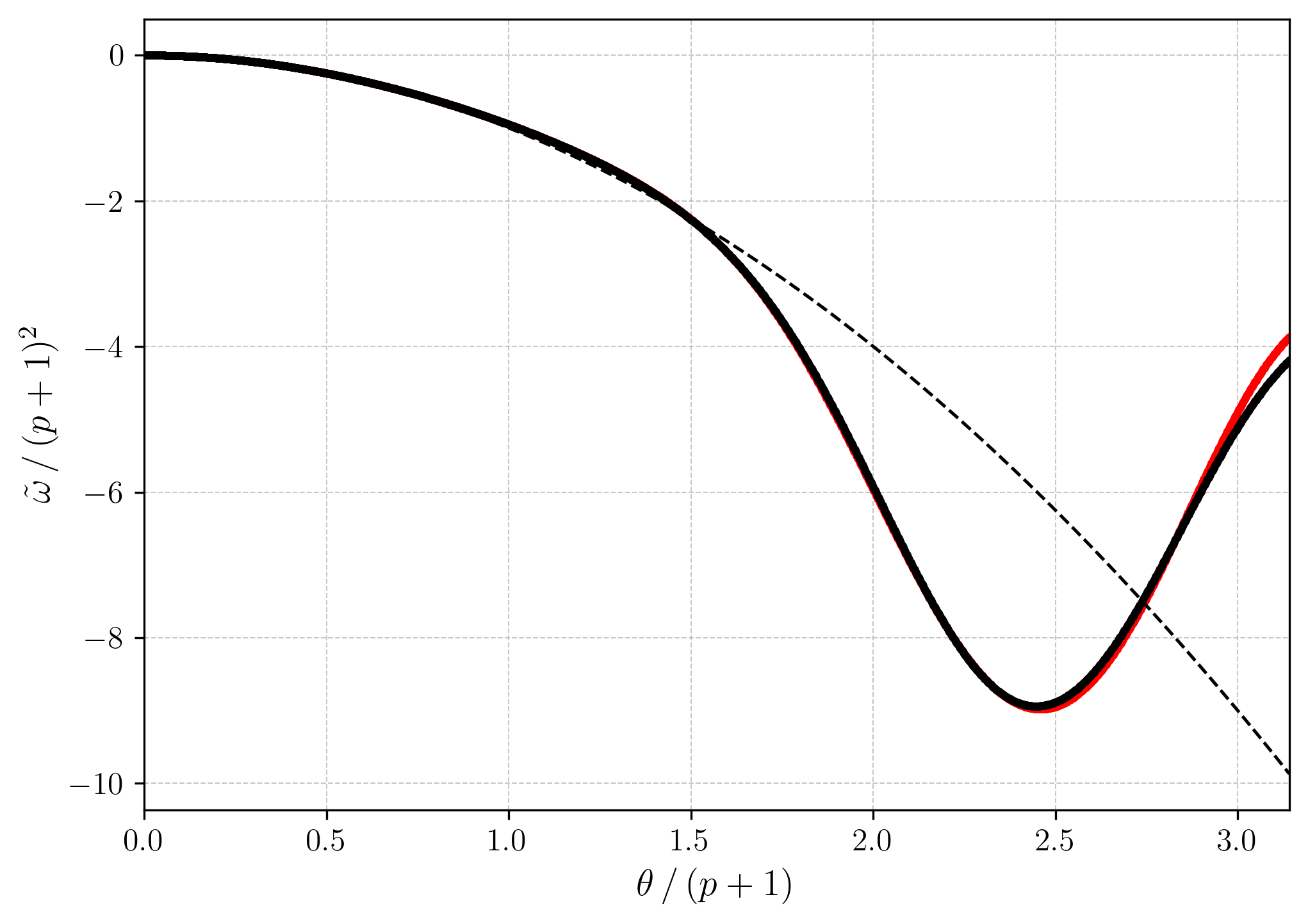}
 \end{subfigure}
 \begin{subfigure}{0.48\textwidth}
     \includegraphics[width=\textwidth]{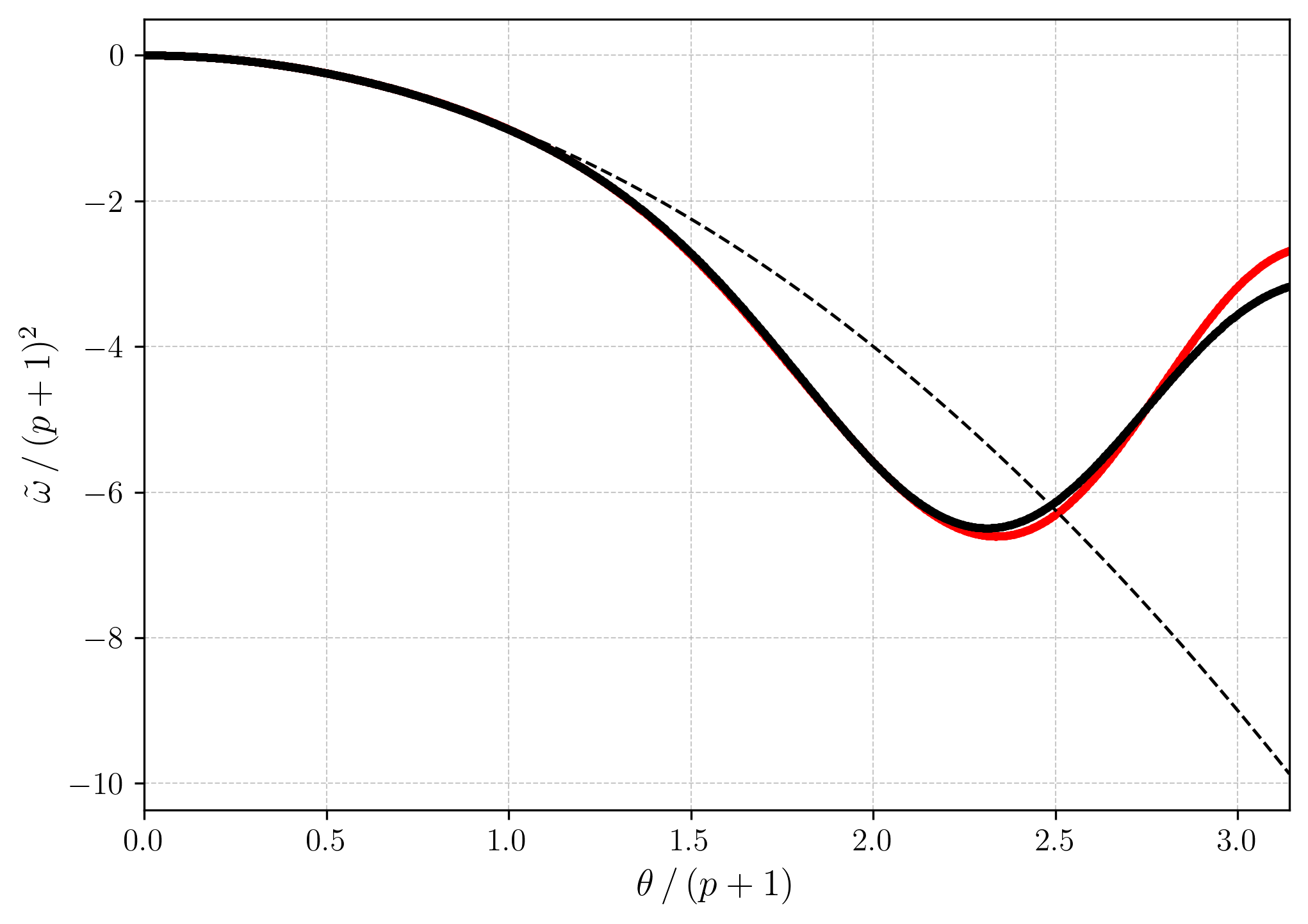}
 \end{subfigure}
 \begin{subfigure}{0.48\textwidth}
     \includegraphics[width=\textwidth]{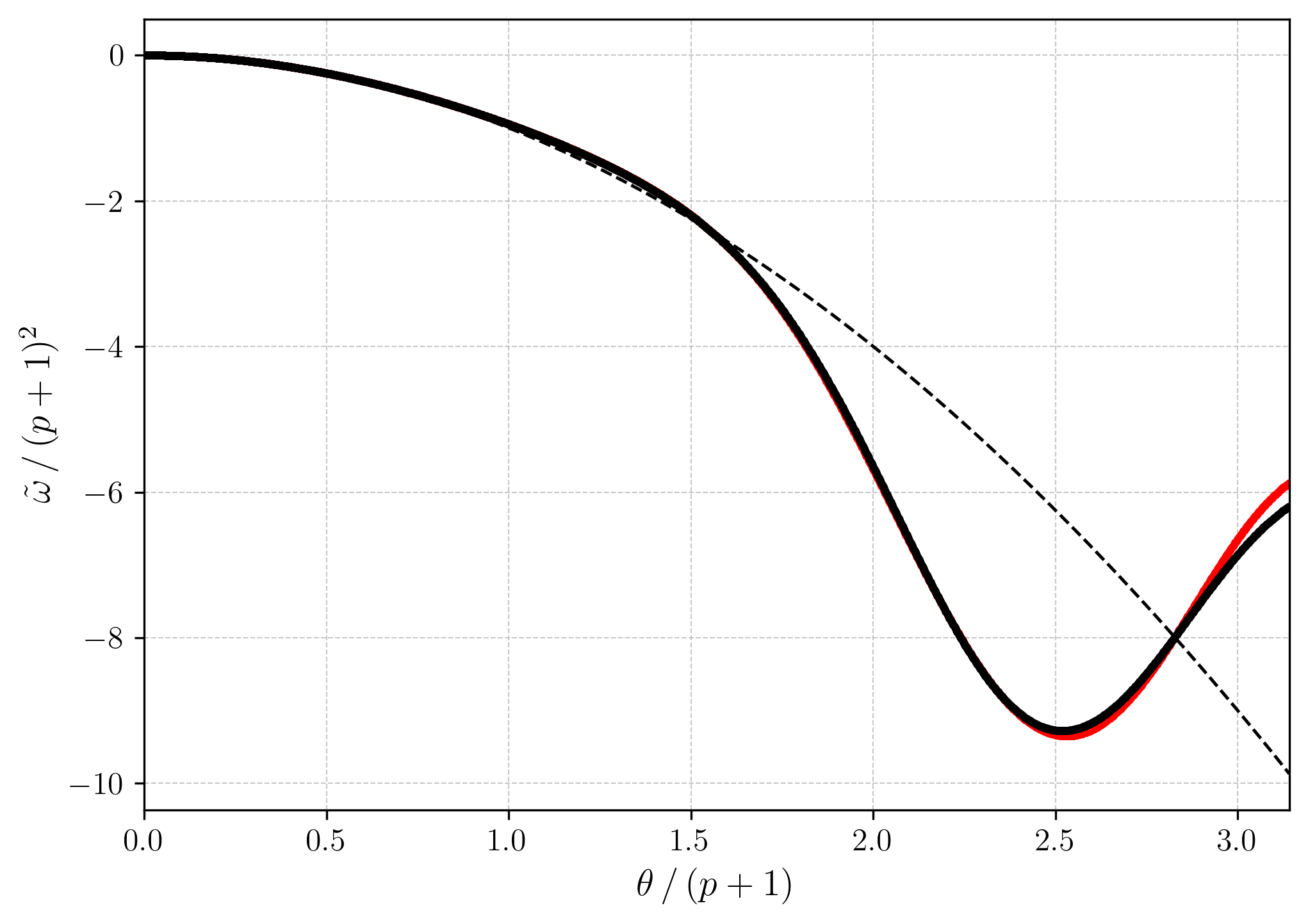}
 \end{subfigure}
  \begin{subfigure}{0.48\textwidth}
     \includegraphics[width=\textwidth]{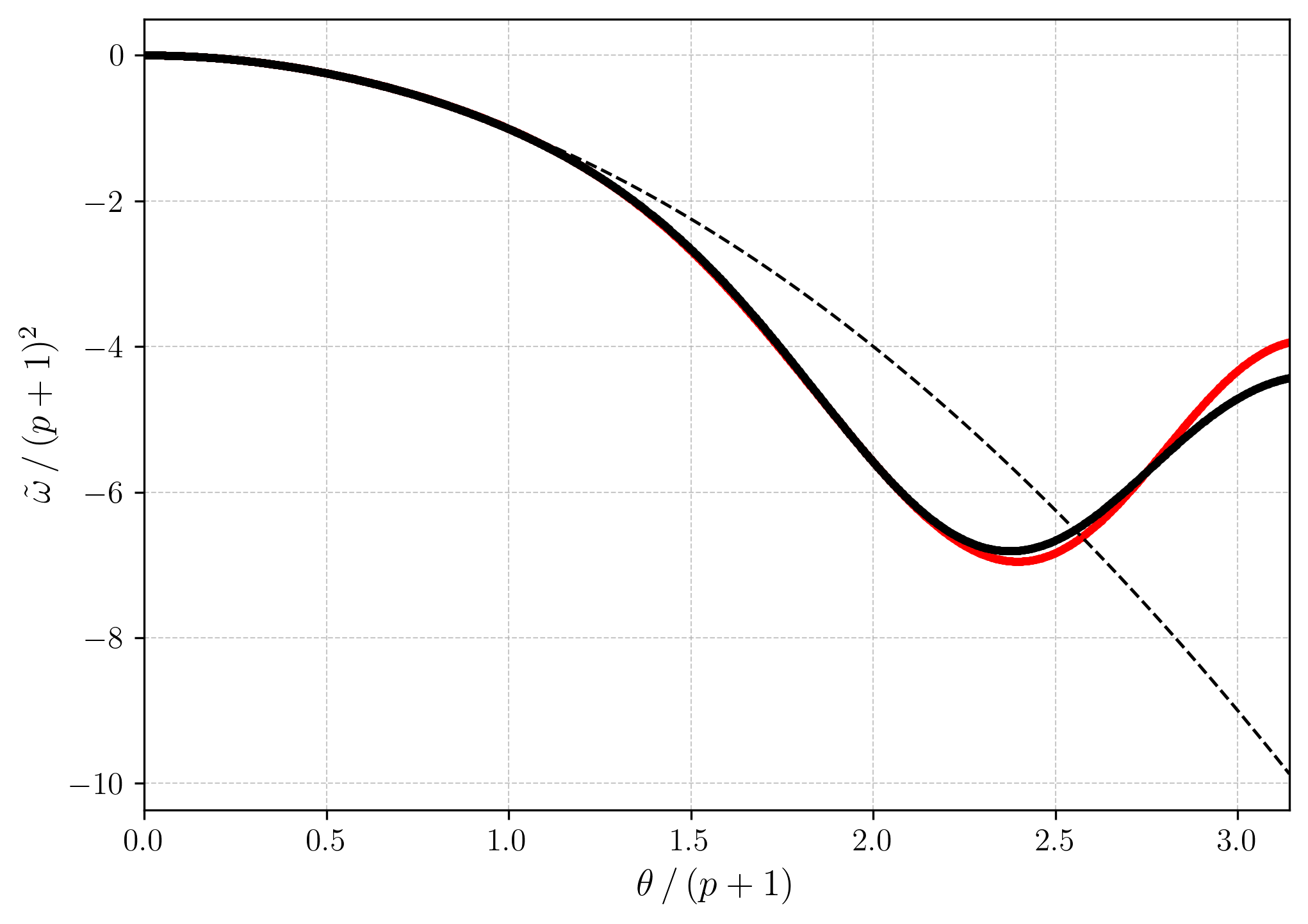}
 \end{subfigure}
 \begin{subfigure}{0.48\textwidth}
     \includegraphics[width=\textwidth]{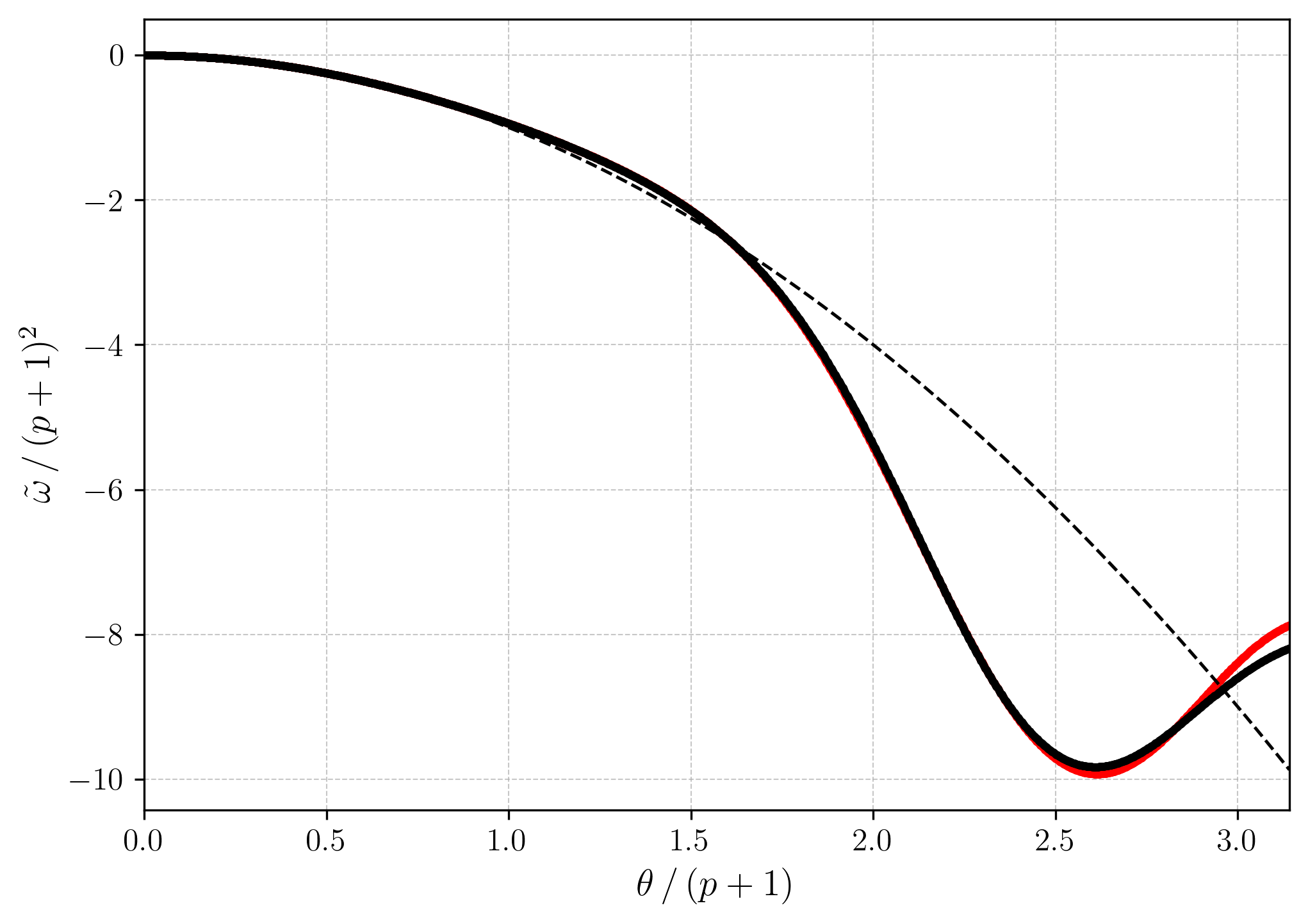}
 \end{subfigure}
   \begin{subfigure}{0.48\textwidth}
     \includegraphics[width=\textwidth]{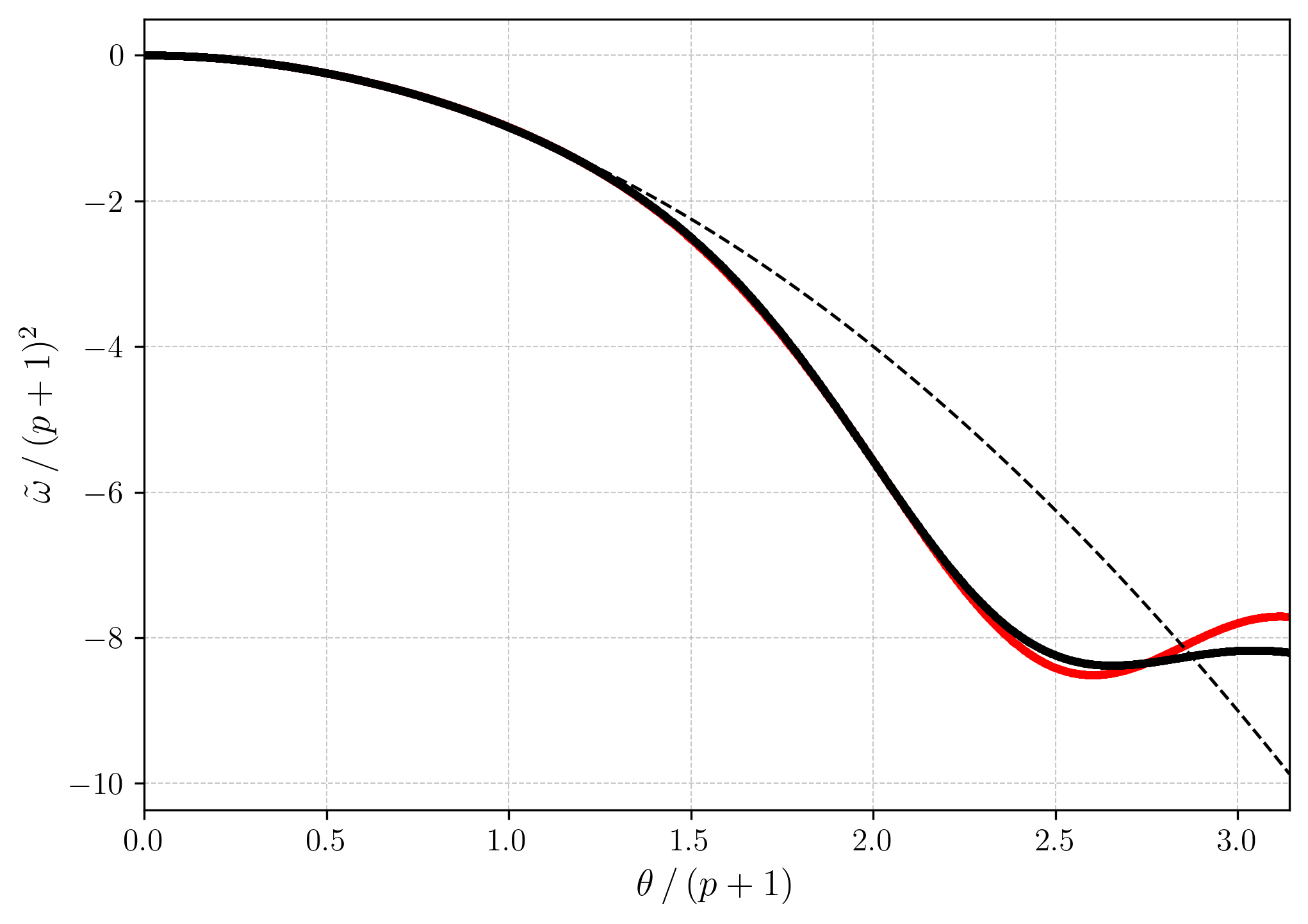}
 \end{subfigure}
 \begin{subfigure}{0.48\textwidth}
     \includegraphics[width=\textwidth]{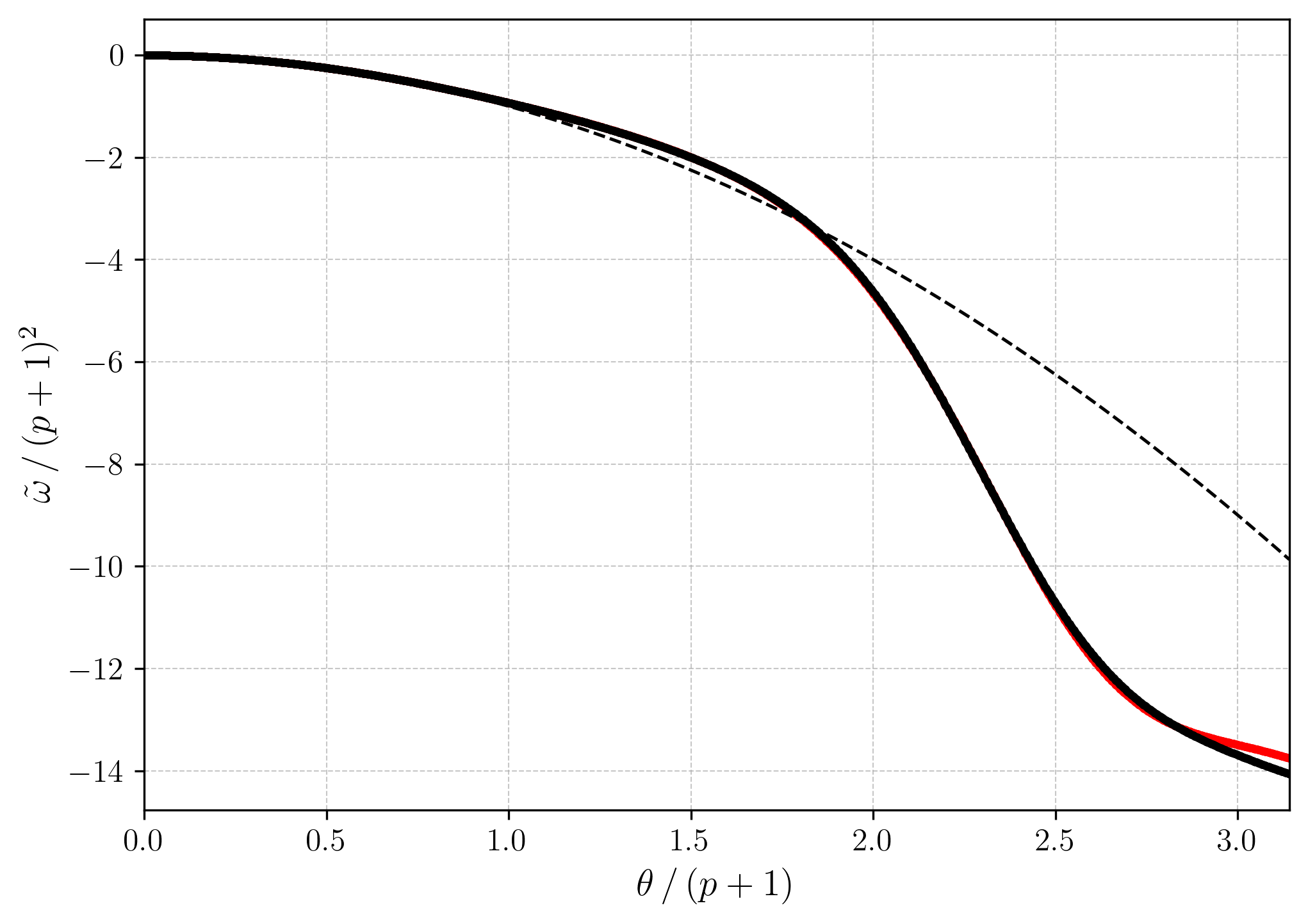}
 \end{subfigure}
\caption{Combined-mode dissipation curves for the SD method for different $p=3$ (left) and $p=4$ (right) and different values of the interior penalty parameter $\eta_{\rm IP}$. From top to bottom $\eta_{\rm IP}=0.0,0.01,0.02,0.05$. Black curves denotes the compact stencil formulation and red curves indicates the standard approach.}
   \label{fig:eigen_diff_tau_single}
\end{figure}

In figure~\ref{fig:eigen_diff_tau_single} similar plots are shown for $p=3$ and $p=4$ with the inclusion of the interior penalty term, where we varied the relevant parameter up to $\eta_{\rm IP}=0.05$ as larger values lead to very negative values of dissipation curves. 
As expected, in this case, we can notice a significant increase in numerical dissipation for high-wavenumbers as the interior penalty parameter increases.
%

\bibliographystyle{ctrnat}
\bibliography{j_abbrv,bibexport}
\end{document}